\documentclass[a4paper,10pt]{article}

\usepackage{amssymb}
\usepackage{amsmath}
\usepackage{amsthm}

\usepackage[margin=2cm]{geometry}
\usepackage[T1]{fontenc}
\usepackage{newtxtext}
\usepackage[varvw]{newtxmath}

\usepackage{graphicx}
\usepackage{subcaption}
\usepackage{booktabs}

\usepackage{authblk}
\usepackage{cite}
\usepackage{hyperref}

\usepackage{placeins}

\usepackage{enumitem}
\setlist[itemize]{leftmargin=*}

\theoremstyle{plain}
\newtheorem{theorem}{Theorem}[section]
\newtheorem{corollary}[theorem]{Corollary}

\newtheorem{proposition}[theorem]{Proposition}

\theoremstyle{definition}

\theoremstyle{remark}
\newtheorem{remark}[theorem]{Remark}

\title{\bf Kinetic modeling for evolutionary opinion games with replicator-driven strategies}

\author[1]{Giacomo Albi}
\author[,1]{Matteo Piu\thanks{Corresponding author.}}
\author[2]{David N Reynolds}
\affil[1]{University of Verona, Department of Computer Science,  Strada Le Grazie 15, 37134 Verona, Italy\\ 
Emails: giacomo.albi@univr.it, matteo.piu@univr.it}
\affil[2]{University of Warsaw, Institute of Applied Mathematics and Mechanics, Banacha 2, 02-097 Warszawa, Poland\\ 
Email: d.reynolds@uw.edu.pl}

\date{}

\begin{document}

\maketitle


\begin{abstract}
\noindent We propose a kinetic model for opinion dynamics in which agents carry a public opinion, a fixed private conviction and a mixed adaptive strategy that sets their propensity to compromise and evolves through a payoff-driven rule of replicator type. In the quasi-invariant regime we derive a Vlasov--Fokker--Planck equation for the joint distribution of opinions, convictions and strategies. For constant payoffs the strategy marginal has an exact projective representation, which allows for a complete analysis of the two-strategy case and of representative three-strategy scenarios, including oscillatory and heteroclinic dynamics. For quadratic state-dependent payoffs we characterize a hyperbolic segregation boundary for the preferred strategy. Formal asymptotics for fast opinions and fast strategies give reduced descriptions, and the latter closure is shown not to reproduce the stationary opinion variance. Numerical simulations support the results.
\end{abstract}



\section{Introduction}

Mathematical models of interacting agents provide a flexible framework for describing collective phenomena emerging from simple individual rules, with applications ranging from biological aggregation and animal swarming to social and economic dynamics, control and optimization; see, among others, \cite{CuckerSmale2007,CarrilloFornasierToscaniVecil2010,BisiLoy2024}. Kinetic and mean-field descriptions are particularly effective in this context, since they connect microscopic interaction mechanisms with the evolution of population distributions and of their macroscopic observables \cite{pareschi2013interacting}. They have been used, in particular, for alignment and consensus mechanisms in flocking and swarming \cite{HaTadmor08,CarrilloFornasierToscaniVecil2010}, for consensus-based optimization \cite{PinnauTotzeckTseMartin2017,CarrilloChoiTotzeckTse2018}, and for the mean-field control of large interacting systems \cite{AlbiAlmiMorandottiSolombrino2022,AlbiCalzolaPiu2027}. Within this framework, opinion formation is one of the most studied examples of collective social dynamics. In many social situations individuals do not simply average the opinions they encounter: they weigh the pressure to conform against the wish to remain consistent with their own convictions, and the expressed opinion may differ from the one privately held \cite{kuran1995private}. Which of the two attitudes prevails is rarely fixed: it depends on what the interaction has rewarded in the past. This motivates models in which the ``degree of compromise'' is not a parameter but an adaptive variable.

Classical approaches describe opinions as discrete or continuous variables updated through repeated averaging or social influence. The DeGroot model \cite{DeGroot1974} provides one of the earliest formulations of consensus through averaging, while the Friedkin--Johnsen model \cite{FriedkinJohnsen1990} introduces persistent individual convictions and heterogeneous susceptibility to social influence. Bounded-confidence mechanisms, such as the Deffuant--Weisbuch and Hegselmann--Krause models \cite{DeffuantEtAl2000,HegselmannKrause2002}, account for selective interactions and may generate clustering and polarization; for broader surveys on social influence and dynamic social networks we refer to \cite{ProskurnikovTempo2017,ProskurnikovTempo2018}. At the mesoscopic level, kinetic models based on binary interactions describe large populations with heterogeneous opinions and random fluctuations \cite{Toscani2006}, and have been extended to account for network effects, data-driven interactions,
and graphon-based models with declustering controls
\cite{AlbiPareschiZanella2017,AlbiCalzolaDimarco2025,during2024breaking}. Kinetic models incorporating heterogeneous convictions and additional individual characteristics have also been investigated in
\cite{BondesanBorsotti2026,BrugnaToscani2015,DuringWolfram2015, PareschiToscaniTosinZanella2019}, while the distinction between publicly expressed opinions and privately held beliefs has been studied in \cite{RoyBiswas2022}. In most of these models, however, the strength of compromise, stubbornness or self-thinking is prescribed through fixed or  state-dependent parameters \cite{VinhNguyen2023DiscreteandContinuousDynamicalSystems,glendinning2025optimalwayliemicroscopic,Nugent_Gomes_Wolfram_2025}. A different modeling question arises when the behavioral attitude that governs the opinion update is itself adaptive and changes according to the outcome of social interactions.

Evolutionary game theory provides a natural framework for selection and adaptation in interacting populations, with classical applications in evolutionary biology \cite{TraulsenClaussenHauert2005,Cooney2020}. In its standard formulation, the relative amount of competing strategies evolves according to their payoff advantage, which leads to the replicator dynamics \cite{TaylorJonker1978,HofbauerSigmund2003,HofbauerSigmund1998}. Its connection with microscopic stochastic interactions and learning processes has been investigated from several viewpoints \cite{TraulsenClaussenHauert2005,MertikopoulosMoustakas2010}. Distributional extensions, in which agents carry heterogeneous mixed strategies rather than sharing a single population strategy, have been developed in kinetic and integro-differential settings \cite{BoccabellaNataliniPareschi2011,PinascoRodriguezCartabiaSaintier2021}, with applications to socio-economic interactions \cite{PinascoRodriguezCartabiaSaintier2018,PinascoSaintierKind2024}.In particular, distributional replicator dynamics may exhibit convergence of the mean strategy to a Nash equilibrium while preserving heterogeneity in the strategy distribution \cite{LahkarSeymour2013}. Closely related structures arise in spatially inhomogeneous evolutionary games, where position and mixed strategy evolve simultaneously \cite{AmbrosioFornasierMorandottiSavare2021,BonafiniFornasierSchmitzer2023}.

The interaction between opinion formation and strategic decision making has received increasing attention. Game-theoretic formulations interpret consensus and disagreement as strategic outcomes and couple social influence with adaptive decisions \cite{GargiuloRamasco2012,AlbiPareschiZanella2019,ZinoYeCao2020}. Dynamical opinion games relate attraction towards other agents and towards individual convictions to Nash equilibria \cite{LearReynoldsShvydkoy2021,ReynoldsTudisco2025}, while coevolutionary models combine the evolution of opinions with decisions or actions on interconnected social layers \cite{GargiuloRamasco2012,ZinoYeCao2020}. Related game-theoretic approaches have also investigated the coexistence of opinions induced by mixed evolutionary strategies \cite{CaoLi2008}. The viewpoint studied here is different and complementary. Each individual carries a continuous public opinion, a fixed private conviction and a \emph{mixed} behavioral strategy; the strategy evolves through a replicator-type rule driven by payoffs that may depend on the opinion--conviction state, and it directly modulates the balance between social compromise and persistence towards the private conviction. Strategy selection and opinion formation are therefore coupled in both directions: the current social state enters the payoff, while the evolving strategic composition determines how an agent reacts to subsequent opinion interactions. To the best of our knowledge, the resulting joint kinetic description of opinions, convictions and mixed strategies, and in particular the interplay between the strategic heterogeneity of the population and the evolution of its opinion distribution, has not been analyzed in this form. The coupling of two processes on different time scales has also been studied for the coevolution of opinions and of the interaction network \cite{HolmeNewman2006,GrossBlasius2008}, and the use of slow-fast expansions in kinetic theory goes back to the Chapman--Enskog procedure \cite{CercignaniIllnerPulvirenti1994,Kuehn2015}. We adapt these approaches here by allowing the slow variable to be either the strategic composition of the population or its opinion distribution, depending on the regime. 

For each agent, we model the public opinion $v$, their private conviction $u$ that does not evolve, and a mixed strategy $\sigma$ over a finite set of behavioral strategies. The strategy determines an effective compromise propensity $\Theta(\sigma)$, which controls the balance between alignment with the interacting agent and relaxation towards the private conviction. The strategy is updated through a payoff-driven binary rule, which is the microscopic counterpart of the replicator dynamics. For a single strategy, or if the strategy is frozen, the model reduces to a kinetic opinion model with compromise and self-thinking terms \cite{Toscani2006} and with attraction towards a private conviction of Friedkin-Johnsen type \cite{FriedkinJohnsen1990}. From these binary interactions we formulate the Boltzmann description and derive, in a quasi-invariant interaction regime \cite{Villani1998,Toscani2006}, a Vlasov--Fokker--Planck equation for the joint distribution of opinions, convictions and strategies. The resulting equation is nonlinear through the mean opinion and the population-averaged payoffs, and the equations for the macroscopic moments are in general not closed, since correlations and higher-order information on the strategy distribution enter the dynamics. The two time scales of the opinion and of the strategy updates play a central role in the analysis.

The analysis is not a direct application of known kinetic theory, for three reasons. First, the strategy variable lives on a simplex where the replicator flow is degenerate at the boundary, so that pure strategies are invariant and the long-time distribution can concentrate on faces and vertices. Second, the equation is nonlinear through the mean opinion and the average payoffs, and the quenched convictions turn it into a family of coupled Fokker-Planck equations indexed by $u$; no finite closure of the moment hierarchy is available in general. Third, two internal time scales, those of opinion and of strategy updates, are present, and their ratio determines which description is appropriate. Our main tool is a representation of the strategy flow that does not rely on moments: for constant payoffs the strategy marginal is transported by a projective map, so that the strategic part of the problem becomes exactly solvable and the opinion dynamics can be studied conditionally on it, at any ratio of the time scales. When this structure is not available, as for state-dependent payoffs, or when we want reduced descriptions for a general number of pure strategies, we resort to the geometry of the payoff gain and to formal asymptotics in the ratio of the time scales.

We begin by introducing the binary model and give explicit conditions under which the opinion and strategy domains are preserved, we derive the Vlasov--Fokker--Planck equation in the quasi-invariant regime, and we identify the structural properties of the replicator part (mass conservation, tangency to the boundary of the simplex, no flux through it) together with the marginal and moment equations.

For constant payoff matrices the strategy marginal admits an exact projective representation \cite{LahkarSeymour2013}: it is the push-forward of the initial distribution through a projective map driven by the cumulative relative gains, which solves an autonomous finite-dimensional ODE. This holds without any monokinetic closure, at arbitrary ratio of the opinion and strategy time scales, and it extends to the joint strategy--conviction marginal. This extension allows us to determine exactly the strategy--conviction correlations governing the mean opinion dynamics. Stationary strategy distributions are not Dirac masses in general and retain information on the initial heterogeneity, impacting the asymptotic opinion distribution.
 
For two strategies, interpreted as compromise and stubbornness, with constant payoff, the exact flow yields a complete classification of the long-time behavior into four regimes (compromise dominance, stubbornness dominance, stable coexistence and bistability) determined by the sign of the payoff gain. It gives closed formulas for the evolution and the limit of the mean opinion, which may retain memory of the initial mean opinion and of transient strategy--conviction correlations, and an explicit description of the asymptotic joint distribution as a mixture of Gibbs densities whose centers depend on the strategy and on the conviction. With a third, moderate strategy we identify dynamical behaviors with no two-strategy analogue (interior coexistence, rock-paper-scissors oscillations and heteroclinic cycles) and their signature on the opinion--conviction distribution.

Considering again the two-strategy case, but for quadratic state-dependent payoffs, where the strategy marginal is no longer autonomous, we characterize the time-dependent segregation boundary, which is a hyperbola in the opinion--conviction plane, and the admissible classes of stationary states. Such states retain coexistence of both stubborn and compromising agents, where the predominant strategy is determined by the parameters of the payoff matrix.

We study the asymptotics under two regimes generated by the time scales. When opinions relax rapidly, the slow strategic dynamics evolve to a conditional Gibbs equilibrium with a self-consistent mean opinion. When strategies adapt rapidly, the selected strategic profile enters the reduced opinion dynamics only through its average compromise propensity. We show that the factorization closure that underlies this reduction reproduces the first moments of the full model but not the stationary opinion variance when the limiting strategy distribution is non-degenerate, and we quantify the difference.
 
We provide numerical support via a Monte Carlo particle scheme to illustrate and test the analytical predictions in all the regimes above.

Beyond the technical results, the analysis shows that strategic heterogeneity can persist in the long-time limit and affects the stationary opinion profiles, which are mixtures over strategies rather than the profile of a single representative agent. Furthermore, when compromise dominates, the mean opinion need not converge to the mean conviction, and it retains memory of the initial mean opinion and of transient strategy--conviction correlations. Moreover, with state-dependent payoffs the stationary states are concentrated on the pure strategies, and the instantaneous payoff advantage separates the opinion--conviction plane into regions favoring compromise or stubbornness; and, finally, cyclic payoff structures induce a periodic or intermittent alternation of the opinion--conviction patterns.

The domain preservation of the binary dynamics, the exact strategy flow and the analysis of the two-strategy model with constant payoffs are established rigorously, up to the convergence of the conditional opinion densities, for which we adapt an existing entropy argument. The derivation of the kinetic equation in the quasi-invariant regime and the reductions in the two asymptotic regimes are formal: they rely on expansions in the ratio of the time scales \cite{Kuehn2015} and, for fast strategies, on explicit additional assumptions, which we state and test numerically. For the mathematical foundation of mean-field evolutionary games driven by replicator-type dynamics in probability measure spaces, we refer to \cite{AmbrosioFornasierMorandottiSavare2021}. The classification of the three-strategy dynamics is limited to representative cases; a complete topological classification of the underlying replicator flows is available in \cite{Bomze1995}. A rigorous justification of the mean-field and quasi-invariant limits, for instance through propagation of chaos \cite{Sznitman1991,JabinWang2017}, and of the formal reductions is not addressed here; the numerical experiments provide supporting evidence.

The paper is organized as follows. In Section~\ref{sec:model} we introduce the binary opinion and strategy dynamics, discuss the admissibility of the microscopic state space and specialize the model to two  strategies. Section~\ref{sec:kinetic} develops the kinetic formulation and the quasi-invariant limit, together with the marginal and macroscopic equations and the exact strategy flow for state-independent payoffs. In Section~\ref{sec:asymptotic} we study the fast-opinion and fast-strategy regimes and derive the corresponding reduced descriptions. Section~\ref{sec:exact-analysis} contains the analysis, beyond scale separation, of the two-strategy model with constant and state-dependent payoffs and of the three-strategy model with constant payoffs. Section~\ref{sec:numerics} is devoted to the numerical experiments, and Section~\ref{sec:conclusions} collects conclusions and perspectives.

\section{Microscopic opinion dynamics with adaptive strategies}\label{sec:model}

We consider binary interactions between two agents characterized by their public opinions $v,v_*\in I:=[-1,1]$, private convictions $u,u_*\in I$, and  behavioral strategies $\sigma,\sigma_*\in\Delta_{K-1}$, where $\Delta_{K-1}\subset \mathbb R^K$ is the standard probability simplex, and $K\in\mathbb{N}\setminus\{0\}$. The values $-1$ and $+1$ represent, respectively, the two extreme positions on a given topic, corresponding to strongly negative and strongly positive opinions or convictions.

We distinguish two  fundamental behavioral attitudes governing the  opinion evolution: \emph{compromise} ($C$), representing the tendency to align with the interacting partner, and \emph{stubbornness} ($S$), reflecting the persistence toward one's own private conviction. Denoting the post-interaction states by $(v',u',\sigma')$ and $(v_*',u_*',\sigma_*')$, we consider the following opinion dynamics:
\begin{equation}\label{eq:binary-opinion}
\begin{cases}
v' = v + \alpha\left(
\Theta(\sigma)(v_*-v)
+
\bigl(1-\Theta(\sigma)\bigr)(u-v)
\right)
+ D(v)\xi,
\\[1mm]
v_*' = v_* + \alpha\left(
\Theta(\sigma_*)(v-v_*)
+
\bigl(1-\Theta(\sigma_*)\bigr)(u_*-v_*)
\right)
+ D(v_*)\xi_*,
\end{cases}
\end{equation}
where $\alpha>0$ is the interaction rate, and $\xi,\xi_*$ are centered random variables with zero mean and variance $\nu^2$, i.e. $\langle\xi\rangle=\langle\xi_*\rangle=0$,  and $\langle\xi^2\rangle=\langle\xi_*^2\rangle=\nu^2$, representing self-thinking fluctuations. Moreover, $\xi$ and $\xi_*$ are assumed to be mutually independent and independent of the agent states $(v,u,\sigma,v_*,u_*,\sigma_*)$. The private convictions act as inherent traits and remain strictly unchanged:
\[
u'=u,
\qquad
u_*'=u_*.
\]
The quantity $\Theta(\sigma)\in[0,1]$ represents the effective \textit{compromise propensity} induced by the behavioral state $\sigma$: values close to $1$ favor $C$, while values close to $0$ favor $S$. More precisely, we assume that  each agent is endowed with a finite set of pure strategies $\Pi=\{\pi_1,\dots,\pi_K\}$, and that its behavioral state is described by a mixed strategy $\sigma=(\sigma_1,\ldots,\sigma_K)$, where  $\sigma_i$ denotes the weight assigned to the pure strategy $\pi_i$. The mixed strategy determines the effective propensity for compromise, i.e.~the effective degree of openness, as
\begin{equation*}
    \Theta(\sigma)
    =
    \sigma^\top\theta
\end{equation*}
where $\theta=(\theta_1,\dots,\theta_K)^\top$ is the vector of compromise
parameters. 

\begin{remark}[Role of the compromise parameters] Note that  the compromise parameters determine how distinct pure strategies are from the viewpoint of opinion dynamics. In particular, if $\theta_i=\widetilde\theta$ for all $i=1,\ldots,K$, then $\Theta(\sigma)=\widetilde\theta$ for every $\sigma\in\Delta_{K-1}$. Hence, all pure strategies induce the same effective compromise level and are indistinguishable with respect to their direct effect on the opinion dynamics. More generally, if $\theta_i=\theta_j=\widetilde\theta$ for some $i\neq j$, then the effective compromise level depends only on the aggregate weight $\sigma_i+\sigma_j$ assigned to these two strategies, although they may remain distinct in terms of their strategic interactions and payoffs. To ensure that the pure strategies are distinct from the viewpoint of opinion dynamics, we therefore assume, after relabeling if necessary,
\begin{equation}\label{theta_order}
1\geq\theta_1>\theta_2>\ldots>\theta_K\geq0.
\end{equation} Under this convention, the strategies are ordered according to their degree of openness to the opinion of the interacting agent, with $\pi_1$ being the most compromise-oriented strategy and $\pi_K$ the most stubborn one.
\end{remark}

The adaptation of the mixed strategy $\sigma_i$ is driven by an interaction payoff. Thus, we introduce the payoff matrix
\[
J=(J_{ij})_{i,j\in\{1,\dots,K\}},
\]
where $J_{ij}=J_{ij}(v,u,v_*,u_*)$ is the \emph{pairwise payoff} of an agent playing $\pi_i$ against an opponent playing $\pi_j$. 

\begin{remark}[Symmetric game] We assume that both interacting agents share the same payoff rule, up to exchange of their states and strategic roles. This corresponds to a symmetric game in the game-theoretic sense, although the payoff matrix $J$ itself need not be symmetric. In particular, non-symmetric payoff matrices, including skew-symmetric ones, are admissible, whereas general bimatrix games are not considered here, and would require a multi-population extension of the model. \end{remark}

Given an agent with mixed strategy $\sigma$ interacting with an agent with mixed strategy $\sigma_*$, we define the \emph{expected payoff of
pure strategy} $\pi_i$, denoted by $P_i$, as
\begin{equation*}
P_i(v,u;v_*,u_*,\sigma_*)
=
\sum_{j=1}^K
\sigma_{*,j}\,
J_{ij}(v,u,v_*,u_*),
\end{equation*}
that is, as the average payoff of $\pi_i$ against the pure strategies
of the opponent, weighted by their frequencies. The corresponding
\emph{expected payoff of the mixed strategy} $\sigma$ is
\[
\overline{P} (v,u,\sigma; v_*,u_*,\sigma_*)
=
\sum_{i=1}^K\sigma_iP_i(v,u;v_*,u_*,\sigma_*).
\]
The strategy update then follows a payoff-driven evolution rule:
\begin{equation}\label{eq:binary-game}
\begin{cases}
\sigma_i'
=
\sigma_i+\beta\,\sigma_i\bigl(P_i-\overline P\bigr),
& i=1,\dots,K,
\\[1mm]
\sigma_{*,j}'
=
\sigma_{*,j}+\beta\,\sigma_{*,j}\bigl(P_{*,j}-\overline P_*\bigr),
& j=1,\dots,K,
\end{cases}
\end{equation}
where $\beta>0$ governs the adaptation rate. Eq.~\eqref{eq:binary-game} can be
seen as a binary discrete-time counterpart of the
standard replicator dynamics in which strategies whose payoff exceeds the current average payoff
increase their relative weight, while strategies with below-average
payoff are reduced, e.g. see~\cite{TaylorJonker1978,HofbauerSigmund1998,Weibull1995,Sandholm2010}.

\begin{remark}[Myopic vs.~anticipatory payoffs]
Note that the payoff may account for different behavioral costs associated with
the current interaction state. A typical \emph{myopic} formulation is:
\begin{equation}\label{eq:payoffs-micro}
    J_{ij}(v,u,v_*,u_*)
    =
    -a_{ij}(v-u)^2
    -b_{ij}(v-v_*)^2
    -c_{ij},
\end{equation}
where $a_{ij},b_{ij},c_{ij}\geq0$. The first term penalizes the discrepancy between public opinion and private conviction, the second penalizes disagreement with the interacting agent, while $c_{ij}$ represents an intrinsic cost associated with the pair of strategies. Hence, different strategies may assign different weights to internal coherence and social alignment. Because \eqref{eq:payoffs-micro} is evaluated before the opinion update, it drives a reactive adaptation, and it does not account for the outcome of the subsequent opinion update. Alternatively, one may consider an \textit{anticipatory} payoff in which each strategy is evaluated through the post-interaction opinion it is expected to generate; see, e.g., \cite{Selten1991,CamererHoChong2002}. This introduces additional modeling and analytical difficulties, since the payoff becomes coupled to the interaction rule determining the post-interaction state. We leave this to a future study.
\end{remark}

Since
\begin{equation}\label{sigma_constr}
\sigma_K=1-\sum_{i=1}^{K-1}\sigma_i,
\end{equation}
taking strategy $\pi_K$ as a reference and denoting by
\(
\widehat{\sigma}
=
(\sigma_1,\dots,\sigma_{K-1})
\)
the independent coordinates in $\Delta_{K-1}$, it is convenient to reformulate
\eqref{eq:binary-game} in terms of the \emph{relative individual  payoff gains}
\[
G_i(v,u;v_*,u_*,\widehat{\sigma}_*)
:=
P_i(v,u;v_*,u_*,\sigma_*)-P_K(v,u;v_*,u_*,\sigma_*),
\qquad i=1,\dots,K-1.
\]
Equivalently, defining $a_{ij}=a_{ij}(v,u,v_*,u_*)$, where
\[
a_{ij}:=J_{ij}-J_{Kj},
\qquad
i=1,\dots,K-1,\quad j=1,\dots,K,
\]
we can write the gains as
\begin{equation}\label{def:Gi}
G_i(v,u;v_*,u_*,\widehat{\sigma}_*)
=
a_{iK}
+
\sum_{j=1}^{K-1}
\sigma_{*,j}
\left(a_{ij}-a_{iK}\right),
\qquad i=1,\dots,K-1.
\end{equation}
The resulting strategy dynamics can then be written as
\begin{equation}\label{eq:binary-game-gain}
\begin{cases}
\sigma_i'
=
\sigma_i+\beta\,\sigma_i
\left(
G_i(v,u;v_*,u_*,\widehat{\sigma}_*)
-
\displaystyle\sum_{\ell=1}^{K-1}
\sigma_\ell G_\ell(v,u;v_*,u_*,\widehat{\sigma}_*)
\right),
& i=1,\dots,K-1,
\\[3mm]
\sigma_{*,j}'
=
\sigma_{*,j}+\beta\,\sigma_{*,j}
\left(
G_{j,*}(v_*,u_*;v,u,\widehat{\sigma})
-
\displaystyle\sum_{\ell=1}^{K-1}
\sigma_{*,\ell}G_{\ell,*}(v_*,u_*;v,u,\widehat{\sigma})
\right),
& j=1,\dots,K-1,
\end{cases}
\end{equation}
where $G_{j,*}$ denotes the relative payoff gain of the interacting agent,
defined analogously in terms of its payoff matrix
$J^*=(J^*_{ij})$ and
\(
a^*_{ij}:=J^*_{ij}-J^*_{Kj}.
\)
The remaining components $\sigma_K$ and $\sigma_{*,K}$ are then
determined by the normalization constraints \eqref{sigma_constr}.

\subsection{Admissibility of the microscopic dynamics}
To guarantee mathematical well-posedness, the opinions $v, v_*$ must remain confined within $I$, and the mixed strategies $\sigma, \sigma_*$ must be preserved within $\Delta_{K-1}$. This is ensured by the following proposition.

\begin{proposition}[Domain preservation]\label{prop:domain_invariance}
Let $v,v_*,u,u_*\in I$ and $\sigma,\sigma_*\in\Delta_{K-1}$. Assume that the diffusion coefficient $D:I\to\mathbb{R}$ satisfies conditions
\begin{enumerate}
    \item[$\mathrm{(D1)}$] $D(\pm1)=0$
    \item[$\mathrm{(D2)}$] $C_D
    :=
    \sup_{|v|<1}
    \frac{|D(v)|}{1-|v|}
    <+\infty$, for some $C_D>0$.
\end{enumerate}
Assume that the random variables $\xi,\xi_*$ are almost surely bounded by $\max\{|\xi|,|\xi_*|\}\leq\xi_{\max}$, and  that $J_{ij}\in[J_{\min},J_{\max}]$ for all $i,j$ and all admissible states.  If
\begin{enumerate}
    \item[(i)] $\alpha \in (0,1]$ and $\xi_{\max} \le \frac{1-\alpha}{C_D}$,
    \item[(ii)] $J_{\max}>J_{\min}$ and $\beta \in \left(0, \frac{1}{J_{\max} - J_{\min}}\right]$,
\end{enumerate}
then $v', v_*' \in I$ and $\sigma', \sigma_*' \in \Delta_{K-1}$ almost surely. 
If $J_{\max}=J_{\min}$ condition \emph{(ii)} is vacuous, and any $\beta > 0$ is admissible.
\end{proposition}
\begin{proof}
By symmetry, it suffices to consider $(v', \sigma')$. Rewriting the opinion update as
\[
v' = (1-\alpha)v + \alpha\left(\Theta(\sigma)v_* + (1-\Theta(\sigma))u\right) + D(v)\xi,
\]
and applying the triangle inequality with
$\Theta(\sigma)\in[0,1]$ and $|v_*|,|u|\leq1$, we obtain
\[
|v'|
\leq
(1-\alpha)|v|+\alpha+|D(v)|\,|\xi|.
\]
For $|v|=1$,  by condition (D1), $|v'|\leq1$.
For $|v|<1$, the invariant bound $|v'|\leq1$  requires
\[
|D(v)|\,|\xi|
\leq
(1-\alpha)(1-|v|).
\]
Hence, a uniform sufficient condition is
\[
|\xi|
\leq
(1-\alpha)
\frac{1-|v|}{|D(v)|},
\]
whenever $D(v)\neq0$. Using (D2), it is therefore
sufficient to require
\[
\xi_{\max}
\leq
(1-\alpha)
\inf_{\substack{|v|<1\\D(v)\neq0}}
\frac{1-|v|}{|D(v)|}
=
\frac{1-\alpha}{C_D},
\] which holds by  (i). Thus $|v'|\leq1$ for every $v\in I$ almost surely.

For the strategic variables, summing $\sigma_i'$  yields:
\[
\sum_{i=1}^K \sigma_i' = \sum_{i=1}^K \sigma_i + \beta \left( \sum_{i=1}^K \sigma_i P_i - \overline{P} \sum_{i=1}^K \sigma_i \right)=1+\beta(\overline{P}-\overline{P})=1,
\] therefore the total mass is conserved. For non-negativity, rewriting the strategy update as:
\[
\sigma_i' = \sigma_i \left( 1 + \beta(P_i - \bar  P) \right).
\]
Since $\sigma_* \in \Delta_{K-1}$, each $P_i$ is a convex combination of payoff entries, so $P_i \in [J_{\min}, J_{\max}]$. The same holds for $\overline{P}$, giving $P_i - \overline{P} \ge J_{\min} - J_{\max}$, hence 
\[
1 + \beta(P_i - \overline{P}) \geq 1 - \beta(J_{\max} - J_{\min})\geq 0.
\]
which is guaranteed by (ii). Thus, $\sigma_i' \ge 0$, and combined with the conservation of mass this gives $\sigma' \in \Delta_{K-1}$.
\end{proof}
\begin{remark}[Examples of diffusion coefficients]
The condition (D2) requires the diffusion coefficient to vanish at least linearly as the opinion approaches the boundary of $I$. An admissible choice corresponds to $D(v)=1-v^2,$ for which $C_D=2$. More generally, coefficients of the form
\[
    D(v)=(1-v^2)^\gamma,
    \qquad \gamma\geq1,
\]
satisfy both conditions (D1) and (D2) and therefore guarantee microscopic
invariance under a suitable uniform bound on the noise. On the other hand, choices with a slower decay at the boundary, i.e. $0\leq\gamma<1$ do not satisfy condition (D2).
\end{remark}

\subsection{The two-strategy binary interaction game}

The general $K$-strategy formulation provides a flexible framework, but the essential mechanisms of behavioral adaptation and opinion coupling already appear in the simplest non-trivial case ($K=2$). 

By setting the compromise parameters $\theta=(1,0)^\top$, it is natural to identify $\pi_C:=\pi_1$ as the \emph{compromise strategy} and $\pi_S:=\pi_2$ as the \emph{stubbornness strategy}, with associated mixed strategy $\sigma\in\Delta_1$ written as $\sigma = (s,1-s)$. 

Given a general $2\times2$ pairwise payoff matrix
\[
J=
\begin{pmatrix}
J_{CC} & J_{CS}\\
J_{SC} & J_{SS}
\end{pmatrix},
\]
and choosing $\pi_S$ as the reference strategy, the relative payoff gain of $\pi_C$ over $\pi_S$, as defined in \eqref{def:Gi}, is given by
\[
G(s_*)=P_C-P_S
=
s_*(J_{CC}-J_{SC})
+
(1-s_*)(J_{CS}-J_{SS}).
\]
The analogous gain $G_*(s)$ for the interacting agent is computed using $J^*_{ij}$, which denotes the payoff entries evaluated from the perspective of the second agent.

With this specification, the effective compromise level, in the corresponding binary opinion
interaction \eqref{eq:binary-opinion} reduces simply to $ \Theta(\sigma) = s$, while the general strategy update
\eqref{eq:binary-game-gain} reduces to the coupled discrete-time replicator equations:
\begin{equation}
\label{eq:two-binary-game}
\begin{cases}
s' = s+\beta\,s(1-s)G(s_*),\\[1mm]
s_*' = s_*+\beta\,s_*(1-s_*)G_*(s).
\end{cases}
\end{equation}

In general, for state-dependent payoffs, $G$ and $G_*$ need not coincide, whereas in the state-independent case one can consider $J^*\equiv J$ and therefore $G_*\equiv G$.   The two-strategy dynamics is thus fully described by the state pair $(s,s_*)$, and the direction of each strategy update is uniquely determined by the sign of the respective relative gain, as discussed in the following remark.

\begin{remark}[Qualitative behavior of two-strategy dynamics]
The qualitative behavior of the strategy dynamics
\eqref{eq:two-binary-game} is governed by the relative gains
$G(s_*)$ and $G(s)$. Since each gain depends
affinely on the opponent's single independent strategy variable, the
resulting dynamics can be characterized by the sign of the gain and,
when it changes sign, by the corresponding threshold at which the gain
vanishes. Indeed, for
$s,s_*\in(0,1)$, the factors $s(1-s)$ and $s_*(1-s_*)$ are positive, so that the directions of the two strategy updates are determined by $G$ and $G_*$, respectively. For the first agent, three distinct regimes can be distinguished:
\begin{itemize}
\item  \emph{Compromise favored.} If $J_{CC} > J_{SC}$, and $J_{CS} > J_{SS},$ then
$G(s_*)>0$ for all $s_*\in(0,1)$. Hence $s$ increases and the compromise strategy becomes dominant.
\item \emph{Stubbornness favored.}
If $J_{CC} < J_{SC}$, $J_{CS} < J_{SS}$, then
$G(s_*)<0$ for all $s_*\in(0,1)$.
Hence $s$ decreases and the stubborn strategy becomes dominant.
\item \emph{Frequency-dependent threshold.} If $(J_{CC}-J_{SC})(J_{CS}-J_{SS})<0$  the strategies have opposing advantages, and the sign of the gain function depends on the opponent's strategy composition. In particular, there exists a  unique threshold
$\widehat s_*\in(0,1)$ given by
\begin{equation*}
\widehat s_* =
\frac{J_{SS}-J_{CS}}
{(J_{CC}-J_{SC})+(J_{SS}-J_{CS})},
\end{equation*}
such that  $G(\widehat s_*)=0$. The response of the agent depends on the sign of $G'(s_*)=(J_{CC}-J_{SC})+(J_{SS}-J_{CS})$, and we have two distinct regimes. If $G'(s_*)>0$, matching the opponent's strategy is highly rewarded ({\em coordination response}). Then $G(s_*)>0$ if $s_*>\widehat s_*$, and $G(s_*)<0$ for $s_*<\widehat s_*$.
The response pushes the system away from $\widehat s_*$ and toward the boundaries.  If  $G'(s_*)<0$, it is advantageous to adopt the minority strategy ({\em anti-coordination response}). Then $G(s_*)<0$ if $s_*>\widehat s_*$, and  $G(s_*)>0$ if $s_*<\widehat s_*$. Here an agent tends to reinforce the strategy that is less represented by the opponent.  
\end{itemize}
Note that the same classification applies to the second agent through the relative gain $G_*(s)$, and for state-dependent payoffs, the two interacting agents may experience
different gain functions, so that $G_*\neq G$. In contrast, when $G\equiv G_*$, namely the two agents share the same payoff structure, the resulting strategy dynamics is symmetric and can be represented in the $(s,s_*)$ plane through the following displacement field
\[
(\Delta s,\Delta s_*)
    =
    \left(
        s(1-s)G(s_*),
        s_*(1-s_*)G(s)
    \right).
\]
Figure~\ref{fig:strategy-vector-fields} illustrates the corresponding qualitative regimes. If $G>0$ on $(0,1)$, both agents are driven toward the compromise state $(1,1)$, whereas if $G<0$ they are driven toward the stubborn state $(0,0)$. In the frequency-dependent case, the nullclines $s=\widehat s$ and $s_*=\widehat s$ divide the square into
regions with opposite strategic responses. For $G(s)=2s-1$, corresponding to $G'>0$, the interaction is of coordination type: the interior state $(1/2,1/2)$ is a saddle and the
symmetric pure states $(0,0)$ and $(1,1)$ are stable equilibria. Conversely, for $G(s)=1-2s$, corresponding to $G'<0$, the interaction is of anti-coordination type: the asymmetric states $(0,1)$ and $(1,0)$ are stable equilibria, while $(1/2,1/2)$ remains a saddle in the full two-agent phase space. Restricted to the invariant diagonal $s=s_*$, however,
the interior threshold is unstable in the coordination case and stable in the anti-coordination case. 

The corresponding population-level classification for constant payoffs, including the asymptotic strategy distribution and its coupling with opinion dynamics, is developed in Section~\ref{sec:2s-constant}.

\begin{figure}[ht]
\centering
\includegraphics[width=\textwidth]{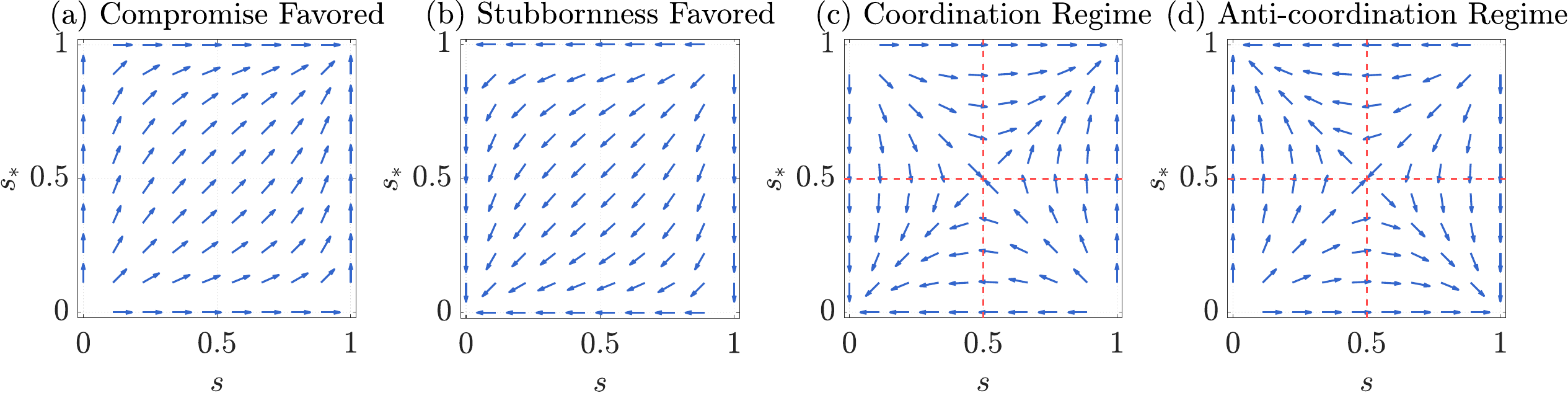}
\caption{Vector-field representation of the two-strategy update in the $(s,s_*)$ plane for symmetric payoffs, $G\equiv G_*$. The arrows show the direction of the discrete displacement induced by the strategy dynamics in the three payoff regimes.}
\label{fig:strategy-vector-fields}
\end{figure}
\end{remark}

\section{Kinetic formulation for opinions with evolutionary strategies}\label{sec:kinetic}

To study the collective dynamics of the agent-based model \eqref{eq:binary-opinion}-\eqref{eq:binary-game}, we move to a kinetic description. We denote by $f(v,\sigma,u,t)$ the continuous joint probability density, which describes the distribution of agents with opinion $v \in I$, behavioral mixed strategy $\sigma \in \Delta_{K-1}$, and conviction $u \in I$ at time $t \ge 0$. 

We assume the density is normalized, such that \begin{equation*}
\int_{\Omega} f(v,\sigma,u,t)\, \mathrm{d}\omega =1
\end{equation*} where $ \mathrm{d}\omega := \mathrm dv\,\mathrm{d}\sigma\,\mathrm du$ denotes the volume element on the state space $\Omega := I\times \Delta_{K-1}\times I$. Here, $\mathrm{d}\sigma$ denotes the $(K-1)$-dimensional Lebesgue 
measure on $\Delta_{K-1}$.

The evolution of $f$ is governed by a Boltzmann-type equation \cite{pareschi2013interacting,Toscani2006} of the  form
\begin{equation}
\label{eq:Boltzmann_full}
    \partial_t f(v,\sigma,u,t)
    =
    \frac{1}{\tau_v}Q^v(f,f)(v,\sigma,u,t)
    +
    \frac{1}{\tau_\sigma}Q^\sigma(f,f)(v,\sigma,u,t),
\end{equation}
where we have introduced two parameters $\tau_v,\tau_\sigma>0$ which dictate the characteristic time scales of opinion formation and strategy adaptation, respectively. The interaction operator $Q^v$ accounts for the variation of the opinion variable induced by binary interactions \eqref{eq:binary-opinion}, while $Q^\sigma$ describes the payoff-driven evolution of the mixed strategy \eqref{eq:binary-game}. The conviction variable $u$ is not modified by the interactions and therefore acts as an intrinsic heterogeneity parameter.

For any smooth test function $\varphi=\varphi(v,\sigma,u)$, the opinion operator is defined in weak form by
\begin{equation*}
\int_\Omega Q^v(f,f)(v,\sigma,u,t)\,\varphi(v,\sigma,u)\,\mathrm{d}\omega =
\iint_{\Omega\times \Omega}
\left\langle
\varphi(v',\sigma,u)-\varphi(v,\sigma,u)
\right\rangle
f(v,\sigma,u,t)f(v_*,\sigma_*,u_*,t)
\,\mathrm{d}\omega\,\mathrm{d}\omega_*,
\end{equation*}
where $v'$ is the post-interaction opinion \eqref{eq:binary-opinion} and the expectation $\langle\cdot\rangle:=\mathbb{E}_\xi[\cdot]$  is taken with respect to the random noise $\xi$, which is independent of the agent states. Analogously, the strategy operator is defined by
\begin{equation*}
\int_\Omega Q^\sigma(f,f)(v,\sigma,u,t)\,\varphi(v,\sigma,u)\,\mathrm{d}\omega  =
\iint_{\Omega\times\Omega}
\left[
\varphi(v,\sigma',u)-\varphi(v,\sigma,u)
\right]
f(v,\sigma,u,t)f(v_*,\sigma_*,u_*,t)
\,\mathrm{d}\omega\,\mathrm{d}\omega_*,
\end{equation*}
where $\sigma'$ is determined by the replicator update \eqref{eq:binary-game}. 

The parameters $\tau_v$ and $\tau_\sigma$ determine the relative speed of the two macroscopic mechanisms:  $\tau_\sigma\ll\tau_v$ implies strategies adapt faster than opinions,  $\tau_v\ll\tau_\sigma$ implies opinions relax faster than strategies, and $\tau_v\sim\tau_\sigma$ corresponds to comparable evolutionary time scales.

\subsection{Quasi-invariant limit and Vlasov--Fokker--Planck model}

The Boltzmann equation \eqref{eq:Boltzmann_full} characterizes the mesoscopic evolution driven by discrete binary jumps. To obtain a tractable continuous formulation for studying the macroscopic and asymptotic behavior, we consider a quasi-invariant (or grazing) interaction regime \cite{pareschi2013interacting,Toscani2006}.

In this regime, each binary interaction produces only a small variation of the state variables, while the interaction frequency increases proportionally.

Therefore, we introduce a small parameter $0<\varepsilon\ll1$ and rescale the interaction parameters in \eqref{eq:binary-opinion} and \eqref{eq:binary-game} as follows
\begin{equation}\label{eq:grazing}
    \alpha \to \varepsilon,
    \qquad
    \beta \to \varepsilon,
    \qquad
    \nu^2 \to \varepsilon \nu^2,
\end{equation}
where the  self-thinking fluctuation is driven by a scaled centered random variable $\xi_\varepsilon=\sqrt\varepsilon\xi$ with variance $\langle \xi_\varepsilon^2\rangle=\varepsilon\nu^2$ and
finite third moment, i.e. $\langle|\xi_\varepsilon|^3\rangle<+\infty$.
Time is appropriately scaled as $t\to t/\varepsilon$ to observe finite macroscopic effects. The time-scaled weak formulation for the scaled density $f_\varepsilon=f_\varepsilon(v,\sigma,u,t)$ reads:
\begin{equation*}
\frac{d}{dt}\int_\Omega \varphi f_\varepsilon\,\mathrm{d}\omega
=
\frac{1}{\varepsilon\tau_v}
\int_\Omega Q^v(f_\varepsilon,f_\varepsilon)\varphi\,\mathrm{d}\omega+
\frac{1}{\varepsilon\tau_\sigma}
\int_\Omega Q^\sigma(f_\varepsilon,f_\varepsilon)\varphi\,\mathrm{d}\omega,
\end{equation*}
for any test function $\varphi(v,\sigma,u)\in C^3_b(\Omega)$ with bounded derivatives up to order three. We analyze the two operators separately.

\paragraph{Grazing limit for the opinion operator.}
Under the scaling \eqref{eq:grazing}, the  opinion dynamics \eqref{eq:binary-opinion} yields:
\begin{equation*}
    v'-v
    =
    \varepsilon A^v(v,\sigma,u;v_*,\sigma_*,u_*)
    +
    \,D(v)\xi_\varepsilon,
\end{equation*} where the interaction drift is given by 
\(    A^v
    =
    \Theta(\sigma)(v_* - v)
    +
    (1-\Theta(\sigma))(u-v).
\) Expanding $\varphi(v',\sigma,u)$ around $(v,\sigma,u)$ gives
\[
\begin{aligned}
\varphi(v',\sigma,u)-\varphi(v,\sigma,u)
&=
\partial_v\varphi(v,\sigma,u)(v'-v)
+
\frac12\partial^2_{vv}\varphi(v,\sigma,u)(v'-v)^2
+
R_\varepsilon^v.
\end{aligned}
\]
Taking 
expectation with respect to $\xi_\varepsilon$, 
we obtain
\begin{equation*}
\left\langle
\varphi(v',\sigma,u)-\varphi(v,\sigma,u)
\right\rangle
=
\varepsilon A^v\,\partial_v\varphi
+
\frac{\varepsilon\nu^2}{2}D(v)^2\partial^2_{vv}\varphi
+
R_\varepsilon^v.
\end{equation*}
Since $\varphi\in C_b^3$, and $\xi_\varepsilon \sim \mathcal{O}(\sqrt{\varepsilon})$ the remainder $R^v_\varepsilon$ can be written as 
\[
\left|\left\langle R_\varepsilon^v\right\rangle\right|\leq \frac{1}{6}|\partial_v^3\varphi(\widehat v,\sigma,u)||v'-v|^3
\]
for some $\widehat v$ between $v$ and $v'$, which gives the bound 
\[
\left|\left\langle R_\varepsilon^v\right\rangle\right|
\leq
\frac{1}{6}\left\|\partial_{v}^3\varphi\right\|_\infty
\left\langle|v'-v|^3\right\rangle
\leq
C_\varphi \left\langle |v'-v|^3\right\rangle
=
\mathcal{O}(\varepsilon^{3/2}).
\]
Therefore, after division by $\varepsilon$, the limiting weak contribution of the opinion operator is given by
\[
\begin{aligned}
&\frac{1}{\varepsilon}
\int_\Omega Q^v(f_\varepsilon,f_\varepsilon)(v,\sigma,u,t)
\varphi(v,\sigma,u)\,\mathrm{d}\omega\\
&\qquad =
\iint_{\Omega\times\Omega}
\left[
 A^v\,\partial_v\varphi
+
\frac{\nu^2}{2}D(v)^2\partial^2_{vv}\varphi
\right]
f_\varepsilon(v,\sigma,u,t)
f_\varepsilon(v_*,\sigma_*,u_*,t)
\,\mathrm{d}\omega\,\mathrm{d}\omega_*
+
\mathcal{O}(\sqrt{\varepsilon}).
\end{aligned}
\]
Passing to the limit $\varepsilon\to0$, we  formally obtain the following operator in weak form
\begin{equation*}
\int_{\Omega}
\left[
\mathcal{A}^v[f]\partial_v\varphi
+
\frac{\nu^2}{2}D(v)^2\partial^2_{vv}\varphi
\right]
f(v,\sigma,u,t)\,\mathrm{d}\omega,
\end{equation*}
where the non-local macroscopic opinion drift $\mathcal{A}^v[f]$ is obtained by integrating $A^v$ over the interacting variables:
\begin{align*}
\mathcal{A}^v[f](v,\sigma,u,t)
    &=
    \Theta(\sigma)(m_v(t)-v)
    +
    (1-\Theta(\sigma))(u-v)\cr
    &=
    \Theta(\sigma)(m_v(t)-u)
    +
    u-v,
\end{align*}
with $m_v(t)=\int_\Omega v f(v,\sigma,u,t)\,\mathrm{d}\omega$ denoting the mean public opinion.

\paragraph{Grazing limit for the strategy operator.}
Similarly, the scaled strategy dynamics  \eqref{eq:binary-game} reads:
\begin{equation}\label{eq:strategy-scaled}
\sigma_i'-\sigma_i
 =
\varepsilon A_i^\sigma(v,\sigma,u;v_*,\sigma_*,u_*),
\qquad i=1,\dots,K,
\end{equation}
where \(  A_i^\sigma =\,\sigma_i(P_i-\overline{P})\). 

Expanding $\varphi(v,\sigma',u)$ around $(v,\sigma,u)$ gives
\[
\begin{aligned}
\varphi(v,\sigma',u)-\varphi(v,\sigma,u)
&=
\sum_{i=1}^K
\partial_{\sigma_i}\varphi(v,\sigma,u)(\sigma_i'-\sigma_i)
+
R_\varepsilon^\sigma .
\end{aligned}
\]
Using \eqref{eq:strategy-scaled}, we obtain
\begin{equation*}
\varphi(v,\sigma',u)-\varphi(v,\sigma,u)
=
\varepsilon
\sum_{i=1}^K
A_i^\sigma\partial_{\sigma_i}\varphi(v,\sigma,u)
+
R_\varepsilon^\sigma.
\end{equation*}
Since $\varphi\in C_b^3$ and because the strategy update is deterministic, and assuming the payoff differences $P_i-\overline{P}$ are bounded, the remainder satisfies
\[
    |R_\varepsilon^\sigma|
    \leq
    C_\varphi |\sigma'-\sigma|^2 = \varepsilon^2\sum_{i=1}^K |A_i^\sigma|^2=
    \mathcal O(\varepsilon^2).
\]

Consequently, the limiting weak contribution of the strategy operator can be written as
\[
\begin{aligned}
&\frac{1}{\varepsilon}
\int_{\Omega} Q^\sigma(f_\varepsilon,f_\varepsilon)(v,\sigma,u,t)
\varphi(v,\sigma,u)\,\mathrm{d}\omega\\
&\qquad =
\iint_{\Omega\times\Omega}
\sum_{i=1}^K
\,\sigma_i(P_i-\overline{P})
\partial_{\sigma_i}\varphi(v,\sigma,u)
f_\varepsilon(v,\sigma,u,t)
f_\varepsilon(v_*,\sigma_*,u_*,t)
\,\mathrm{d}\omega\,\mathrm{d}\omega_*
+
\mathcal O(\varepsilon).
\end{aligned}
\]
Passing to the limit $ \varepsilon\to 0 $, we formally obtain the following operator in weak form:
\[
   \int_{\Omega}
    \sum_{i=1}^K
    \mathcal A_i^\sigma[f](v,\sigma,u,t)
    \partial_{\sigma_i}\varphi(v,\sigma,u)
    f(v,\sigma,u,t)
    \,\mathrm{d}\omega,
\]
where the nonlocal replicator drift is:
\begin{equation}\label{eq:Asigma-drift}
   \mathcal A_i^\sigma[f](v,\sigma,u,t)
    =
 \,\sigma_i
    \left(
    \mathcal P_i[f](v,,u,t)
    -
    \overline{\mathcal P}[f](v,u,t)
    \right),
    \qquad i=1,\dots,K,
\end{equation}
with the macroscopic expected payoffs given by:
\begin{equation*}
    \mathcal P_i[f](v,u,t)
    =
    \int_{\Omega}
    P_i(v,u;v_*,u_*,\sigma_*)
    f(v_*,\sigma_*,u_*,t)\,\mathrm{d}\omega_*, \qquad
    \overline{\mathcal P}[f](v,\sigma,u,t)
    =
    \sum_{i=1}^K
    \sigma_i\mathcal P_i[f](v,u,t).
\end{equation*}
\paragraph{Vlasov--Fokker--Planck model.}
Combining the opinion and strategy limits, we obtain the weak form
\begin{equation}\label{eq:limit-weak}
\frac{d}{dt}
\int_\Omega \varphi(v,\sigma,u)f(v,\sigma,u,t)\,\mathrm{d}\omega
=
\frac{1}{\tau_v}
\int_\Omega
\left[
\mathcal{A}^v[f]\partial_v\varphi
+
\frac{\nu^2}{2}D(v)^2\partial^2_{vv}\varphi
\right]
f\,\mathrm{d}\omega+
\frac{1}{\tau_\sigma}
\int_\Omega
\sum_{i=1}^K
\mathcal{A}^\sigma_i[f]\partial_{\sigma_i}\varphi
f\,\mathrm{d}\omega.
\end{equation}

Integrating by parts, and using the zero-flux conditions at $v=\pm1$ and on $\partial\Delta_{K-1}$, the boundary terms vanish and we recover the strong formulation of the Vlasov--Fokker--Planck (VFP)  governing the joint density $f$:
\begin{equation}\label{eq:FP_final}
\partial_t f
=
-\frac{1}{\tau_v}\partial_v\left(\mathcal{A}^v[f]f
-
\frac{\nu^2}{2} \partial_{v} (D(v)^2f)\right)
-
\frac{1}{\tau_\sigma}
\sum_{i=1}^K
\partial_{\sigma_i}
\left(
\mathcal{A}^\sigma_i[f]f
\right),
\end{equation}
supplemented with the initial condition
\begin{equation*}
f(v,\sigma,u,0)=f_0(v,\sigma,u),
\qquad
(v,\sigma,u)\in\Omega.
\end{equation*}

Equation \eqref{eq:FP_final} is complemented by zero-flux boundary conditions in the opinion variable and on the
boundary of the strategy simplex. In the opinion variable, this reads
\[
\left.
\left(
\frac{1}{\tau_v}\mathcal{A}^v[f]f
-
\frac{\nu^2}{2\tau_v}
\partial_v(D(v)^2f)
\right)
\right|_{v=\pm1}
=0.
\] Since $D(\pm1)=0$ by condition (D1) of Proposition~\ref{prop:domain_invariance}, the diffusion flux vanishes at $v=\pm1$, and the zero-flux condition reduces to $\mathcal{A}^v[f]f|_{v=\pm1}=0$.
On the boundary $\partial\Delta_{K-1}$ of the simplex, one imposes
\begin{equation}\label{eq:0_f_cond}
\left(\mathcal{A}^\sigma[f]f\right)\cdot n_\sigma=0.
\end{equation}
On the face $\{\sigma_i=0\}$ of $\partial\Delta_{K-1}$, the outward normal has $i$-th component $-1$, and  $\mathcal{A}_i^\sigma[f]\propto\sigma_i$ vanishes identically there. Hence the zero-flux condition is automatically satisfied by the replicator structure of the drift, for each $i=1,\ldots,K$.

\paragraph{Formulation in terms of relative payoff gains.} 
For analytical tractability, it is convenient to reformulate the strategic dynamics using the $K-1$ independent coordinates $\widehat{\sigma}$ in $\Delta_{K-1}$. Taking strategy $K$ as the reference, we define the relative payoff gains:
\begin{equation}\label{eq:macrogain}
\mathcal G_i[f](v,u,t)
:=
\mathcal P_i[f](v,u,t)
-
\mathcal P_K[f](v,u,t),
\qquad
i=1,\ldots,K-1.
\end{equation}
Rewriting the average payoff  as
\[
\overline{\mathcal P}[f](v,\widehat \sigma, u,t)
=
\mathcal P_K[f]
+
\sum_{j=1}^{K-1}\sigma_j\mathcal G_j[f],
\]
the first $K-1$ components of the replicator field \eqref{eq:Asigma-drift} can then be expressed in reduced coordinates $\widehat{\sigma}$, as
\begin{equation*}
\mathcal A_i^\sigma[f](v,\widehat \sigma, u,t)
=
\sigma_i
\left(
\mathcal G_i[f](v,u,t)
-
\sum_{j=1}^{K-1}\sigma_j\mathcal G_j[f](v,u,t)
\right).
\end{equation*}
In this coordinate system, the probability density becomes $f=f(v,\widehat{\sigma},u,t)$, and the VFP equation \eqref{eq:FP_final} takes the equivalent form:
\begin{equation}\label{eq:FP_gain}
\partial_t f
=
-\frac{1}{\tau_v}
\partial_v
\left(
\mathcal A^v[f]f
-
\frac{\nu^2}{2}
\partial_v\left(D(v)^2f\right)
\right)
-
\frac{1}{\tau_\sigma}
\sum_{i=1}^{K-1}
\partial_{\sigma_i}
\left(
\sigma_i
\left(
\mathcal G_i[f]
-
\sum_{j=1}^{K-1}
\sigma_j\mathcal G_j[f]
\right)f
\right).
\end{equation}

\begin{remark}[Reduced formulation in terms of payoff gains]
The reduced formulation \eqref{eq:FP_gain} can be obtained directly by
starting from the replicator dynamics written in terms of the payoff
gains, as in \eqref{eq:binary-game-gain}, and applying the same
scaling and kinetic-limit procedure described above. This yields the
Vlasov--Fokker--Planck equation directly in reduced strategy
coordinates $\widehat \sigma$, without the need to first derive the equation in the
full simplex and subsequently rewrite it in reduced coordinates.
\end{remark}

\subsection{Marginal and macroscopic dynamics}
To characterize the collective dynamics, we study the marginal distributions of the opinion, conviction, and strategy variables
\begin{equation*}
\begin{aligned}
h(v,t)
=
\iint_{\Delta_{K-1} \times I}
f(v,\sigma,u,t)\,\mathrm{d}\sigma\,\mathrm{d}u,\, k(u,t) = \iint_{I\times\Delta_{K-1}}
f(v,\sigma,u,t)\,\mathrm dv\,\mathrm{d}\sigma,
\, g(\sigma,t)
=
\iint_{I^2}
f(v,\sigma,u,t)\,\mathrm dv\,\mathrm du.
\end{aligned}
\end{equation*}

We also introduce the macroscopic averages: the mean opinion $m_v(t)$, the mean conviction $m_u(t)$, and the vector of mean strategies $m_\sigma(t) = \left(m_{\sigma_1}(t),\ldots,m_{\sigma_K}(t)\right)^\top$
as
\begin{equation*}
m_v(t) = \int_\Omega v f(v,\sigma,u,t)\,\mathrm{d}\omega, \quad m_u(t) = \int_\Omega u f(v,\sigma,u,t)\,\mathrm{d}\omega, \quad
m_{\sigma_i}(t) = \int_\Omega \sigma_i f(v,\sigma,u,t)\,\mathrm{d}\omega. 
\end{equation*}

At the level of the binary interaction rules, the conviction $u$ is assumed to remain unchanged  and therefore acts as a quenched variable. The following proposition shows that this property is retained by the resulting kinetic equation, so that both the conviction marginal and its mean remain constant in time.

\begin{proposition}[Stationarity of the conviction marginal and
conservation of its mean]
\label{prop:conviction-stationary}
Let $f$ be a solution of \eqref{eq:FP_final} satisfying the weak
formulation \eqref{eq:limit-weak}. Then the conviction marginal is independent of time, that is,
\[
k(u,t)=k(u):=k(u,0),\qquad \forall t\geq0,
\]
and the mean conviction is also conserved:
\[
m_u(t)=m_u:=m_u(0),
\qquad \forall t\geq0.
\]
\end{proposition}

\begin{proof}
{Since $k(u,t)$ does not depend on $\sigma$ nor $v$ we} let $\psi\in C^\infty(I)$ and choose the test function
$\varphi(v,\sigma,u)=\psi(u)$ in the weak formulation
\eqref{eq:limit-weak}. Since $\psi$ is independent of $v$
and $\sigma$, all the derivatives appearing on the right-hand
side vanish. Hence,
\[
\frac{d}{dt}\int_I \psi(u)k(u,t)\,\mathrm{d}u=0.
\]
Since this identity holds for every test function
$\psi$, the conviction marginal is constant in time, in the
sense of measures. In particular, $k(u,t)=k_0(u)$ almost
everywhere whenever the marginal admits a density. Finally, choosing $\psi(u)=u$ gives
\[
m_u(t)
=\int_I u\,k(u,t)\,\mathrm{d}u
=\int_I u\,k(u)\,\mathrm{d}u
=m_u(0),
\]
which proves the conservation of the mean conviction.
\end{proof}

\subsubsection{Opinion marginal and mean dynamics}
Integrating the VFP equation \eqref{eq:FP_final} with respect to $\sigma$ and $u$, the strategy flux term vanishes
due to the zero-flux boundary condition \eqref{eq:0_f_cond},  yielding the continuity equation for the opinion marginal:
\begin{equation}\label{eq:opinion-marginal}
\partial_t h(v,t)
+
\frac{1}{\tau_v}\partial_v
\left(
\mathcal V^{v}[f](v,t)h(v,t)
\right)
=
\frac{\nu^2}{2\tau_v}
\partial_{vv}\left(D(v)^2h(v,t)\right),
\end{equation}
where for $h(v,t)>0$,  the effective macroscopic opinion velocity is given by the conditional expectation of the drift:
\begin{align*}
\mathcal V^{v}[f](v,t)
&:=
\frac{1}{h(v,t)}
\iint_{\Delta_{K-1}\times I}
\mathcal A^v[f](v,\sigma,u,t)
f(v,\sigma,u,t)\,\mathrm{d}\sigma\,\mathrm{d}u
\nonumber\\
&=
\frac{1}{h(v,t)}
\iint_{\Delta_{K-1}\times I}
\left[
\Theta(\sigma)(m_v(t)-u)
+u-v
\right]
f(v,\sigma,u,t)\,\mathrm{d}\sigma\,\mathrm{d}u.
\end{align*} Although~\eqref{eq:opinion-marginal} has the structure of a one-dimensional Fokker--Planck equation, the effective velocity $\mathcal{V}^v[f]$ depends on the conditional statistics of the full joint density $f$ and cannot be expressed in terms of $h$ alone; the marginal equation is therefore not closed.

The evolution of the first moment is retrieved  multiplying \eqref{eq:FP_final} by $v$ and integrating over $\Omega$, the diffusion term and the strategy flux contribute zero: the former vanishes after double integration by parts in $v$ using the zero-flux boundary conditions at $v=\pm1$, while the latter vanishes because $v$ is independent of $\sigma$. Integration by parts of the opinion drift gives
\begin{align*}
    \frac{d}{dt}m_v(t) 
    &= \frac{1}{\tau_v} \int_\Omega \mathcal{A}^v[f] f \,\mathrm{d}\omega 
    \\
    &= \frac{1}{\tau_v} \int_\Omega \left[ \Theta(\sigma)(m_v(t) - u) + u - v \right] f \,\mathrm{d}\omega \\
    &= \frac{1}{\tau_v} \left( m_v(t)\sum_{i=1}^K\theta_i m_{\sigma_i}(t) - \sum_{i=1}^K \theta_i \left(\int_\Omega \sigma_i u f \,\mathrm{d}\omega\right) + m_u - m_v(t)\right).
\end{align*} 

Introducing the strategy--conviction covariance vector: 
\begin{equation}\label{eq:cov_sigma_u}
    \Sigma_{\sigma u}(t)
    :=
    \int_\Omega
    \bigl(\sigma-m_\sigma(t)\bigr)
    \bigl(u-m_u\bigr)f\,\mathrm{d}\omega
    \in\mathbb R^K,
\end{equation}
we obtain the following equation for the mean opinion:
\begin{equation}\label{eq:mean_opinion}
    \frac{d}{dt}m_v(t) = \frac{1}{\tau_v} \left( \left(1 - {\theta}^\top m_\sigma(t)\right)\left(m_u - m_v(t)\right) - \theta^\top\Sigma_{\sigma u}(t) \right).
\end{equation}
Alternatively, the evolution of $m_v(t)$ can be obtained by choosing $\varphi(v,\sigma,u)=v$ in the weak formulation \eqref{eq:limit-weak}.

Equation \eqref{eq:mean_opinion} is not closed since it depends on the correlation between strategies and intrinsic convictions. A simplification occurs when $\sigma$ and $u$ are statistically independent, in this case the covariance term vanishes and the mean opinion relaxes towards $m_u$ at a rate determined by the population's average stubbornness $(1-{\theta}^\top m_\sigma(t))$. This holds, for example, for constant payoffs; more generally, independence between $\sigma$ and $u$ is preserved under~\eqref{eq:FP_final} whenever it holds initially and the gains $\mathcal G_i[f]$ are independent of $(v,u)$, as for constant payoffs. If they depended on $v$ alone, the conditional law of $v$ given $(\sigma,u)$, which depends on $u$ through the attraction center, would correlate $\sigma$ and $u$.
For state-dependent payoffs, instead, the strategy dynamics may generate non-trivial correlations with $u$, producing a systematic displacement of the mean opinion from $m_u$, or drive persistent macroscopic oscillations of $m_v(t)$.

\subsubsection{Strategy marginal and mean dynamics}
Integrating \eqref{eq:FP_final} with respect to $v$ and $u$, we obtain the evolution equation for the strategy marginal:
\begin{equation}\label{eq:strategy-marginal}
\partial_t g(\sigma,t)
+
\frac{1}{\tau_\sigma}\sum_{i=1}^K
\partial_{\sigma_i}
\left(
\mathcal V^\sigma_i[f](\sigma,t)g(\sigma,t)
\right)
=
0,
\end{equation}
where, for $g(\sigma,t)>0$, 
introducing the conditional payoffs
\[
\mathcal P_{i|\sigma}[f](t)
:=
\frac{1}{g(\sigma,t)}
\iint_{I^2}
\mathcal P_i[f]\,f\,\mathrm dv\,\mathrm du,
\qquad
\overline{\mathcal P}_{|\sigma}[f]
=
\sum_{j=1}^K\sigma_j\mathcal P_{j|\sigma}[f],
\]
the advective velocity fields takes the replicator form
\[
\mathcal V^\sigma_i[f](\sigma,t)
=
\sigma_i
\left(
\mathcal P_{i|\sigma}[f]
-
\overline{\mathcal P}_{|\sigma}[f]
\right).
\]

Equation \eqref{eq:strategy-marginal} nevertheless is
not closed in $g$, since the conditional payoffs depend on the full joint density $f$, and cannot be expressed in terms of $g$ alone without further assumptions. 

\begin{remark}[Replicator structure of the strategy marginal]
The replicator form of the velocity field implies that $\mathcal{V}^\sigma_i[f]=0$ on $\{\sigma_i=0\}$, 
so the simplex faces are forward-invariant. Furthermore,  $\sum_{i=1}^K\mathcal{V}^\sigma_i[f]=0$, 
so the flow is tangent to $\Delta_{K-1}$. 
Consequently, an initial distribution supported in 
$\Delta_{K-1}$ remains supported there for all $t\geq0$.
Moreover, under the
corresponding no-flux condition on the boundary of the simplex, we have conservation of the total density
\[
    \frac{d}{dt}
    \int_{\Delta_{K-1}}
    g({\sigma},t)\,\mathrm{d}\sigma
    =
    0.
\]\end{remark}

The evolution of the first moment is recovered multiplying \eqref{eq:FP_final} by $\sigma_i$ and integrating over $\Omega$. The opinion flux vanishes (since $\sigma_i$ is independent of $v$) and
integration by parts in $\sigma_i$ with zero-flux boundary conditions on
$\partial\Delta_{K-1}$ gives
\begin{equation}\label{eq:mean_strategy}
    \frac{d}{dt}m_{\sigma_i}(t) 
    = \frac{1}{\tau_\sigma} \int_\Omega \mathcal{A}^\sigma_i[f] f \,\mathrm{d}\omega 
    = \frac{1}{\tau_\sigma} \int_\Omega \sigma_i \Big( \mathcal{P}_i[f](v,\sigma,u,t) - \overline{\mathcal{P}}[f](v,\sigma,u,t) \Big) f \,\mathrm{d}\omega,\quad i=1,\ldots,K.
\end{equation}

The evolution of the mean \eqref{eq:mean_strategy}, which defines a macroscopic replicator-type dynamics for the mean strategy,  is generally not closed because of its dependence on the joint distribution.

\paragraph{Constant payoff case.} 
A fundamental simplification arises when the payoff matrix is constant, $J\in\mathbb{R}^{K\times K}$, implying that the interactions are independent of the states of the agents. In this case, the expected payoff vector depends linearly on the strategic mean \(m_\sigma(t)\), while the average payoff is linear in the individual strategy \(\sigma\) and depends on \(m_\sigma(t)\): 
\[ \mathcal P_i[f](t)=(Jm_\sigma(t))_i,\qquad
    \overline{\mathcal P}[f]({\sigma},t)
    =
    {\sigma}^{\top}
    J m_\sigma(t).
\] The velocity field is then given by
\[
    \mathcal V^\sigma[f](\sigma,t)
    =
    M({\sigma})
    J m_\sigma(t),
    \qquad
    M(\sigma)
    :=
    \operatorname{diag}({\sigma})
    -
    {\sigma}{\sigma}^{\top},
\]
with $M({\sigma})$ denoting the replicator (or mobility) matrix on the
simplex. Hence the marginal equation \eqref{eq:strategy-marginal}  yields a closed equation for $g(\sigma,t)$ given by

\begin{equation}\label{eq:strategy-marginal-closed}
\partial_t g(\sigma,t)
+
\frac{1}{\tau_\sigma}\nabla_\sigma\cdot \left( M(\sigma)Jm_\sigma(t)g(\sigma,t) \right) = 0,
\end{equation}

Since \(M(\sigma)\mathbf 1=0\), the velocity is tangent to the hyperplane \(\sum_{i=1}^K\sigma_i=1\). Moreover, \((\mathcal V^\sigma [g])_i=0\) whenever \(\sigma_i=0\), so the vector field preserves the simplex.

Furthermore, the evolution of the average strategy \eqref{eq:mean_strategy} is
\begin{equation}\label{eq:mean_strategy_constant}
\frac{d}{dt}m_{\sigma_i}(t)
=
\frac{1}{\tau_\sigma}
\left[
m_{\sigma_i}(t)(Jm_\sigma(t))_i
-
\sum_{j=1}^K
\left(\int_\Omega \sigma_i\sigma_j f \,\mathrm{d}\omega\right) (Jm_\sigma(t))_j
\right],
\qquad i=1,\ldots,K,
\end{equation}
Equivalently, introducing the
strategy covariance matrix
\[
\Sigma_\sigma(t) := \int_\Omega (\sigma-m_\sigma(t))(\sigma -m_\sigma(t))^\top f\,\mathrm{d}\omega\in\mathbb{R}^{K\times K},
\]
equation \eqref{eq:mean_strategy_constant} can be written in vector form as
\[
\frac{d}{dt}m_\sigma(t)
=
\frac{1}{\tau_\sigma}
\left[
M(m_\sigma(t))
-
\Sigma_\sigma(t)
\right]Jm_\sigma(t).
\]
In the absence of the covariance term, the equation formally reduces to the classical replicator equation for the strategy vector \(m_\sigma\). In the present kinetic setting, however, \(m_\sigma\) represents the population mean of the individual strategy distributions, rather than the fractions of agents adopting the corresponding pure strategies. The covariance term accounts for the effect of strategy heterogeneity on the evolution of the mean.

\begin{remark}[Moment hierarchy and closure]
The equations \eqref{eq:mean_opinion}-\eqref{eq:mean_strategy} for $m_v$ and $m_\sigma$ show that the macroscopic
dynamics is not closed at the level of first moments. In particular, there is a dependency on the covariance matrix
$\Sigma_{\sigma}(t)$, and vector $\Sigma_{\sigma u}(t)$, and this dependence is not solved even for state-independent payoffs in \eqref{eq:mean_strategy_constant}. 
In general the kinetic model generates a hierarchy of correlations between opinion, conviction, and strategy variables, and a closed macroscopic description requires additional assumptions. 
However, this does not imply nonclosure of the full strategy marginal: as shown below, for constant payoffs the marginal $g$ admits an exact autonomous transport equation, without any moment closure approximation, and an explicit projective representation.
\end{remark}

\subparagraph{Exact representation of the strategy marginal.}\label{sec:strategy-flow}

Equation \eqref{eq:strategy-marginal-closed} is  a non-linear transport equation in $g$, and its solution provides an exact representation of the strategy marginal and of all its moments in terms of the initial distribution and the cumulative payoff gains. To derive this representation, it is useful to rewrite equation~\eqref{eq:strategy-marginal-closed} in terms of the relative payoff gains introduced in \eqref{eq:macrogain}. 

Since $J$ is constant, $G_i(\widehat\sigma_*)$ is independent 
of $(v_*,u_*)$, and integration of $f$ over $(v_*,u_*)\in I^2$ 
yields the strategy marginal $g$:
\begin{align*}
\mathcal G_i[f](t)
&= \int_\Omega G_i(\widehat\sigma_*)\,f(v_*,\widehat \sigma_*,u_*,t)
\,\mathrm{d}\omega_*
= \int_{\Delta_{K-1}}G_i(\widehat\sigma_*)\,g(\widehat\sigma_*,t)
\,d\widehat\sigma_*.
\end{align*}
Moreover, since in the constant-payoff setting $G_i$ is affine in the opponent's reduced strategy  $\widehat\sigma_*$, its integral 
against a probability measure equals its value at the mean:
\[
\int_{\Delta_{K-1}}G_i(\widehat\sigma_*)\,g(\widehat\sigma_*,t)
\,d\widehat\sigma_* = G_i(\widehat m_\sigma(t)),
\]
    where $\widehat m_\sigma(t)=(m_{\sigma_1}(t),\ldots,m_{\sigma_{K-1}}(t))^\top$ is the reduced mean strategy. Thus, the relative payoff gain is an affine function of $\widehat m_\sigma(t)$. In particular, since $\mathcal G_i[f](t)$ is independent of $(v,u)$ and depends on $f$ only through  $\widehat m_\sigma(t)$, the strategy marginal $g(\widehat\sigma,t)$ 
evolves autonomously from the opinion and conviction variables. Consequently, \eqref{eq:strategy-marginal-closed}  can be rewritten in reduced coordinates as
\begin{equation}\label{eq:strategy-marginal-gain}
\partial_t g(\widehat \sigma,t)+
\frac{1}{\tau_\sigma}\nabla_{\widehat \sigma}\cdot \left(  M(\widehat \sigma)G(\widehat m_\sigma(t))g(\widehat\sigma,t) \right) = 0.
\end{equation}
Here,
$G =(G_1, \ldots, G_{K-1})^\top$ denotes the vector of individual gain functions, evaluated at  $\widehat m_\sigma(t)$,  and  $ M(\widehat \sigma)=\operatorname{diag}({\widehat \sigma})-{\widehat \sigma}{\widehat \sigma}^{\top}$ is the reduced replicator matrix. This formulation separates the geometry of the simplex, encoded by the reduced replicator matrix $M(\widehat \sigma)$, from the selective force, represented by the relative payoff gain vector $ G$. 

In the general state-dependent setting, the same reduction to independent simplex coordinates remains valid. However, the relative gains may depend on the individual strategy $\widehat{\sigma}$, on $(u,v)$, and on the full joint density $f$, rather than solely on its strategy marginal $g$. In that case, the resulting transport equation for $g$ is generally not closed.

The autonomous structure of~\eqref{eq:strategy-marginal-gain} allows its solution to be determined exactly via the method of characteristics, and analogously to the approach developed in
\cite{LahkarSeymour2013}, as formalized in the following theorem.

\begin{theorem}[Exact projective representation of the strategy 
marginal]\label{thm:projective}
Let $J\in\mathbb{R}^{K\times K}$ be a constant payoff matrix,
and let $g_0$ be a probability measure on $\Delta_{K-1}$
in reduced coordinates. The unique solution of
\eqref{eq:strategy-marginal-gain} in the sense of measures,
with initial datum $g_0$, is given by the push-forward
\begin{equation}\label{eq:strategy-projective-pushforward}
    g(\cdot,t)=\bigl(\Phi_{\mathcal{T}(t)}\bigr)_\#g_0,
\end{equation}
where $\Phi_{\mathcal{T}}:\Delta_{K-1}\to\Delta_{K-1}$ is
the homeomorphism defined component-wise, for
$i=1,\ldots,K-1$, by
\begin{equation}\label{eq:strategy-projective-map}
    \Phi_{\mathcal{T}}^i(\widehat{\sigma})
    :=
    \frac{\mathrm{e}^{\mathcal{T}_i}\sigma_i}
    {\sigma_K+\displaystyle\sum_{j=1}^{K-1}
    \mathrm{e}^{\mathcal{T}_j}\sigma_j},
\end{equation}
where $\sigma_K:=1-\sum_{j=1}^{K-1}\sigma_j$ and we write
$\sigma_i$ for the $i$-th component of $\widehat{\sigma}$.
The projective map satisfies $\Phi_{\mathcal{T}}^{-1}=\Phi_{-\mathcal{T}}$,
and $g(\cdot,t)$ is a probability measure on $\Delta_{K-1}$
for all $t\geq0$.

The cumulative relative payoff gain vector
$\mathcal{T}(t)=(\mathcal{T}_1(t),\ldots,\mathcal{T}_{K-1}(t))
^\top\in\mathbb{R}^{K-1}$, defined by
\begin{equation}\label{cumulative_gains}
    \mathcal{T}_i(t)
    :=
    \frac{1}{\tau_\sigma}
    \int_0^t G_i\bigl(\widehat{m}_\sigma(s)\bigr)\,\mathrm{d}s,
    \qquad i=1,\ldots,K-1,
\end{equation}
with $G_i$ as in~\eqref{def:Gi}, is the unique global
solution of the autonomous ODE system
\begin{equation}\label{eq:strategy-cumulative-gains}
    \mathcal{T}_i'(t)
    =
    \frac{1}{\tau_\sigma}
    G_i\bigl(\mathcal{E}(\mathcal{T}(t))\bigr),
    \qquad
    \mathcal{T}_i(0)={0},
    \qquad i=1,\ldots,K-1,
\end{equation}
where the reduced mean strategy map
$\mathcal{E}:\mathbb{R}^{K-1}\to\Delta_{K-1}$ is defined by
\begin{equation}\label{def:E}
    \mathcal{E}(\mathcal{T})
    :=
    \int_{\Delta_{K-1}}
    \Phi_{\mathcal{T}}(\widehat{\sigma})\,
    g_0(\mathrm{d}\widehat{\sigma}).
\end{equation}
In particular,
\(\widehat{m}_\sigma(t)=\mathcal{E}(\mathcal{T}(t))\).
\end{theorem}

\begin{proof}
By the analysis of the preceding paragraph,
$\mathcal{G}_i[f](t)=G_i(\widehat{m}_\sigma(t))$
depends on $f$ only through $\widehat{m}_\sigma(t)$, so
equation~\eqref{eq:strategy-marginal-gain} is an autonomous
transport equation in $g$ with time-dependent velocity
$M(\widehat{\sigma})G(\widehat{m}_\sigma(t))$.
We solve it via the method of characteristics.

\medskip
\noindent\textit{Step 1: Characteristic flow and projective map.}
The characteristic system for~\eqref{eq:strategy-marginal-gain}
is
\[
    \frac{dS}{dt}
    =
    \frac{1}{\tau_\sigma}M(S)\,G(\widehat{m}_\sigma(t)),
    \qquad
    S(0;\widehat{\sigma}_0)=\widehat{\sigma}_0,
\]
where $S(t;\widehat{\sigma}_0)\in\Delta_{K-1}$ denotes the
characteristic starting from $\widehat \sigma_0=(\sigma_{0,1},\ldots,\sigma_{0,K-1})^\top\in\Delta_{K-1}$.
In components, for $i=1,\ldots,K-1$:
\begin{equation}\label{csystem}
    \frac{dS_i}{dt}
    =
    \frac{1}{\tau_\sigma}S_i
    \biggl(G_i(\widehat{m}_\sigma(t))
    -\sum_{j=1}^{K-1}S_j G_j(\widehat{m}_\sigma(t))\biggr),
    \qquad
    S_i(0;\widehat{\sigma}_0)=\sigma_{0,i}.
\end{equation}
We set $\sigma_{0,K}:=1-\sum_{j=1}^{K-1}\sigma_{0,j}$.
The reference quantity $S_K(t;\widehat{\sigma}_0)
:=1-\sum_{i=1}^{K-1}S_i$ is a derived quantity, not a
component of the reduced characteristic; differentiating
and using~\eqref{csystem} gives
\[
    \frac{dS_K}{dt}
    =
    -\frac{1}{\tau_\sigma}S_K
    \sum_{j=1}^{K-1}S_jG_j(\widehat{m}_\sigma(t)).
\]
For characteristics starting in the interior of the simplex
($\sigma_{0,i}>0$ for all $i$), we compute the evolution
of the log-ratio $\ln(S_i/S_K)$:
\[
    \frac{d}{dt}\ln\frac{S_i}{S_K}
    =
    \frac{\dot{S}_i}{S_i}-\frac{\dot{S}_K}{S_K}
    =
    \frac{1}{\tau_\sigma}G_i(\widehat{m}_\sigma(t)),
    \qquad i=1,\ldots,K-1.
\]
Integrating over $[0,t]$ and using the
definition~\eqref{cumulative_gains}:
\begin{equation*}
\frac{S_i(t;\widehat{\sigma}_0)}{S_K(t;\widehat{\sigma}_0)}
    =
    \frac{\sigma_{0,i}}{\sigma_{0,K}}\,\mathrm{e}^{\mathcal{T}_i(t)},
    \qquad i=1,\ldots,K-1.
\end{equation*}
Solving for $S_i$ using the constraint
$\sum_{i=1}^K S_i=1$ yields the projective
map~\eqref{eq:strategy-projective-map}:
\[
    S_i(t;\widehat{\sigma}_0)
    =
    \Phi_{\mathcal{T}(t)}^i(\widehat{\sigma}_0)
    =
    \frac{\mathrm{e}^{\mathcal{T}_i(t)}\sigma_{0,i}}
    {\sigma_{0,K}+\displaystyle\sum_{j=1}^{K-1}
    \mathrm{e}^{\mathcal{T}_j(t)}\sigma_{0,j}}.
\]
The map $\Phi_{\mathcal{T}}$ is well-defined and continuous
on all of ${\Delta_{K-1}}$. The denominator satisfies
\begin{equation}\label{eq:denom-bound}
    \sigma_K+\sum_{j=1}^{K-1}\mathrm{e}^{\mathcal{T}_j}\sigma_j
    \geq
    \min\bigl\{1,\mathrm{e}^{\mathcal{T}_1},\ldots,
    \mathrm{e}^{\mathcal{T}_{K-1}}\bigr\}
    \cdot\underbrace{\biggl(\sigma_K+
    \sum_{j=1}^{K-1}\sigma_j\biggr)}_{=1}
    >0
\end{equation}
for every finite $\mathcal{T}\in\mathbb{R}^{K-1}$, and the formula extends continuously to
$\partial\Delta_{K-1}$. In particular, if $\sigma_{0,K}=0$, the reconstructed component $S_K(t;\widehat{\sigma}_0)$ remains zero, consistently with the invariance of the corresponding face under the
replicator flow. The inverse map $\Phi_{\mathcal{T}}^{-1}
=\Phi_{-\mathcal{T}}$ is verified by direct substitution,
and is also continuous by~\eqref{eq:denom-bound} applied
with $-\mathcal{T}$. Therefore $\Phi_{\mathcal{T}}$ is a
homeomorphism of ${\Delta_{K-1}}$ onto itself.

\medskip
\noindent\textit{Step 2: Push-forward satisfies the transport equation.}
We verify that $g(t,\cdot)=(\Phi_{\mathcal{T}(t)})_\#g_0$
satisfies~\eqref{eq:strategy-marginal-gain} in the weak
sense. Since $g_0$ may have mass on $\partial\Delta_{K-1}$,
we test against $\varphi\in C^1({\Delta_{K-1}})$.
For any such $\varphi$:
\begin{align*}
    \frac{d}{dt}\int_{\Delta_{K-1}}\varphi\,dg(t)
    &=
    \frac{d}{dt}\int_{\Delta_{K-1}}
    \varphi\bigl(\Phi_{\mathcal{T}(t)}(\widehat{\sigma}_0)\bigr)
    \,g_0(\mathrm{d}\widehat{\sigma}_0)
    \\
    &=
    \frac{1}{\tau_\sigma}
    \int_{\Delta_{K-1}}
    \nabla\varphi\bigl(\Phi_{\mathcal{T}(t)}(\widehat{\sigma}_0)\bigr)
    \cdot M\bigl(\Phi_{\mathcal{T}}(\widehat{\sigma}_0)\bigr)
    G(\widehat{m}_\sigma(t))\,g_0(\mathrm{d}\widehat{\sigma}_0)
    \\
    &=
    \frac{1}{\tau_\sigma}
    \int_{\Delta_{K-1}}
    \nabla\varphi(\widehat{\sigma})\cdot
    M(\widehat{\sigma})\,G(\widehat{m}_\sigma(t))\,
    g(\mathrm{d}\widehat{\sigma},t),
\end{align*}
where the second line uses the chain rule and
the characteristic equation, and the third follows by the
change of variables
$\widehat{\sigma}=\Phi_{\mathcal{T}(t)}(\widehat{\sigma}_0)$. The characteristic velocity is tangent to every face of the simplex, so its normal component vanishes on the boundary. Hence, the above identity gives the weak formulation of
\eqref{eq:strategy-marginal-gain}.

\medskip
\noindent\textit{Step 3: Self-consistency and well-posedness.}
Since $\widehat{m}_\sigma(t)=\mathcal{E}(\mathcal{T}(t))$
by definition~\eqref{def:E}, differentiating~\eqref{cumulative_gains}
and using $\mathcal{G}_i[f](t)=G_i(\mathcal{E}(\mathcal{T}(t)))$
directly yields the autonomous system~\eqref{eq:strategy-cumulative-gains}.

The right-hand side of~\eqref{eq:strategy-cumulative-gains}
is globally Lipschitz in $\mathcal{T}$: direct differentiation
of the projective map gives
\[
    \frac{\partial\Phi_{\mathcal{T}}^i}{\partial\mathcal{T}_j}
    =
    \Phi_{\mathcal{T}}^i\bigl(\delta_{ij}-\Phi_{\mathcal{T}}^j\bigr),
    \qquad i,j=1,\ldots,K-1,
\]
so that
$\bigl|\partial_{\mathcal{T}_j}\Phi_{\mathcal{T}}^i\bigr|
=\Phi_{\mathcal{T}}^i|\delta_{ij}-\Phi_{\mathcal{T}}^j|\leq1/4$ (indeed, for $i=j$ this follows from $x(1-x)\leq1/4$ for $x\in[0,1]$, whereas for $i\neq j$ it follows from
$\Phi_{\mathcal{T}}^i+\Phi_{\mathcal{T}}^j\leq1$ and the arithmetic--geometric mean inequality)
uniformly on $\Delta_{K-1}$, making $\mathcal{E}$ Lipschitz
in $\mathcal{T}$; since $G_i$ is affine, the
composition is Lipschitz. By the Picard--Lindel\"of theorem,
the system~\eqref{eq:strategy-cumulative-gains} admits a
unique global solution, which determines $\mathcal{T}(t)$
and hence the push-forward~\eqref{eq:strategy-projective-pushforward}
uniquely.

\medskip
\noindent\textit{Step 4: Uniqueness of the measure-valued solution.}
Conversely, any measure-valued solution $g(\cdot,t)$
of~\eqref{eq:strategy-marginal-gain} is transported by its
own characteristic flow, which satisfies the same ODE
system~\eqref{csystem} with $\widehat{m}_\sigma(t)
=\int\widehat{\sigma}\,g(d\widehat{\sigma},t)$. The
cumulative gains defined via~\eqref{cumulative_gains}
therefore satisfy~\eqref{eq:strategy-cumulative-gains}.
By uniqueness of the ODE solution established in Step~3,
these cumulative gains coincide with those of the
push-forward solution, so $g(\cdot,t)
=(\Phi_{\mathcal{T}(t)})_\#g_0$ for all $t\geq0$.

Since $\Phi_{\mathcal{T}}$ maps $\Delta_{K-1}$ to itself,
$g(\cdot,t)$ remains a probability measure on $\Delta_{K-1}$
for all $t\geq0$.
\end{proof}

\begin{remark}[Full simplex coordinates]
The projective map~\eqref{eq:strategy-projective-map} 
can be extended to the full $K$-component representation by setting
\[
\Phi^K_{\mathcal T}(\widehat\sigma)
=
\frac{\sigma_K}{\sigma_K+\sum_{j=1}^{K-1}
\mathrm{e}^{\mathcal T_j}\sigma_j},
\]
Together with the first $K-1$ components defined in
\eqref{eq:strategy-projective-map}, this gives the full strategy
vector  $(\Phi^1_{\mathcal T},\ldots,\Phi^K_{\mathcal T}):
\Delta_{K-1}\to\Delta_{K-1}$, satisfying
\(
\sum_{i=1}^K\Phi_{\mathcal T}^i(\widehat{\sigma})=1.\)
Similarly, the full mean-strategy vector $(\mathcal E_1,\ldots,\mathcal E_K)$, is obtained by adding the last component to the reduced mean-strategy map  
$$\mathcal E_K(\mathcal T)=1-\sum_{j=1}^{K-1}
\mathcal E_j(\mathcal T).$$
The full mean-strategy vector is then recovered by setting
\[
m_{\sigma_K}(t)
=
1-\sum_{j=1}^{K-1}\mathcal{E}_j(\mathcal{T}(t)),
\]
so that
\(
m_\sigma(t)
=
\left(
\mathcal{E}_1(\mathcal{T}(t)),\ldots,
m_{\sigma_K}(t)
\right)^\top.
\)
\end{remark}

Theorem~\ref{thm:projective} treats $g_0$ as a general probability measure, which is not necessarily endowed with a density. When $g_0$ is absolutely continuous, we can derive the explicit formula for the density at any given time via the Jacobian determinant, as stated in the following proposition.

\begin{proposition}[Density formula for the strategy marginal]\label{prop:density}
Under the assumptions of Theorem~\ref{thm:projective}, suppose additionally that $g_0$ is absolutely continuous with respect to the $(K-1)$-dimensional Lebesgue measure on $\Delta_{K-1}$. Then $g(\cdot,t)$ is absolutely continuous for all $t\geq0$, and its density is given by
\begin{equation}\label{eq:strategy-projective-density}
g(\widehat{\sigma},t) = g_0\!\left( \Phi_{-{\mathcal{T}}(t)}(\widehat{\sigma})\right) \frac{ \exp\!\left( -\sum_{j=1}^{K-1}\mathcal{T}_j(t) \right) }{ \left( 1-\sum_{j=1}^{K-1}\sigma_j + \sum_{j=1}^{K-1} \mathrm{e}^{-\mathcal{T}_j(t)}\sigma_j \right)^K }.
\end{equation}  
\end{proposition}

\begin{proof}
By Theorem~\ref{thm:projective}, the map
$\Phi_{\mathcal{T}}$ is a homeomorphism of $\Delta_{K-1}$
onto itself, with inverse $\Phi_{-\mathcal{T}}$.
Moreover, the explicit formulas for $\Phi_{\mathcal{T}}$
and its inverse show that both maps are smooth in the
reduced coordinates, since their denominators are strictly
positive on the closed simplex. Hence,
$\Phi_{-\mathcal{T}(t)}$ is a $C^\infty$-diffeomorphism
in the reduced coordinates. The push-forward of an
absolutely continuous measure is therefore absolutely
continuous, and the change-of-variables formula gives
\[
    g(\widehat{\sigma},t)
    =
    g_0\!\left(\Phi_{-\mathcal{T}(t)}(\widehat{\sigma})\right)
    \left|\det D\Phi_{-\mathcal{T}(t)}(\widehat{\sigma})\right|.
\]
It remains to compute the Jacobian determinant. 
Set $D:=\sigma_K+\sum_{j=1}^{K-1}\mathrm{e}^{-\mathcal{T}_j}\sigma_j$.
The $(K-1)\times(K-1)$ matrix $D\Phi_{-\mathcal{T}}$ 
has entries
\[
    [D\Phi_{-\mathcal{T}}]_{il}
    =
    \frac{\mathrm{e}^{-\mathcal{T}_i}}{D}\delta_{il}
    -
    \frac{\mathrm{e}^{-\mathcal{T}_i}\sigma_i}{D^2}
    (\mathrm{e}^{-\mathcal{T}_l}-1),
\]
which can be written in matrix form as
\[
    D\Phi_{-\mathcal{T}}
    =
    \frac{1}{D}
    \operatorname{diag}(\mathrm{e}^{-\mathcal{T}})
    \left[
        I - \frac{1}{D}\,\widehat{\sigma}\,
        (\mathrm{e}^{-\mathcal{T}}-\mathbf{1})^\top
    \right].
\]
By the matrix determinant lemma, and using 
$\sum_j(\mathrm{e}^{-\mathcal{T}_j}-1)\sigma_j=D-1$:
\[
    \det\!\left[
        I - \frac{1}{D}\,\widehat{\sigma}\,
        (\mathrm{e}^{-\mathcal{T}}-\mathbf{1})^\top
    \right]
    =
    1 - \frac{D-1}{D}
    =
    \frac{1}{D}.
\]
Therefore
\[
    \left|\det D\Phi_{-\mathcal{T}(t)}(\widehat{\sigma})\right|
    =
    \left(\frac{1}{D}\right)^{K-1}
    \mathrm{e}^{-\sum_j\mathcal{T}_j}
    \cdot\frac{1}{D}
    =
    \frac{\mathrm{e}^{-\sum_j\mathcal{T}_j}}{D^K},
\]
and substituting into the change-of-variables formula 
yields~\eqref{eq:strategy-projective-density},  where $\sigma_K:=1-\sum_{j=1}^{K-1}\sigma_j$.
\end{proof}

\begin{remark}[Asymptotic concentration of the strategy marginal]
The above proposition ensures that absolute continuity is preserved at every finite time, but does not imply that the limiting distribution remains absolutely continuous as $t\to+\infty$. Indeed, the projective flow may become asymptotically degenerate, leading to concentration of the probability mass at the boundary of the simplex.
\end{remark}

\begin{remark}[Stationary strategy marginals and equilibrium structure]
Any finite equilibrium ${\mathcal{T}}_\infty$ of
\eqref{eq:strategy-cumulative-gains} yields a reduced mean strategy
\[
\widehat m_\sigma^\infty
=
\mathcal{E}({\mathcal{T}}_\infty)
\]
and the corresponding full mean strategy
\[
m_\sigma^\infty
=
\left(
\widehat m_\sigma^\infty,\,\,
1-\sum_{j=1}^{K-1}m_{\sigma,j}^\infty
\right)^\top
\]
satisfies the equal-payoff condition
\begin{equation}\label{eq:strategy-equal-payoff}
(Jm_\sigma^\infty)_1
=
\ldots
=
(Jm_\sigma^\infty)_K.
\end{equation}
Consequently, $m_\sigma^\infty$ is a symmetric
Nash equilibrium of the game associated with $J$, which is
interior whenever the initial strategy marginal is absolutely
continuous.
Conversely, whenever a probability distribution has a mean strategy
satisfying \eqref{eq:strategy-equal-payoff}, all relative payoff
gains vanish and its strategy transport velocity is identically
zero. Hence, every distribution with such a mean is stationary
for the strategy transport equation, and stationary strategy
marginals need not be unique. The equal-payoff condition is sufficient, but not necessary, for
stationarity of general measure-valued marginals. In particular,
the vertices of the simplex are stationary states of the
replicator flow regardless of the payoff values, and additional
stationary distributions may be supported on invariant boundary
faces. For $K\geq 3$, the equilibrium structure is therefore
richer than in the two-strategy case, and degenerate payoff
matrices may also give rise to continuous families of equilibrium
means. Convergence to a Nash equilibrium is not
guaranteed in general; however, under suitable dissipativity
conditions on $J$, the mean strategy may converge to a Nash
equilibrium, while the strategy marginal can converge to a
non-Dirac stationary distribution
\cite{HofbauerSandholm2009,LahkarSeymour2013}.
For classifications of replicator equilibria, we refer to
\cite{HofbauerSigmund1998,HofbauerSigmund2003}.
\end{remark}

\subparagraph{Exact representation of the strategy--conviction marginal.}
The projective representation derived above can be extended to the joint strategy--conviction marginal. Indeed, under constant payoffs, the strategy dynamics depends on the distribution only through the mean strategy and is independent of the conviction
variable $u$, which remains unchanged over time. We therefore introduce the strategy--conviction marginal (in reduced coordinates) of $f$ as:
\begin{equation*}
    \eta(\widehat\sigma,u,t)
    :=
    \int_I f(v,\widehat\sigma,u,t)\,\mathrm dv.
\end{equation*}

\begin{corollary}[Exact evolution of the strategy--conviction marginal]
\label{cor:eta-pushforward}
Under the assumptions of Theorem~\ref{thm:projective}, let $\eta_0(\widehat \sigma,u)$ denote the initial strategy--conviction marginal. Then:
\begin{itemize}
\item[(i)] \emph{(Measure-valued solution.)}
$\eta(\cdot,\cdot,t)$ is a probability measure on
$\Delta_{K-1}\times I$ for all $t\geq0$, and it is given
by the push-forward
\begin{equation}\label{eq:eta-pushforward}
    \eta(\cdot,\cdot,t)
    =
    \bigl(\Phi_{\mathcal{T}(t)},\mathrm{Id}\bigr)_\#\eta_0,
\end{equation}
where $\mathcal{T}(t)$ is the same cumulative gain as in
Theorem~\ref{thm:projective} and $\mathrm{Id}$ denotes the
identity map on $I$.

\item[(ii)] \emph{(Absolutely continuous case.)}
If $\eta_0$ is absolutely continuous with respect to the
Lebesgue measure on $\Delta_{K-1}\times I$, then
$\eta(\cdot,\cdot,t)$ is absolutely continuous for all
$t\geq0$, with density
\begin{equation}\label{eq:eta-density-general}
    \eta(\widehat\sigma,u,t)
    =
    \eta_0\!\bigl(\Phi_{-\mathcal{T}(t)}(\widehat\sigma),\,u\bigr)
    \cdot
    \frac{
        \exp\!\bigl(-\sum_{j=1}^{K-1}\mathcal{T}_j(t)\bigr)
    }{
        \Bigl(\sigma_K
        +\sum_{j=1}^{K-1}\mathrm{e}^{-\mathcal{T}_j(t)}\sigma_j
        \Bigr)^K
    }.
\end{equation}
\end{itemize}

In particular, the $\widehat\sigma$-marginal of $\eta$ coincides
with the strategy marginal $g(\cdot,t)$ of
Theorem~\ref{thm:projective}, while the $u$-marginal satisfies
$\int_{\Delta_{K-1}}\eta(\widehat\sigma,u,t)\,d\widehat\sigma
=k(u)$ consistently with
Proposition~\ref{prop:conviction-stationary}.
\end{corollary}

\begin{proof}
Integrating the VFP equation~\eqref{eq:FP_gain} over $v\in I$ and using the zero-flux conditions at $v=\pm1$,  $\eta$ satisfies
\begin{equation}\label{eq:eta-transport-general}
    \partial_t\eta
    +
    \frac{1}{\tau_\sigma}
    \sum_{i=1}^{K-1}
    \partial_{\sigma_i}\!
    \bigl(
    \mathcal{A}_i^\sigma[g](\widehat\sigma,t)\,
    \eta(\widehat\sigma,u,t)
    \bigr)
    =
    0,
\end{equation}
where $\mathcal{A}_i^\sigma[g](\widehat\sigma,t)
=\sigma_i(G_i(\widehat m_\sigma(t))-\sum_{j=1}^{K-1}\sigma_j
G_j(\widehat m_\sigma(t)))$ depends on $\eta$ only through
$\widehat m_\sigma(t)=\mathcal{E}(\mathcal{T}(t))$,
already determined by Theorem~\ref{thm:projective}.
In particular, this velocity field is independent of $u$,
so the characteristics of~\eqref{eq:eta-transport-general}
in $\Delta_{K-1}\times I$ are
\[
    \frac{d\widehat\sigma}{dt}
    =
    \frac{1}{\tau_\sigma}M(\widehat\sigma)
    G(\widehat m_\sigma(t)),
    \qquad
    \frac{du}{dt}=0.
\]
The characteristic flow in $\widehat\sigma$ is exactly
$\Phi_{\mathcal{T}(t)}$ by Theorem~\ref{thm:projective},
while $u$ evolves by the identity. The push-forward
representation~\eqref{eq:eta-pushforward} follows. For~(ii), the density formula~\eqref{eq:eta-density-general}
follows from the change-of-variables formula applied to the
push-forward: the Jacobian of $\Phi_{-\mathcal{T}(t)}$ with
respect to $\widehat\sigma$ is given by
Proposition~\ref{prop:density}, and the Jacobian of the
identity in $u$ is $1$.
\end{proof}

\begin{remark}[Strategy--conviction covariance and constant payoff]
Corollary~\ref{cor:eta-pushforward} implies that the
strategy--conviction covariance vector $\Sigma_{su}(t)$ \eqref{eq:cov_sigma_u} is determined exactly by $\eta_0$ and $\mathcal{T}(t)$ for all $t\geq0$, without solving the full VFP equation. 
In particular, in the case of initial independence $\eta_0(\widehat\sigma,u)=g_0(\widehat\sigma)
k(u)$, then
$\eta(\widehat\sigma,u,t)=g(\widehat\sigma,t)k(u)$
for all $t\geq0$ and $\Sigma_{\sigma u}(t)\equiv0$.
\end{remark}

\section{Asymptotic dynamics and closure regimes}\label{sec:asymptotic}

The VFP equation \eqref{eq:FP_final} couples the opinion and the strategy dynamics through two intrinsic time scales, the opinion relaxation time $\tau_v$ and the strategy adaptation  time $\tau_\sigma$. We can rewrite \eqref{eq:FP_final} as: 
\begin{equation}\label{eq:FP_timescales}
\partial_t f  = \frac{1}{\tau_v}\mathcal L_v[f] +
\frac{1}{\tau_\sigma}\mathcal L_\sigma[f],
\end{equation}
where we introduce the opinion and strategy operators:
\begin{equation*} \mathcal L_v[f] :=  -\partial_v\!\left(\mathcal A^v[f]f\right)  +  \frac{\nu^2}{2} \partial_{vv}\!\left(D(v)^2f\right), \qquad  \mathcal L_\sigma[f] := -\sum_{i=1}^K  \partial_{\sigma_i}  \left(    \mathcal A_i^\sigma[f]f \right). \end{equation*} 

The relative magnitude of $\tau_v$ and $\tau_\sigma$ determines which component of the dynamics  evolves on the fast time scale,  thereby suggesting two complementary asymptotic regimes. If $\tau_v\ll\tau_\sigma$, opinions evolve faster, while strategies and convictions provide a slowly varying background.  If $\tau_\sigma\ll\tau_v$, strategies adapt rapidly while the opinion--conviction state varies slowly. Under suitable assumptions, both regimes yield reduced kinetic descriptions of lower dimensionality by integrating out the fast variables.

As detailed in the following, these two limits provide complementary reduction  mechanisms. Following a standard zeroth-order Chapman--Enskog type approximation, the fast variable is eliminated through  a quasi-equilibrium approximation, \cite{ChapmanOnTL,Enskog}. Despite the apparent symmetry in the time-scaling, the mathematical nature of the two closures differs profoundly. The fast-opinion limit (Section \ref{sec:fast-opinion}) relies on the unique strictly positive stationary state (the Gibbs distribution) of the Fokker-Planck operator, requiring an algebraic self-consistency condition to close the macroscopic mean. In contrast, the fast-strategy limit (Section \ref{sec:fast-strategy}) involves the highly degenerate null space of the replicator operator, necessitating explicit equilibrium selection principles based on the exact projective flow (for constant payoff). In both cases, however, the resulting reduced system inherits non-local  couplings through either effective interaction payoffs or averaged compromise propensities. Numerical tests of the two asymptotic regimes are illustrated in Section~\ref{sec:numerics}.

We refer to \cite{Kuehn2015} and references therein for a comprehensive background on fast-slow dynamics.

\subsection{Fast opinion relaxation}
\label{sec:fast-opinion}
We first consider the regime $ \tau_v\ll\tau_\sigma$. Setting the scale parameter $\epsilon_v:= {\tau_v}/{\tau_\sigma}\ll1$
and introducing the slow time variable
$ \tau={t}/{\tau_\sigma}$ equation \eqref{eq:FP_timescales} rescales to
\begin{equation}\label{eq:fast-opinion-scaled}
\partial_\tau f=\frac{1}{\epsilon_v}\mathcal L_v[f]+\mathcal L_\sigma[f].
\end{equation}
We formally expand $f=f^0+\epsilon_v f^1+\mathcal O(\epsilon_v^2)$, and collect terms by order in $\epsilon_v$.

Since the opinion drift $\mathcal A^v[\cdot]$ is non-linear in $f$ through the mean $m_v$, the expansion of the opinion operator to first order in $\epsilon_v$ reads,
\begin{equation}\label{Lv_eps}
    \mathcal{L}_v[f] = \mathcal{L}_v[f^0] - \epsilon_v \partial_v \left(  \mathcal{A}^v[f^0]f^1 + \Theta(\sigma)m_v^{(1)}f^0 - \frac{\nu^2}{2} \partial_v \left( D(v)^2f^1 \right) \right) + \mathcal{O}(\epsilon_v^2),
\end{equation}
where we introduce the notation $m_v^{(k)}=\int_{\Omega}vf^{k}\mathrm{d}\omega$. In particular $m_v=m_v^{(0)}+\epsilon_v m_v^{(1)}+\mathcal O (\epsilon_v^2)$. Conversely, the strategy drift $\mathcal{A}_i^{\sigma}[\cdot]$ depends linearly on $f$ through the payoff integral, therefore the strategy operator expansion reads
\begin{equation}\label{Lsigma_eps}
 \mathcal{L}_\sigma[f] =  \mathcal{L}_\sigma[f^0] -\epsilon_v \sum_{i=1}^K \partial_{\sigma_i} \left( \mathcal{A}_i^{\sigma}[f^0] f^1 + \mathcal{A}_i^{\sigma}[f^1] f^0  \right) + \mathcal{O}(\epsilon_v^2).
 \end{equation}

\paragraph{Leading order $O(\varepsilon^{-1})$:  conditional  quasi-equilibrium.}
At leading order, $f^0$,  satisfies the nonlinear stationary equation, \begin{equation}\label{eq:lead_order}\mathcal L_v[f^0]=0,\end{equation} and therefore belongs to the equilibrium manifold of the fast opinion dynamics. Equation \eqref{eq:lead_order} is a stationary Fokker--Planck equation in $v$ with the variables$(\sigma,u,\tau)$ acting as frozen parameters. We define the strategy--conviction marginal of $f^0$,
\begin{equation*}
\eta^0(\sigma,u,\tau)
:=\int_I f^0(v,\sigma,u,\tau)\,\mathrm dv,
\end{equation*}
and, for $\eta(\sigma,u,\tau)>0$,  the conditional quasi-equilibrium  density,
\[
\rho^{\rm eq}(v\mid \sigma,u;\tau)
:=
\frac{f^0(v,\sigma,u,\tau)}
{\eta^0(\sigma,u,\tau)}, \qquad
    \int_I\rho^{\rm eq}(v\mid \sigma,u;\tau)\,\mathrm dv=1,
\] so that $f^0=\rho^{\rm eq}\eta^0$ by construction. 
Substituting this decomposition into \eqref{eq:lead_order} and
dividing by $\eta^0>0$ yields, for each class $(\sigma,u)$, the zero-flux stationary equation in $v$ supplemented by the zero-flux boundary condition at $v=\pm 1$,
\begin{equation}\label{eq:stationary-opinion-FP}
    \mathcal{A}^v[f^0](v,\sigma,u;\tau)\,\rho^{\rm eq}
    -
    \frac{\nu^2}{2}\partial_v\bigl(D(v)^2\rho^{\rm eq}\bigr)
    =
   0.
\end{equation}

Defining the local opinion target,
\begin{equation}\label{eq:local-target}
    \hat V[\eta^0](\sigma,u;\tau)
    :=
    \Theta(\sigma)\,m_v^{\rm eq}[\eta^0](\tau)
    +(1-\Theta(\sigma))\,u,
\end{equation}
as a convex combination of the global quasi-equilibrium mean, 
\[
m_v^{\rm eq}[\eta^0](\tau)
:=
m_v^{(0)}=\iint_{\Delta_{K-1}\times I}
    \left[
        \int_I v\,\rho^{\rm eq}(v\mid \sigma,u;\tau)\,\mathrm dv
    \right]\eta^0(\sigma,u,\tau)\,\mathrm{d}\sigma\,\mathrm{d}u,
\]
and the private conviction, $u$, the unique normalized solution
of~\eqref{eq:stationary-opinion-FP} is the Gibbs density,
\begin{equation}\label{eq:rho-eq}
    \rho^{\rm eq}(v\mid \sigma,u;\tau)
    =
    \frac{C(\sigma,u,\tau)}{D(v)^2}
    \exp\!\left(
        \frac{2}{\nu^2}
        \int^v
        \frac{\hat V[\eta^0](\sigma,u;\tau)-z}{D(z)^2}\,\mathrm{d}z
    \right),
\end{equation}
where $C(\sigma,u,\tau)>0$ is the normalization constant, which is finite for any $D(v)$ satisfying condition (D1) of Proposition~\ref{prop:domain_invariance}. Note that $\hat V[\eta^0]$ depends implicitly on $m_v^{\rm eq}[\eta^0](\tau)$, which is itself unknown and must be determined via a self-consistency condition.
 
The marginal $\eta^0$ is unaffected by the fast opinion operator, which conserves mass for each fixed pair $(\sigma,u)$ under vanishing opinion flux at the boundary of $I$. It therefore evolves on the slow strategy time scale. Its evolution equation is obtained by integrating the $\mathcal O(1)$ equation of the asymptotic expansion with respect to $v$.

\paragraph{Self-consistency for the mean opinion.}

Integrating the zero-flux equation \eqref{eq:stationary-opinion-FP} over $I$ and using $D(\pm1)=0$ (condition (D1) of Proposition~\ref{prop:domain_invariance}) to eliminate boundary terms gives:
 \[
\int_I v\,\rho^{\rm eq}(v\mid \sigma,u;\tau)\,\mathrm{d}v
=
\hat V[\eta^0](\sigma,u;\tau).
\]
We consequently obtain the following self-consistency
condition
\begin{align}
    m_v^{\rm eq}[\eta^0](\tau)
    &=
    \iint_{\Delta_{K-1}\times I}
    \hat V[\eta^0](\sigma,u;\tau)
    \eta^0(\sigma,u,\tau)\,\mathrm{d}\sigma\,du\cr
   & = m_v^{\rm eq}[\eta^0](\tau)\theta^\top m_\sigma^{(0)}[\eta^0](\tau) + (1-\theta^\top m_\sigma^{(0)}[\eta^0](\tau))m_u^{(0)} - \theta^\top\Sigma_{\sigma u}[\eta^0](\tau),
\label{eq:self-consistency-mv}
\end{align}
where we have defined the moments:
\[ m^{(0)}_\sigma[\eta^0](\tau)
    :=
    \iint_{\Delta_{K-1}\times I}
    \sigma\,\eta^0(\sigma,u,\tau)\,\mathrm{d}\sigma\,du,
    \qquad
    m^{(0)}_u[\eta^0](\tau)
    :=
    \iint_{\Delta_{K-1}\times I}
    u\,\eta^0(\sigma,u,\tau)\,\mathrm{d}\sigma\,du,
\] and the strategy--conviction covariance vector:
\[    
\Sigma^{(0)}_{\sigma u}[\eta^0](\tau) := \int_{\Delta_{K-1}\times I}
    (\sigma -m^{(0)}_\sigma(\tau) )(u-m_u^{(0)})
    \eta^0(\sigma,u,\tau)\,\mathrm{d}\sigma\,du
\] 
As will be shown below, the slow evolution equation for $\eta^0$ is a conservation law in the strategy variable $\sigma$, with no flux across the boundary of $\Delta_{K-1}$. Consequently, the marginal distribution of the conviction variable $u$ is preserved, and in particular $m^{(0)}_u[\eta^0](\tau)=m^{(0)}_u[\eta^0](0)$ is constant in $\tau$.

Rearranging~\eqref{eq:self-consistency-mv} as a linear equation in $m_v^{\rm eq}[\eta^0]$, and provided that $\theta^\top m_\sigma^{(0)}[\eta^0](\tau)<1$, the quasi-equilibrium mean opinion is uniquely determined by the self-consistency relation
\begin{equation}\label{eq:fast-opinion-mean}
    m_v^{\rm eq}[\eta^0](\tau)
    =
    m_u^{(0)}[\eta^0]
    -
    \frac{
    \theta^\top\Sigma^{(0)}_{\sigma u}[\eta^0](\tau)
    }{
        1-\theta^\top m^{(0)}_\sigma[\eta^0](\tau)
    },
\end{equation}
which represents the quasi-stationary counterpart of the exact macroscopic  equation \eqref{eq:mean_opinion}.

Two observations follow. First, if $\Sigma_{\sigma u}[\eta^0]=0$ (i.e.~strategy and conviction are uncorrelated under $\eta^0$), then $m_v^{\rm eq}[\eta^0]=m_u^{(0)}[\eta^0]$ regardless of the strategic mean $m_\sigma^{(0)}[\eta^0]$. Second, if $\theta^\top m_\sigma^{(0)}[\eta^0](\tau)=1$, the equation is completely degenerate: this implies $\Theta(\sigma)=1$ almost everywhere with respect to $\eta^0$, meaning the opinion dynamics has no restoring force toward individual convictions, and $m_v^{\rm eq}$ is conserved from the initial data.

\begin{remark} (Beta quasi-equilibrium)\label{rmk:beta_quasi_eq}
In the uncorrelated, non-degenerate case discussed above, the local opinion target reduces to $$\widehat V[\eta^0](\sigma,u;\tau)
=
\Theta(\sigma)m_u^{(0)}[\eta^0]
+
(1-\Theta(\sigma))u.$$
For the specific diffusion choice  $D(v)=\sqrt{1-v^2}$, the conditional quasi-equilibrium is a rescaled Beta distribution on $I$, explicitly given by
\[
\rho^{\rm eq}(v\mid \sigma,u;\tau)
=
C(\sigma,u;\tau)
(1+v)^{\frac{1+\widehat V[\eta^0](\sigma,u;\tau)}{\nu^2}-1}
(1-v)^{\frac{1-\widehat V[\eta^0](\sigma,u;\tau)}{\nu^2}-1}.
\]
For $|\widehat V[\eta^0](\sigma,u;\tau)|<1$, its mean and variance are, respectively:
\[
\int_I v \rho^{\rm eq}(v\mid \sigma,u;\tau)\ \mathrm{d}v=\widehat V[\eta^0](\sigma,u;\tau),
\,
\int_I (v-\widehat V[\eta^0](\sigma,u))^2 \rho^{\rm eq}(v\mid \sigma,u;\tau)\ \mathrm{d}v=
\frac{\nu^2}{2+\nu^2}
\left(1-\widehat V[\eta^0](\sigma,u;\tau)^2\right).
\]
Thus, in this scenario, the conditional quasi-equilibrium is explicitly parametrized by the state pair $(\sigma,u)$, the compromise propensity $\Theta(\sigma)$, and the conserved population mean conviction $m_u^{(0)}[\eta^0]$.
\end{remark}

\paragraph{Next order $\mathcal O (1)$: slow strategy dynamics.}
On the slow time scale, the strategy--conviction marginal $\eta^0$ evolves in the quasi-equilibrium opinion background $\rho^{\rm eq}[\eta^0]$. Collecting the $\mathcal{O}(1)$ terms in~\eqref{eq:fast-opinion-scaled} yields:
\begin{equation*}
\partial_\tau f^0 =   -  \partial_v \left(  \mathcal{A}^v[f^0]f^1 + \Theta(\sigma)m_v[f^1]f^0 - \frac{\nu^2}{2} \partial_v \left( D(v)^2 f^1 \right) \right) 
+\mathcal{L}_\sigma[f^0] .
\end{equation*}
To derive the evolution equation for $\eta^0$, we integrate over $v\in I$. The opinion divergence vanishes due to the zero-flux boundary conditions, obtaining:
\begin{equation}\label{eta_int}
\partial_\tau\eta^0
=-
\sum_{i=1}^K
\partial_{\sigma_i}
\left( \eta^0\int_I \mathcal A^\sigma_i[f^0]\rho^{\rm eq}[\eta^0]\,\mathrm{d}v
\right).
\end{equation}

By averaging the original interaction payoff over the conditional quasi-equilibrium opinion distributions of both interacting agents, we define the effective expected payoff,
\begin{equation*}
\mathcal P_i^{\rm eq}[\eta^0](\sigma,u,\tau)
=
\sum_{j=1}^K
\iint_{\Delta_{K-1}\times I}
\sigma_{*,j}\,
\mathcal J_{ij}^{\rm eq}[\eta^0]
(\sigma,u;\sigma_*,u_*,\tau)
\eta^0(\sigma_*,u_*,\tau)\,
\mathrm{d}\sigma_*\,\mathrm{d}u_*,
\end{equation*}
where the effective payoff kernel is given by
\begin{equation*}
\mathcal J_{ij}^{\rm eq}[\eta^0]
(\sigma,u;\sigma_*,u_*,\tau)
:=
\iint_{I^2}
J_{ij}(v,u;v_*,u_*)
\rho^{\rm eq}[\eta^0](v\mid \sigma,u;\tau)
\rho^{\rm eq}[\eta^0](v_*\mid \sigma_*,u_*;\tau)
\,\mathrm{d}v\,\mathrm{d}v_*.
\end{equation*}
Introducing the effective average payoff 
\[
\overline{\mathcal P}^{\rm eq}[\eta^0](\sigma,u,\tau)
=
\sum_{i=1}^K
\sigma_i\mathcal P_i^{\rm eq}[\eta^0](\sigma,u,\tau),
\]
equation \eqref{eta_int} rewrites concisely as
\begin{equation*}
\partial_\tau\eta^0
+
\sum_{i=1}^K
\partial_{\sigma_i}
\left[
\sigma_i
\left(
\mathcal P_i^{\rm eq}[\eta^0]
-
\overline{\mathcal P}^{\rm eq}[\eta^0]
\right)
\eta^0
\right]
=0.
\end{equation*}

\paragraph{The reduced kinetic system and its closure.} Collecting the preceding relations, the leading-order dynamics in the non-degenerate regime $\theta^\top m_\sigma^{(0)}[\eta^0](\tau)<1$ is governed by the following closed system:
\begin{equation*}
\begin{cases}
\displaystyle
m_v^{\rm eq}[ \eta^0](\tau)
=
   m_u^{(0)}[\eta^0]
    -
    \frac{
        \theta^\top\Sigma^{(0)}_{\sigma u}[ \eta^0](\tau)
    }{
        1-\theta^\top m_\sigma^{(0)}[ \eta^0](\tau)
    },
\\[5mm]
\displaystyle
\rho^{\rm eq}[ \eta^0](v\mid \sigma,u,\tau)
=
\frac{C(\sigma,u,\tau)}{D(v)^2}
\exp\left[
\frac{2}{\nu^2}
\int^v
\frac{
\Theta(\sigma)m_v^{\rm eq}[ \eta^0](\tau)
+
(1-\Theta(\sigma))u-z
}{
D(z)^2
}\,\mathrm{d}z
\right],
\\[5mm]
\displaystyle
\partial_\tau \eta^0
+
\sum_{i=1}^K
\partial_{\sigma_i}
\left[
\sigma_i
\left(
\mathcal P_i^{\rm eq}[ \eta^0]
-
\sum_{j=1}^K\sigma_j
\mathcal P_j^{\rm eq}[ \eta^0]
\right)
 \eta^0
\right]
=0.
\end{cases}
\end{equation*}

These three equations play distinct interconnected roles: the first is an algebraic self-consistency relation determining $m_v^{\rm eq}[ \eta^0]$; the second is the explicit Gibbs formula for $\rho^{\rm eq}[ \eta^0]$; the third is the PDE governing the slow evolution of $ \eta^0$. All auxiliary quantities, $m_\sigma[ \eta^0]$, $m_u[ \eta^0]$, $\Sigma_{\sigma u}[ \eta^0]$, are moments of $ \eta^0$, computable by direct integration at each $\tau$ without additional evolution equations. The dependence chain
\[
     \eta^0
    \;\to\;\left(
    m_\sigma[ \eta^0],\,\Sigma_{\sigma u}[ \eta^0],\,m_u[ \eta^0]\right)
    \;\to\;
    m_v^{\rm eq}[ \eta^0]
    \;\to\;
    (\hat V[ \eta^0],\,\rho^{\rm eq}[ \eta^0])
    \;\to\;
    \mathcal{J}_{ij}^{\rm eq}[ \eta^0]
    \;\to\;
    \mathcal{P}_i^{\rm eq}[ \eta^0]
\]
is sequential and non-recursive given $ \eta^0$, confirming that the system is self-contained and closed at the kinetic level within the reduced state space $\Delta_{K-1}\times I$. The system does not, however, automatically yield a finite-dimensional closure of the moment equations for $m_\sigma[ \eta^0]$ or $\Sigma_{\sigma u}[ \eta^0]$. Rather, these evolve through the full PDE for $ \eta^0$, and no moment hierarchy truncation is universally available for general state-dependent payoff structures.

\subsection{Fast strategy dynamics}
\label{sec:fast-strategy}
We now consider the complementary regime $ \tau_\sigma\ll\tau_v$. Setting the scale parameter $\epsilon_\sigma:= {\tau_\sigma}/{\tau_v}\ll1$ and defining the slow time as $ \tau={t}/{\tau_v}$ we rescale \eqref{eq:FP_timescales} as follows
\begin{equation}\label{eq:fast-strategy-scaled}
    \partial_\tau f
    =
   \mathcal L_v[f]
    +
     \frac{1}{\epsilon_\sigma}\mathcal L_\sigma[f].
\end{equation}
We formally expand $f=f^0+\epsilon_\sigma f^1+\mathcal O(\epsilon_\sigma^2)$, and collect terms by order in $\epsilon_\sigma$.  The expansions of the operators $\mathcal L_v[f]$ and $\mathcal L_\sigma[f]$ are derived analogously to the fast-opinion regime, with $\epsilon_v$ replaced by $\epsilon_\sigma$. For brevity, we omit the details and refer back to \eqref{Lv_eps} and \eqref{Lsigma_eps} for the corresponding first-order expansions.

\paragraph{Leading order $\mathcal O(\epsilon_\sigma^{-1})$: characterization of $f^0$.} The leading-order distribution equation is a stationary transport equation in $\sigma$ with $(v,u,\tau)$ frozen, given by
\begin{equation}\label{eq:leading-order-strategy}
\mathcal L_\sigma[f^0]=0,\end{equation} 
thus $f^0$ belongs to the null space of $\mathcal{L}_\sigma$ at each slow time $\tau$.

Since $\mathcal L_\sigma$ acts exclusively on the strategy variable and the replicator flux vanishes on $\partial \Delta_{K-1}$, equation \eqref{eq:leading-order-strategy} conserves mass for each fixed $(v,u)$, and the opinion--conviction marginal of $f^0$, defined by
\begin{equation*}
\rho^0(v,u,\tau) := \int_{\Delta_{K-1}} f^0(v,\sigma,u,\tau)\,\mathrm{d}\sigma,
\end{equation*}
is therefore invariant under the fast strategy relaxation and evolves entirely on the slow opinion time scale.

For $\rho^0(v,u,\tau)>0$, we define the conditional  distribution
\begin{equation*}
    \eta^{\rm qs}(\sigma\mid v,u;\tau)=\frac{f^0(v,\sigma,u,\tau)}{\rho^0(v,u,\tau)}\,,
    \qquad
    \int_{\Delta_{K-1}} \eta^{\rm qs}(\sigma\mid v,u;\tau)\,\mathrm{d}\sigma=1,
\end{equation*}
so that $f^0=\eta^{\rm qs}\rho^0$ by construction.  The superscript ``qs'' stands for \emph{quasi-stationary}: $\eta^{\rm qs}$ is stationary on the fast time scale but may in principle evolve on the slow time scale $\tau$.

The null space of $\mathcal{L}_\sigma$ is infinite-dimensional; any distribution supported exclusively on replicator equilibria is stationary under the leading-order strategy dynamics. By contrast, a distribution transported along nonstationary replicator orbits generally evolves in time and does not belong to the null space. Therefore, $\eta^{\rm qs}$ is not uniquely determined by \eqref{eq:leading-order-strategy} alone. This differs profoundly from the fast-opinion regime, where, for a given global opinion mean the stationary equation explicitly identifies the conditional equilibrium density in $v$. 

To uniquely identify $\eta^{\rm qs}$, one must analyze the dynamics on the fast time scale and impose a matching condition. This requires specific assumptions on the payoff structure, or on the equilibrium selection, as discussed below.

\paragraph{Identification of $\eta^{\rm qs}$: fast time scale and constant payoff.}
On the fast time scale $\zeta:=t/\tau_\sigma$, with $(v,u)$ treated as frozen parameters, the leading dynamics is governed by the replicator equation $\partial_\zeta f\approx\mathcal{L}_\sigma[f]$. The fast dynamics thus transports the initial strategy distribution toward an element of the null space of $\mathcal{L}_\sigma$.

For general state-dependent payoffs, the relative payoff gains depend on the frozen state $(v,u)$: each opinion--conviction class follows a distinct replicator trajectory, and the limiting distribution (if it exists) depends structurally on $(v,u)$. In such cases, $\eta^{\rm qs}$ cannot be reconstructed from $\rho^0$ alone, and a separate equilibrium selection principle (such as selection of an evolutionarily stable state for each frozen $(v,u)$) is required. We do not pursue this direction here.

Instead, we focus on the case of constant payoffs ($J \in \mathbb{R}^{K \times K}$), for which the analysis simplifies substantially. We introduce the following assumptions.

\noindent\textbf{(A1) Initial independence.} The strategy is initially independent of the opinion--conviction pair, i.e.~$f(v,\sigma,u,t=0)=\rho(v,u,t=0)\,g(\sigma,t=0)$.

Under (A1), the conditional initial
strategy distribution is the same across all classes: $\eta(\sigma,\zeta=0\mid v,u)=g(\sigma,\zeta=0)$. By Theorem~\ref{thm:projective},  the projective flow acts identically on every class:
\begin{equation*}
    \eta(\sigma,\zeta\mid v,u)
    =
    \bigl(\Phi_{\mathcal{T}(\zeta)}\bigr)_\#g(\sigma,\zeta=0),
\end{equation*} with cumulative gains $\mathcal{T}(\zeta)$ independent of $(v,u)$ and determined by~\eqref{eq:strategy-cumulative-gains}.

\noindent\textbf{(A2) Convergence of the fast dynamics.}
We assume that all relative cumulative gains $\mathcal{T}_i(\zeta)-\mathcal{T}_j(\zeta)$, for $i,j\in\{1,\ldots,K\}$, converge in $[-\infty,+\infty]$ as $\zeta\to+\infty$, where $\mathcal{T}_K:=0$ denotes the cumulative gain of the reference strategy. This condition ensures  pointwise convergence of the projective flow
\[
    \Phi_\infty(\sigma)
    :=
    \lim_{\zeta\to+\infty}
    \Phi_{\mathcal{T}(\zeta)}(\sigma),
    \qquad \sigma\in\Delta_{K-1}.
\]
Consequently, the transported strategy distribution converges weakly as a probability measure on $\Delta_{K-1}$,  and
\begin{equation*}
    g^{\rm eq}(\sigma)
    :=
    \lim_{\zeta\to+\infty}
    \bigl(\Phi_{\mathcal{T}(\zeta)}\bigr)_\#g(\sigma,\zeta=0).
\end{equation*}
The condition on the relative cumulative gains is sufficient because the projective flow depends on $\mathcal{T}$ only through pairwise differences of its components. 

Under (A2), the fast dynamics selects the same limiting distribution for every opinion--conviction class, yielding the limiting mean:
\begin{equation*}
    m_\sigma[g^{\rm eq}]:=\int_{\Delta_{K-1}}\sigma g^{\rm eq}(\sigma)\mathrm{d}\sigma
    =\lim_{\zeta\to+\infty}\mathcal{E}(\mathcal{T}(\zeta)),
\end{equation*} 
If $\mathcal{T}(\zeta)\to\mathcal{T}_\infty \in\mathbb{R}^{K-1}$, then $m_\sigma[g^{\rm eq}]=\mathcal{E}(\mathcal{T}_\infty)$ and the limiting mean satisfies the equal-payoff condition \eqref{eq:strategy-equal-payoff}. If some relative cumulative gains diverge, the limiting distribution may instead be supported on a proper face, possibly a vertex, of $\Delta_{K-1}$, determined by the asymptotic relative gains and the support of the initial distribution, corresponding to the dominant strategy or strategies.

The matching condition between the fast and slow time scales then reads:
\begin{equation}\label{eq:matching}
\eta^{\rm qs}(\sigma\mid v,u;\tau=0)
 =
 \lim_{\zeta\to+\infty}
\eta(\sigma,\zeta\mid v,u)
= g^{\rm eq}(\sigma).
\end{equation}

\begin{remark}[Non-convergent fast dynamics]
Assumption~(A2) fails in two critical scenarios, such as conservative oscillations and heteroclinic cycles, in which the strategy distribution $g(\sigma,\zeta)$ does not converge. In these cases, the effective compromise propensity may oscillate on the fast time scale. The response of the opinion marginal $h(v,t)$ depends on the coupling between strategy and opinion dynamics. The fast-strategy closure developed below does not apply to these non-convergent regimes.
\end{remark}

\paragraph{Factorization closure and its validity.}

The matching condition~\eqref{eq:matching} identifies the post-initial-layer strategy profile at $\tau=0$. To obtain a closed description over the entire slow time scale, we introduce a third assumption.

\noindent\textbf{(A3) Factorization closure.} The quasi-stationary conditional strategy distribution remains permanently equal to $g^{\rm eq}$ for all $\tau\geq0$:
\begin{equation*}
\eta^{\rm qs}(\sigma\mid v,u;\tau)
\approx
g^{\rm eq}(\sigma),
\end{equation*}
so that the leading-order distribution is approximated by
\begin{equation*}
f^0(v,\sigma,u,\tau)
\approx
g^{\rm eq}(\sigma)\,\rho^0(v,u,\tau).
\end{equation*}
Assumption (A3) amounts to neglecting correlations between the strategy variable and the opinion--conviction pair $(v,u)$ on the slow time scale. Its validity depends on the dependence of the opinion dynamics on $\sigma$. In particular, if $\Theta(\sigma)$ is constant on the support of $g^{\rm eq}$ (e.g.~ for pure-strategy equilibria), the opinion operator is independent of $\sigma$ on that support and no strategy--opinion correlations are generated, provided the factorization holds at $\tau=0$. In this case, (A3) is preserved by the reduced dynamics. Conversely, when $\Theta$ varies across the equilibrium strategy profile, the slow opinion dynamics generally generates correlations between $\sigma$ and $(v,u)$, and (A3) should be regarded as an averaged decorrelation closure.

We stress that the consistency of (A3) with the leading-order fast-strategy equation \eqref{eq:leading-order-strategy}  is a distinct issue: whenever $g^{\rm eq}$ is supported on replicator equilibria, the leading-order strategy flux vanishes for the factorized ansatz. This condition ensures consistency with the fast strategy dynamics, but does not by itself guarantee that the factorization is preserved by the slow opinion dynamics.

Under assumption (A1)-(A3),  we therefore consider the following factorized approximation for the leading-order
distribution:
\[
    f^0(v,\sigma,u,\tau)
    :=
    g^{\rm eq}(\sigma)\,\rho^0(v,u,\tau),
\] which is now fully characterized. Within this reduced description, the strategy distribution is {frozen} at $g^{\rm eq}$ while the remaining slow evolution takes is described by $\rho^0(v,u,\tau)$.

\paragraph{Order $\mathcal O (1)$: slow opinion dynamics and the reduced equation.} Collecting the $\mathcal O (1)$ terms in \eqref{eq:fast-strategy-scaled} yields:

\begin{equation*}
    \partial_\tau f^0
    =
    \mathcal{L}_v[f^0]
    -  \sum_{i=1}^K \partial_{\sigma_i} \left( \mathcal{A}_i^{\sigma}[f^0] f^1 + \mathcal{A}_i^{\sigma}[f^1] f^0  \right) .
\end{equation*}

Integrating over $\Delta_{K-1}$, the strategy-divergence term vanishes by the divergence theorem and the zero-flux condition on $\partial\Delta_{K-1}$, while the averaged opinion drift on $f^0$ integrated against $g^{\rm eq}$ reads
\begin{align*}
\int_{\Delta_{K-1}}
    \mathcal{A}^v[f^0]\,g^{\rm eq}(\sigma)\,\mathrm{d}\sigma
&= \overline\Theta^{\rm eq}(m_v[\rho^0]-v)+(1-\overline\Theta^{\rm eq})(u-v),\\ &=
\overline\Theta^{\rm eq}(m_v[\rho^0]-u)
+u-v,
\end{align*}
where \begin{equation*}
    m_v^{(0)}[\rho^0](\tau):=\int_{I^2} v\rho^0(v,u,\tau)\,\mathrm{d}v\,\mathrm{d}u,
\end{equation*}
and the  averaged compromise propensity is a  constant scalar:
\begin{equation*}
    \overline\Theta^{\rm eq}
    :=
    \int_{\Delta_{K-1}}
    \Theta(\sigma)g^{\rm eq}(\sigma)\,\mathrm{d}\sigma
    =
    \theta^\top m_\sigma[g^{\rm eq}]
    =
    \theta^\top\mathcal{E}(\mathcal{T}_\infty)
    \in[0,1].
\end{equation*}
Since $\Theta(\sigma)$ is linear,  $\overline\Theta^{\rm eq}=\Theta(m_\sigma[g^{\rm eq}])$.
Thus, within the factorized approximation, the strategy distribution enters the reduced opinion dynamics only through its averaged equilibrium compromise propensity $\overline\Theta^{\rm eq}$.

We therefore obtain the following closed Fokker--Planck equation for the reduced density $\rho^0(v,u,\tau)$,

\begin{equation}\label{eq:slow-opinion-fast-strategy}
    \partial_\tau\rho^0
    =
    -\partial_v
    \left(
        \left(
        \overline\Theta^{\rm eq}\left(m_v^{(0)}[\rho^0]-u\right)
        +u-v
        \right)\rho^0
    \right)    +
    \frac{\nu^2}{2}
\partial_{vv}\left(D(v)^2\rho^0\right).
\end{equation}
The original VFP equation in three variables $(v,\sigma,u)$ has been effectively reduced to  two variables $(v,u)$, with the entire strategic complexity encoded into the single parameter $\overline\Theta^{\rm eq}$.

Note that, since since the operator in~\eqref{eq:slow-opinion-fast-strategy}  contains only $v$-derivatives, integrating over $v\in I$ gives $\partial_\tau k^0(u,\tau)=0$, and hence the conviction marginal $k^0(u):=\int_I\rho^0(v,u,\tau)\,\mathrm{d}v$ is conserved. In particular, $m_u^{(0)}:=\int_{I^2}u\,\rho^0\,\mathrm{d}v\,du$ is constant in $\tau$.

The first moment equation is obtained integrating \eqref{eq:slow-opinion-fast-strategy} against $v$ resulting in
\begin{equation}\label{eq:mean-opinion-fast-strategy}
    \frac{d}{d\tau} m_v^{(0)}[\rho^0](\tau)
    =
    \left(1-\overline\Theta^{\rm eq}\right)\left(m_u^{(0)}-m_v^{(0)}[\rho^0](\tau)\right).
\end{equation}
Provided $\overline\Theta^{\rm eq}<1$, its solution is
\begin{equation}\label{sol_mv}m_v^{(0)}[\rho^0](\tau)=m_u^{(0)}+(m_v^{(0)}[\rho^0](0)-m_u^{(0)}) {\rm e}^{-(1-\overline\Theta^{\rm eq})\tau},\end{equation}
which converges exponentially to $ m_u^{(0)}$, with rate $(1-\overline\Theta^{\rm eq})/\tau_v$ in the original time variable.  In the degenerate case $\overline\Theta^{\rm eq}=1$, the mean opinion is conserved, $ m_v^{(0)}[\rho^0](\tau) \equiv m_v^{(0)}[\rho^0](0).$

Assuming convergence to a zero-flux stationary state as  $\tau\to+\infty$, the stationary joint opinion--conviction density factorizes as $\rho^0_\infty(v,u)= \rho^0_\infty(v\mid u)k^0(u)$, where the conditional density, is given by the Gibbs form:
\begin{equation*}
    \rho^0_\infty(v\mid u)
    =
    \frac{C(u)}{D(v)^2}
    \exp\left[
        \frac{2}{\nu^2}
        \int^v
        \frac{
            V_{\infty}(u)-z
        }{D(z)^2}\,\mathrm{d}z
    \right],
\end{equation*}
where $C(u)$ is determined by normalization, and the conditional attraction center is defined as:
$$
V_{\infty}(u)=\begin{cases}
\overline\Theta^{\rm eq}m_u^{(0)}
            +(1-\overline\Theta^{\rm eq})u, & \overline\Theta^{\rm eq}<1\\
m_v^{(0)}[\rho^0](0), & \overline\Theta^{\rm eq}=1.
\end{cases}
$$
For \(\overline\Theta^{\rm eq}<1\), the attraction center interpolates smoothly between individual stubbornness (\(\overline\Theta^{\rm eq}=0\), center \(u\)) and population consensus (\(\overline\Theta^{\rm eq}\to1^-\), center \(m_u^{(0)}\), independent of \(u\)). At the endpoint \(\overline\Theta^{\rm eq}=1\), instead, the opinion mean is conserved, and the attraction center is determined by its initial value.

Thus, in the fast-strategy regime, the projective dynamics first selects the stationary profile $g^{\rm eq}$, while the subsequent reduced opinion dynamics relaxes towards the equilibrium associated with its averaged compromise propensity $\overline\Theta^{\rm eq}$.

\paragraph{Moment equations and closure.}
Within the factorized fast-strategy model, we can derive a reduced description of the macroscopic opinion dynamics. We introduce the opinion variance, the conviction variance and their covariance:

\[  \Sigma_v^{(0)}(\tau)
    :=
    \int_{I^2}\left(v-m_v^{(0)}[\rho^0]\right)^2\rho^0(v,u,\tau)\,\mathrm{d}v\,\mathrm{d}u,\qquad
    \Sigma_u^{(0)}
    :=
    \int_I \left(u-m_u^{(0)}\right)^2k^0(u)\,\mathrm{d}u,
   \]
\[
    \Sigma_{vu}^{(0)}(\tau)
    :=
    \int_{I^2}
    \left(v-m_v^{(0)}[\rho^0]\right)\left(u-m_u^{(0)}\right)\rho^0(v,u,\tau)\,\mathrm{d}v\,\mathrm{d}u.
\]
By conservation of the conviction marginal, $\Sigma_u^{(0)}$ is constant in time.  The equation for $m_v^{(0)}[\rho^0]$ is given by \eqref{eq:mean-opinion-fast-strategy}, whereas the remaining second-order moments of \eqref{eq:slow-opinion-fast-strategy} satisfy
\begin{equation}\label{eq:fast-strategy-moments-general}
\begin{aligned}
    \frac{d }{d\tau}\Sigma_{vu}^{(0)}
    &=
    -\Sigma_{vu}^{(0)}
    +\left(1-\overline\Theta^{\rm eq}\right)\Sigma_u^{(0)},
    \\
    \frac{d }{d\tau}\Sigma_v^{(0)}
    &=
    -2\Sigma_v^{(0)}
    +2(1-\overline\Theta^{\rm eq})\Sigma_{vu}^{(0)}
    +\nu^2
\int_{I^2}D(v)^2\rho^0(v,u,\tau)\,\mathrm{d}v\,\mathrm{d}u .
\end{aligned}
\end{equation}
Thus, although the fast-strategy ansatz removes the moments and correlations involving the strategy variable, the second-order macroscopic system is not closed for a general diffusion coefficient $D$, because the last term in \eqref{eq:fast-strategy-moments-general} is not necessarily determined by the first and second moments of $\rho^0$.

For the particular diffusion choice $ D(v)=\sqrt{1-v^2}$, we have
$$
\int_{I^2}D(v)^2\rho^0(v,u,\tau)\,\mathrm{d}v\,\mathrm{d}u = 1-\int_{I^2}v^2\rho^0(v,u,\tau)\,\mathrm{d}v\,\mathrm{d}u =1-\Sigma_v^{(0)}(\tau)-\left(m_v^{(0)}[\rho^0](\tau)\right)^2,
$$ therefore, this algebraic relation elegantly closes the hierarchy and \eqref{eq:mean-opinion-fast-strategy}-\eqref{eq:fast-strategy-moments-general} reduces to a closed system. Hence, in this particular case the reduced Fokker--Planck equation closes exactly at second order, admitting a closed finite-dimensional macroscopic description without any moment closure approximation beyond the three hypotheses above, bypassing entirely the moment hierarchy generated by the full joint distribution $f$.

Since the equation for $\Sigma_{vu}^{(0)}$ has explicit solution
\begin{equation}\label{sol_svu}
\Sigma_{vu}^{(0)}(\tau)
=\left(1-\overline\Theta^{\rm eq}\right)\Sigma_u^{(0)}
+\left(\Sigma_{vu}^{(0)}(0)-\left(1-\overline\Theta^{\rm eq}\right)\Sigma_u^{(0)}\right)\mathrm{e}^{-\tau}, 
\end{equation} the equation for $\Sigma_v^{(0)}$ reduces to the scalar ODE:
\begin{equation*}
\begin{aligned}
\frac{\mathrm{d}\Sigma_v^{(0)}}{d\tau}
={}&-(2+\nu^2)\Sigma_v^{(0)}
\\
&+2\left(1-\overline\Theta^{\rm eq}\right)\left[
\left(1-\overline\Theta^{\rm eq}\right)\Sigma_u^{(0)}
+\left(\Sigma_{vu}^{(0)}(0)-\left(1-\overline\Theta^{\rm eq}\right)\Sigma_u^{(0)}\right)\mathrm{e}^{-\tau}
\right]
\\
&+\nu^2\left[
1-
\left(
m_u^{(0)}
+(m_v^{(0)}(0)-m_u^{(0)})\mathrm{e}^{-\left(1-\overline\Theta^{\rm eq}\right)\tau}
\right)^2
\right].
\end{aligned}   
\end{equation*}

In this specific setting, the  stationary moments of the factorized reduced model can be characterized explicitly.

\begin{proposition}[Stationary moments under fast-strategy closure] \label{prop:fast-strategy-stationary}
For $D(v)=\sqrt{1-v^2}$ and under assumptions (A1)--(A3), the closed moment system~\eqref{eq:mean-opinion-fast-strategy}-\eqref{eq:fast-strategy-moments-general} admits  the following unique stationary moment limits, for $\tau\to+\infty$.

\noindent{If $\overline\Theta^{\rm eq}\in[0,1)$:}
\begin{equation*}
    m_v^{(0)}[\rho^0](\tau)\to m_u^{(0)},
    \qquad
    \Sigma_{vu}^{(0)}(\tau)\to(1-\overline\Theta^{\rm eq})\Sigma_u^{(0)},
    \qquad
    \Sigma_v^{(0)}(\tau)
    \to
    \frac{2\left(1-\overline\Theta^{\rm eq}\right)^2\Sigma_u^{(0)}+\nu^2\left(1-\left(m_u^{(0)}\right)^2\right)}
    {2+\nu^2}.
\end{equation*}

\noindent{If $\overline\Theta^{\rm eq}=1$:}
\begin{equation*}
m_v^{(0)}[\rho^0](\tau)\equiv m_v^{(0)}[\rho^0](0), \qquad \Sigma_{vu}^{(0)}(\tau)\to0, \qquad
    \Sigma_v^{(0)}(\tau)
    \to
    \frac{\nu^2\left(1-\left(m_v^{(0)}[\rho^0](0)\right)^2\right)}{2+\nu^2}.
\end{equation*}
\end{proposition}

\begin{proof}
The result is obtained by solving the moment equations sequentially; the mean equation \eqref{eq:mean-opinion-fast-strategy} determines the limiting opinion mean, and substitution into the covariance and variance equations \eqref{eq:fast-strategy-moments-general}, followed by elementary algebra, yields the stated stationary limits.
\end{proof}

Thus, at the macroscopic level, the separation-of-scales assumption (A3) removes the strategy-induced nonclosure, in particular the fast strategy dynamics is encoded only through $ \overline\Theta^{\rm eq}$. Consequently, once $g^{\rm eq}$ is selected by the projective flow, any remaining lack of moment closure is due solely to the structure of the diffusion coefficient $D$.

We emphasize, however, that this exact closure property holds for the reduced equation obtained under (A3). It does not imply that the resulting moment equations reproduce exactly the second-order moments of the original kinetic model when $\Theta(\sigma)$ varies across the support of $g^{\rm eq}$. The accuracy of this additional factorization assumption
is discussed below.

\begin{remark}[Accuracy of the fast-strategy closure for
second-order moments]
\label{rem:A3-defect}

The closure (A3) may fail to reproduce the second-order moments of the full kinetic model, even in the fast-strategy limit. This can be quantified explicitly for $K=2$ with constant payoffs in the coexistence regime. Throughout, moments without superscript refer to the full distribution $f$, while moments with superscript $(0)$ refer to the reduced model under~(A3).

\emph{Setting.} In this case, $\Theta(s)=s$ and
\(
    g^{\rm eq}
    =
    (\Phi_{\mathcal T_\infty})_\#g_0
\)
is non-degenerate, with $\overline\Theta^{\rm eq}=\hat m_s\in(0,1)$.
By assumption (A1) we have initial independence between strategies and convictions,
\(    \eta_0(s,u)=g_0(s)k(u),
\) which is preserved by the strategy dynamics for constant payoffs (Corollary~\ref{cor:eta-pushforward}), hence $\Sigma_{su}(t)\equiv0$. We denote
\[
\Sigma_s^{\rm eq}    :=
 \int_{\Delta_1}
    \left(s-\overline\Theta^{\rm eq}\right)^2
    g^{\rm eq}(s)\,ds,
\]

\emph{Exact stationary moments of the full model.} For the full kinetic model, the accumulated effective compromise $\mathcal I(t)$ diverges and the exact mean opinion equation gives $m_v^\infty=m_u$ (Section~\ref{sec:2s-constant}). The asymptotic distribution of the full model is therefore
\[
    f^\infty(v,s,u)
    =
    \rho^\infty(v\mid s,u)
    g^{\rm eq}(s)k(u),
\]
where $\rho^\infty(v\mid s,u)$ is the Gibbs density
with attraction center
\(
    \lambda(s,u)=s\,m_u+(1-s)u.
\)
In particular, for $D(v)=\sqrt{1-v^2}$,
\[
     \int_{I} v\,\rho^\infty(v\mid s,u)\,\mathrm{d}v=\lambda(s,u),
    \qquad
     \int_{I}
    \left(v-\lambda(s,u)\right)^2
    \rho^\infty(v\mid s,u)\,\mathrm{d}v
    =
    \frac{\nu^2}{2+\nu^2}
    \left(1-\lambda(s,u)^2\right).
\]
Since $\lambda(s,u)-m_u=(1-s)(u-m_u)$, the laws of total expectation, covariance and variance give
\begin{equation}\label{eq:Sv-exact-coexistence}
\begin{aligned}
    \Sigma_v^\infty
    &=
    \frac{
        2\left[
            \left(1-\overline\Theta^{\rm eq}\right)^2
            +\Sigma_s^{\rm eq}
        \right]\Sigma_u
        +\nu^2(1-m_u^2)
    }{2+\nu^2},\qquad 
    \Sigma_{vu}^\infty
    &=
    \left(1-\overline\Theta^{\rm eq}\right)\Sigma_u.
\end{aligned}
\end{equation}
These values hold for every ratio $\tau_v/\tau_\sigma\in(0,\infty)$.

\emph{Comparison with the reduced model.} Since the conviction marginal is conserved in both the full and reduced dynamics, the initial independence assumption implies that the two models share the same conviction marginal. Hence, their conviction moments
coincide, namely
\[
    m_u=m_u^{(0)},\qquad
    \Sigma_u=\Sigma_u^{(0)}.
\] Therefore, comparing \eqref{eq:Sv-exact-coexistence} with
Proposition~\ref{prop:fast-strategy-stationary}, with
$\Sigma_v^{(0),\infty}:=\lim_{\tau\to+\infty}\Sigma_v^{(0)}(\tau)$, the closure reproduces the exact stationary mean and opinion--conviction covariance, but underestimates the opinion
variance:
\begin{equation}\label{eq:Sv-defect}
    \Sigma_v^\infty-\Sigma_v^{(0),\infty}
    =
    \frac{2}{2+\nu^2}\Sigma_s^{\rm eq}\Sigma_u.
\end{equation}
Thus, at the level of the stationary variance, the closure
is exact if and only if
$\Sigma_s^{\rm eq}\Sigma_u=0$. In the coexistence regime,
$\Sigma_s^{\rm eq}>0$, and this condition reduces to
$\Sigma_u=0$.

\emph{Origin of the defect.} The origin of this discrepancy can be traced back to the
mixed second-order moment
\[
    E(\tau)
    := \int_\Omega
s\left(v-m_u\right)
\left(u-m_u\right)
f
\,\mathrm{d}\omega, 
\] and $\Sigma_{sv}$, which (A3) replaces by $\overline\Theta^{\rm eq}\Sigma_{vu}$ and $0$, respectively.
At stationarity $\Sigma_{sv}^\infty=0$
by independence of $s$ and $u$. The mixed moment, instead, equals
\[
    \lim_{\tau\to+\infty}E(\tau)=E^\infty
    =
\left(\int_{\Delta_1}s(1-s)g^{\rm eq}(s)\,ds\right)
\left(\int_I(u-m_u^{(0)})^2k(u)\,du\right)
    =
    \left[
        \overline\Theta^{\rm eq}
        (1-\overline\Theta^{\rm eq})
        -\Sigma_s^{\rm eq}
    \right]\Sigma_u,
\]
whereas (A3) replaces this quantity with
\[
E_{\mathrm{A3}}^\infty
    :=
\overline\Theta^{\rm eq}\Sigma_{vu}^\infty
    =
\overline\Theta^{\rm eq}
    (1-\overline\Theta^{\rm eq})\Sigma_u.
\]
Hence, the discrepancy
in the mixed moment is precisely
\[
    E^\infty-E_{\mathrm{A3}}^\infty
    =
    -\Sigma_s^{\rm eq}\Sigma_u.
\]
At stationarity the variance equation reads $(2+\nu^2)\Sigma_v^\infty=2\Sigma_{vu}^\infty-2E^\infty+\nu^2(1-m_u^2)$, so that $\Sigma_v^\infty-\Sigma_v^{(0),\infty}=-2(E^\infty-E^\infty_{\rm A3})/(2+\nu^2)$, which is~\eqref{eq:Sv-defect}: the whole variance defect comes from the mixed moment.

\emph{The defect is not small in $\epsilon_\sigma$.} This discrepancy does not vanish as $\epsilon_\sigma=\tau_\sigma/\tau_v\to0$. Indeed, for constant payoffs, the limiting strategy distribution remains non-degenerate in the coexistence regime: each agent retains a strategy $s^\infty=\Phi_{\mathcal T_\infty}(s_0)$, and the corresponding opinion equilibrium depends on this individual strategy through $\lambda(s,u)$. Therefore, although the fast-strategy closure captures the stationary opinion mean and opinion--conviction covariance exactly, it introduces a persistent error in the stationary opinion variance. In the coexistence regime, (A3) is consequently an $\mathcal O(1)$ approximation for this second-order moment, rather than an $\mathcal O(\epsilon_\sigma)$ approximation. We stress that this is a limitation of the closure (A3) on
$f^0$, not of the Chapman--Enskog expansion itself.

\emph{Dominance endpoints.} For $\overline\Theta^{\rm eq}\in\{0,1\}$ the limiting strategy distribution is a Dirac mass, so $\Sigma_s^{\rm eq}=0$ and there is no $\mathcal O(1)$ defect. For $\overline\Theta^{\rm eq}=0$ the exact limits coincide with those of Proposition~\ref{prop:fast-strategy-stationary}. For $\overline\Theta^{\rm eq}=1$ the exact stationary distribution is a Gibbs density centered at $m_v^\infty=m_u+(m_v(0)-m_u)\mathrm{e}^{-\mathcal I_\infty/\tau_v}$ with $\mathcal I_\infty=\mathcal O(\tau_\sigma)$ (Section~\ref{sec:2s-constant}), whereas the reduced model conserves $m_v^{(0)}(0)$. The two differ by $\mathcal O(\varepsilon_\sigma)$, because of the initial layer of duration $\mathcal O(\tau_\sigma)$ in which $m_s<1$; the same holds for the stationary variance.
\end{remark}

\section{Analysis of the two- and three-strategy models}
\label{sec:exact-analysis}

The asymptotic analysis of Section~\ref{sec:asymptotic} provides approximate reduced descriptions of the coupled opinion-strategy dynamics in the limiting regimes $\tau_v\ll\tau_\sigma$ and $\tau_\sigma\ll\tau_v$. The present section complements that analysis by investigating specific model configurations at comparable time scales $\tau_v\sim\tau_\sigma$, for which exact representations or detailed qualitative characterizations are available.

The key observation is that for constant payoffs $J\in\mathbb{R}^{K\times K}$, Theorem~\ref{thm:projective} provides an exact closed representation of the strategy marginal $g(\sigma,t)$ that holds for \emph{all} values of $\tau_v/\tau_\sigma$, without any scale separation assumption. Combined with the exact evolution of the strategy--conviction marginal $\eta$ and the moment equation for the mean opinion $m_v(t)$, this representation allows us to classify the asymptotic strategic regimes and characterize their consequences for the macroscopic opinion dynamics.

For state-dependent payoffs, the strategy marginal no longer evolves autonomously and exact solvability is lost. In this case we focus on a qualitative characterization of the time-dependent segregation boundary and the structure of candidate stationary states.

The section is organized as follows. In Section~\ref{sec:2s-constant} we analyze the two-strategy model ($K=2$) with constant payoffs: the projective flow classifies all possible asymptotic behaviors of $g(s,t)$ into four regimes determined by the sign structure of the gain $G$, and the resulting dynamics of the mean opinion $m_v(t)$ is characterized exactly. Section~\ref{sec:state-dependent} considers the two-strategy model with a benchmark state-dependent payoff matrix; the autonomous closure for $g$ is no longer available, but the geometry of the time-dependent segregation boundary and the structure of stationary states can be characterized explicitly. Section~\ref{sec:three-strategy} extends the exact constant-payoff analysis to the three-strategy model ($K=3$), where the richer geometry of the two-dimensional simplex generates new dynamical phenomena, including stable interior coexistence and cyclic competition, with no direct analogue in the two-strategy setting. Numerical illustrations of the representative cases are provided in Section~\ref{sec:numerics}.

\subsection{Two-strategy model with constant payoffs}
\label{sec:2s-constant}

We first consider a  $K=2$ setting in which the two strategies are identified with the two behavioral responses considered in the model, setting \( \Pi=\{\pi_C,\pi_S\},\) where $\pi_C$ and $\pi_S$ denote the compromise and stubbornness strategies, respectively, with compromise parameters
\(
(\theta_C,\theta_S)=(1,0)
\).

The simplex $\Delta_1=[0,1]$ is
parametrized, in reduced coordinates, by the scalar compromise frequency
$s:=\sigma_C$, with frequency of the stubborn strategy $\sigma_S=1-s$ determined by the normalization constraint. 

The joint  density $f=f(v,s,u,t)$ defined on \( \Omega=I\times\Delta_1\times I \) satisfies the VFP equation \eqref{eq:FP_gain}, which reduces to:
\begin{equation}\label{eq:FP-two}
\partial_t f
=
-\frac{1}{\tau_v}
\partial_v
\left(
\left(
s(m_v-v)+(1-s)(u-v)
\right)f
-\frac{\nu^2}{2}
\partial_v\left(D(v)^2f\right)
\right)
-\frac{1}{\tau_\sigma}
\partial_s
\left(
s(1-s)\mathcal{G}[f]f
\right).
\end{equation}
where
\(\mathcal{G}[f]
:=
\mathcal{P}_C[f]-\mathcal{P}_S[f]
\)
represents the relative payoff advantage of $\pi_C$ over $\pi_S$.

\paragraph{Constant payoff matrix.} We consider the baseline scenario where the game interaction is state-independent, considering a constant payoff matrix $J\in\mathbb{R}^{2\times2}$ with  entries $J_{CC}, J_{CS}, J_{SC}, J_{SS}$.

In this case, the gain  $\mathcal{G}[f](t)=G(m_s(t))$, since $G$ is affine in $m_s$ and the integration of $f$ against the opponent's strategy yields only $m_s$ (as shown in Section~\ref{sec:strategy-flow}):
\begin{equation*}
G(m_s(t))=
m_s(t)(J_{CC}-J_{SC})
+(1-m_s(t))(J_{CS}-J_{SS}),\qquad m_s(t):=\int_{\Delta_1} \,sg(s,t)\,\mathrm{d}s
\end{equation*} where $g(s,t)$ denotes the strategy marginal.

\paragraph{Exact strategy dynamics.} Since we consider constant payoffs, Theorem~\ref{thm:projective} applies directly,  yielding the exact push-forward representation
\begin{equation}\label{eq:g-pushforward-two}
    g(s,t)
    =
    \bigl(\Phi_{\mathcal{T}(t)}\bigr)_\#g_0(s),
    \qquad \mathcal{T}(t)=\frac{1}{\tau_\sigma}
\int_0^t G(m_s(\tau))\,\mathrm{d}\tau, \qquad
    m_s(t)=\mathcal{E}(\mathcal{T}(t)).
\end{equation}
The projective map specializes to
\begin{equation}\label{eq:two-strategy-flow}
    \Phi_{\mathcal{T}}(s_0)
    =
    \frac{\mathrm{e}^{\mathcal{T}}s_0}{1-s_0+\mathrm{e}^{\mathcal{T}}s_0},
\end{equation}
and the cumulative gain $\mathcal{T}(t)$ satisfies the system \eqref{cumulative_gains} which reduces to a scalar ODE, given by
\begin{equation}\label{eq:cumulative-gain-equation-two}
    \mathcal T'(t)
    =\frac{1}{\tau_\sigma}
    G\bigl(\mathcal E(\mathcal T(t))\bigr),
    \qquad
    \mathcal T(0)=0.
\end{equation} Moreover,  Corollary \ref{cor:eta-pushforward} extends this push-forward structure to the strategy--conviction marginal, yielding
\begin{equation}\label{eq:eta-pushforward-two}
    \eta(s,u,t)
    =
    \bigl(\Phi_{\mathcal{T}(t)},\mathrm{Id}\bigr)_\#\eta_0(s,u).
\end{equation} In particular, the strategy--conviction covariance \(\Sigma_{su}(t)\) can be determined exactly at every time from this representation.

\paragraph{Coupling with opinion dynamics.}
The macroscopic moment equations for $m_v$ \eqref{eq:mean_opinion} and $m_s$ \eqref{eq:mean_strategy_constant} reduce to
\begin{align}
\frac{d}{dt} m_v(t) &= \frac{1}{\tau_v} \left( \left(1 - m_s(t)\right) \left(m_u - m_v(t)\right) -\Sigma_{su}(t)\right), \label{eq:moment-v}\\
\frac{d}{dt} m_s(t) &= \frac{1}{\tau_\sigma} {G}(m_s(t)) \left( m_s(t)-m_s^{(2)}(t)\right), \label{eq:moment-s} 
\end{align}
where  $m_s^{(2)}(t):=\int_{\Delta_1} s^2 g(s,t)ds$ is the second-order moment of the strategy distribution.

Introducing the strategy variance $\Sigma_s(t):=m_s^{(2)}(t)-m_s(t)^2$, we have $m_s-m_s^{(2)} =m_s(1-m_s)-\Sigma_s\geq0$ since $s^2\leq s$ on $[0,1]$, which gives $m_s^{(2)}\leq m_s$. Therefore, equation~\eqref{eq:moment-s} shows that the direction of  evolution of the mean strategy (whenever $g$ is not concentrated at $s\in\{0,1\}$) is determined by the sign of $G(m_s)$, while a larger variance $\Sigma_s(t)$ reduces the effective rate of strategic adaptation. Note that $m_s^{(2)}(t)$ is itself determined exactly by the push-forward representation~\eqref{eq:g-pushforward-two}, so that equation~\eqref{eq:moment-s} is not an independent equation
but a consequence of the exact strategy dynamics.

Introducing $ \mathcal I(t) := \int_0^t \left(1- m_s(\tau)\right)\,d\tau$, equation~\eqref{eq:moment-v} can be integrated exactly, yielding:
\begin{equation}\label{eq:mv-exact}
m_v(t) = m_u + \bigl(m_v(0)-m_u\bigr)\mathrm{e}^{-\mathcal{I}(t)/\tau_v} -\frac{1}{\tau_v} \mathrm{e}^{-\mathcal{I}(t)/\tau_v} \int_0^t {\Sigma_{su}(\tau)} \mathrm{e}^{\mathcal{I}(\tau)/\tau_v}  \mathrm{d}\tau.\end{equation} Hence, the asymptotic mean opinion $m_v^\infty := \lim_{t\to\infty} m_v(t)$ (whenever the limit exists)  depends on the long-time behavior of $\mathcal{I}(t)$, and on the accumulated contribution of $\Sigma_{su}(\tau)$.

\paragraph{Asymptotic regimes classification.}
The long-time behavior of $g$ is determined by  the asymptotic behavior of the cumulative gain $\lim_{t\to\infty}\mathcal{T}(t)$, which in turn is governed by the sign structure of the affine gain $G$. Setting
\[
    G(0)=J_{CS}-J_{SS},
    \qquad
    G(1)=J_{CC}-J_{SC},
    \qquad
    G'=G(1)-G(0),
\]
we distinguish four non-degenerate cases according to the signs of \(G(0)\) and \(G(1)\).

When ${G}$ has a constant sign on $[0,1]$, the cumulative gain \(\mathcal T(t)\) diverges to \(+\infty\) or \(-\infty\), and the strategy distribution converges to a pure strategic state \(s^\infty\in\{0,1\}\), provided the initial distribution has no mass at the repelling boundary. The limiting strategy--conviction marginal then factorizes as $$\eta^\infty(s,u)=g^\infty(s)k(u) \qquad g^\infty(s)=\delta(s-s^\infty).$$ This factorization follows from the asymptotic concentration of the strategy distribution at a single pure state; at finite times, in contrast, it holds only if strategy and conviction are initially independent. We thus obtain the  following two dominance regimes.

\subparagraph{Regime~1: Compromise dominance.} If
$G(0)>0$ and $G(1)>0$, then $G(m_s)>0$ for all $m_s\in[0,1]$. Consequently,
$\mathcal{T}(t)\to+\infty$ and
\[
    g^\infty=\delta(s-1),\qquad m_s^\infty=1.
\]
Since $m_s(t)\to1$ exponentially,  the integral
$\mathcal{I}$ converges to a finite limit $\mathcal{I}_\infty\in\mathbb{R}$. Hence, the asymptotic mean opinion retains memory of the initial condition and, in general, of the transient strategy--conviction correlation, and its value depends on $\tau_\sigma$ through $\mathcal{I}_\infty$. More precisely, $m_v^\infty:=m_v^{\mathrm{comp}}\in I$ where, 
\begin{equation}\label{eq:mv-infty-comp}
    m_v^{\rm comp}
    :=
    m_u
    +\bigl(m_v(0)-m_u\bigr)\mathrm{e}^{-\mathcal{I}_\infty/\tau_v}
    -\frac{1}{\tau_v}\mathrm{e}^{-\mathcal{I}_\infty/\tau_v}
    \int_0^\infty\Sigma_{su}(\tau)\,\mathrm{e}^{\mathcal{I}(\tau)/\tau_v}\,\mathrm{d}\tau.
\end{equation}
The integral in \eqref{eq:mv-infty-comp} converges absolutely: since
$\mathrm{e}^{\mathcal{I}(\tau)/\tau_v}\leq \mathrm{e}^{\mathcal{I}_\infty/\tau_v}
<+\infty$ and $|\Sigma_{su}(\tau)|\leq\sqrt{\Sigma_s(\tau)\Sigma_u}
\to0$ exponentially, the integrand is integrable on $[0,+\infty)$.

\subparagraph{Regime~2: Stubbornness dominance.}
If $G(0)<0$ and $G(1)<0$, then $G(m_s)<0$ for all $m_s$, so
$\mathcal{T}(t)\to-\infty$ and
\[
    g^\infty=\delta(s),\qquad m_s^\infty=0.
\] Since $m_s(t)\to0$ exponentially, $\mathcal{I}(t)\sim t\to+\infty$. By~\eqref{eq:mv-exact}, both the initial-condition term and the covariance integral are asymptotically suppressed, and the opinion dynamics relaxes to the mean conviction: $m_v^\infty=m_u$.

\begin{remark}[Transient role of $\Sigma_{su}$]
In both dominance regimes, $\Sigma_{su}(t)\to0$ as $t\to+\infty$, since $\Sigma_s(t)\to0$ and $|\Sigma_{su}|^2\leq \Sigma_s\cdot\Sigma_u$ by the Cauchy--Schwarz inequality. However, in Regime~1 the finiteness of $\mathcal{I}_\infty$ means that the transient contribution of $\Sigma_{su}$ is not damped by an increasingly large integrating factor in~\eqref{eq:mv-exact} and may persist in $m_v^{\rm comp}$, instead in Regime~2, $\mathcal{I}(t)\to+\infty$ asymptotically damps the entire transient.
\end{remark}

When $G(0)$ and $G(1)$ have opposite signs, $G$ has a unique zero
\begin{equation}\label{eq:mhat}
\hat m_s
=
\frac{J_{SS}-J_{CS}}{(J_{CC}-J_{SC})+(J_{SS}-J_{CS})}
\in(0,1),
\end{equation}
which is stable if $G'<0$ and unstable if $G'>0$. These two cases give rise to the following regimes.

\subparagraph{Regime~3: Stable coexistence (anti-coordination).} If
$G(0)>0$ and $G(1)<0$, then $G'<0$, and $G$ is strictly decreasing with a stable zero $\hat m_s$, defined in \eqref{eq:mhat}. For any initial distribution with $g_0(1)<\hat m_s<1-g_0(0)$ (in particular if $g_0$ is absolutely continuous) then $\mathcal E(\mathcal T)=\hat m_s$ has a unique solution $\mathcal T_\infty\in\mathbb R$ and $m_s(t)\to\hat m_s$,
$\mathcal{T}(t)\to\mathcal{T}_\infty$. Consequently,
\[
    g^\infty(s)
    =
    \bigl(\Phi_{\mathcal{T}_\infty}\bigr)_\#g_0(s),
    \qquad
    \eta^\infty(s,u)
    =
    \bigl(\Phi_{\mathcal{T}_\infty},\mathrm{Id}\bigr)_\#\eta_0(s,u).
\]
The asymptotic strategy distribution is therefore non-degenerate and retains the qualitative structure of the initial distribution under the projective transformation. In particular, the asymptotic strategy--conviction covariance converges to a finite limit $\Sigma_{su}^\infty$. Since $m_s^\infty=\hat m_s<1$, we have $\mathcal{I}(t)\sim(1-\hat m_s)t\to+\infty$ and the initial-condition contribution in \eqref{eq:mv-exact} is damped, giving
\begin{equation*}
    m_v^\infty
    =
    m_u-\frac{\Sigma_{su}^\infty}{1-\hat m_s}.
\end{equation*}
If the initial mean satisfies $m_s(0)=\hat m_s$ exactly, then $\mathcal{T}(t)\equiv0$ by uniqueness of the solution of  the ODE~\eqref{eq:cumulative-gain-equation-two}. Hence, $g(s,t)=g_0(s)$, and $\Sigma_{su}(t)=\Sigma_{su}(0)$ for all $t$, yielding $m_v^\infty=m_u-\Sigma_{su}(0)/(1-\hat m_s)$. 

In particular, the atoms at $\{0,1\}$, which are fixed points in the strategy flow, prevent the mean from reaching $\hat m_s$. For instance, if $g_0(1)\ge\hat m_s$ then $\mathcal T\to-\infty$ and $g^\infty=g_0(1)\delta_1+(1-g_0(1))\delta_0$, and symmetrically if $g_0(0)\ge1-\hat m_s$.

\subparagraph{Regime~4: Bistability (coordination).}
If $G(0)<0$ and $G(1)>0$, then $G'>0$, and $G$ is strictly increasing with an unstable zero $\hat m_s$, given by \eqref{eq:mhat}. This threshold  separates the two basins of attraction, and the asymptotic behavior  depends on the initial mean strategy $m_s(0)$.

\emph{Off-threshold initial conditions.} If $m_s(0)>\hat m_s$ and $g_0(0)=0$, then $\mathcal T(t)\to+\infty$ and $m_s(t)\to1$, so the asymptotic strategic state is the same as in Regime~1, and $m_v^\infty=m_v^{\rm comp}$. Conversely, if $m_s(0)<\hat m_s$ and $g_0(1)=0$,  then $\mathcal T(t)\to-\infty$ and $m_s(t)\to0$, as in Regime~2, giving $m_v^\infty=m_u$. Otherwise, if $g_0$ admits a finite mass at $\{0,1\}$ the limiting distribution are given in the two cases by $g^\infty=g_0(0)\delta_0+(1-g_0(0))\delta_1$, and $g^\infty=(1-g_0(1))\delta_0+g_0(1)\delta_1$ respectively.

\emph{Threshold initial condition.} If $m_s(0)=\hat m_s$, again $\mathcal{T}(t)\equiv0$ by uniqueness, so that $g(s,t)=g_0(s)$, and $\Sigma_{su}(t)=\Sigma_{su}(0)$ for all $t$. Since $m_s(t)\equiv \hat m_s$, the mean opinion converges to $m_v^\infty=m_u-\Sigma_{su}(0)/(1-\hat m_s)$.

\begin{remark}[Convergence rates]
\label{rem:convergence-rates}
The exact projective representation determines the convergence rate of $m_s(t)$ in each regime. In Regimes~1 and~2, if the initial strategy distribution is supported away from the repelling boundary, then
\[
1-m_s(t)=\mathcal O\!\left(\mathrm{e}^{-\mu_s t}\right),
\qquad
\mu_s=\frac{G(1)}{\tau_\sigma},
\quad\text{(Regime~1),}
\]
and
\[
m_s(t)=\mathcal O\!\left(\mathrm{e}^{-\mu_s t}\right),
\qquad
\mu_s=\frac{|G(0)|}{\tau_\sigma},
\quad\text{(Regime~2).}
\]
Without this additional assumption, convergence to the
attracting pure strategy still holds, provided there is no
initial mass at the repelling boundary, but the convergence
rate need not be exponential.

In Regime~3, linearizing the scalar
ODE~\eqref{eq:cumulative-gain-equation-two} around $\mathcal T_\infty$ yields exponential convergence with rate
\[
|m_s(t)-\hat m_s|
=\mathcal O\!\left(\mathrm{e}^{-\mu_s t}\right),
\qquad
\mu_s
=
\frac{|G'(\hat m_s)|\,\mathcal{E}'(\mathcal{T}_\infty)}
{\tau_\sigma},
\]
provided $\mathcal E'(\mathcal T_\infty)>0$. In Regime~4,
the corresponding estimates apply to off-threshold initial
conditions, whereas $m_s(t)\equiv\hat m_s$ at the threshold.

The exact moment equation~\eqref{eq:moment-v} then yields
exponential convergence of $m_v(t)$ whenever the strategy
moments and the covariance $\Sigma_{su}(t)$ converge
exponentially. In particular, this holds in Regime~3 and
at the threshold of Regime~4. In the dominance regimes,
it holds under the additional support assumption stated
above. The corresponding limiting values are those
characterized in the preceding analysis.
\end{remark}

\begin{remark}[Initial independence]
If $\eta_0(s,u)=g_0(s)k(u)$, then $\Sigma_{su}(t)\equiv0$ for all $t\geq0$ by the push-forward representation~\eqref{eq:eta-pushforward-two}, since $\int(u-m_u)k(u)\,du=0$. In this case $m_v^\infty=m_u$ in all regimes except Regime~1 and Regime~4 (compromise basin), where $m_v^{\rm comp}=m_u+(m_v(0)-m_u)\mathrm{e}^{-\mathcal{I}_\infty/\tau_v}$.
\end{remark}

Table~\ref{tab:regimes_summary} summarizes the four regimes.

\begin{remark}[Degenerate payoff conditions]
The classification above assumes $G(0)\neq0$ and $G(1)\neq0$. If $G\equiv0$ ($J_{CC}=J_{SC}$ and $J_{CS}=J_{SS}$), the strategic dynamics is completely frozen and $g(s,t)=g_0(s)$ for all $t$. In the remaining boundary cases (exactly one of $G(0),G(1)$ equals zero), the qualitative attractor is unchanged but the convergence to the attracting pure strategy may become algebraic, causing $\mathcal{I}(t)$ to grow like $\log t$ rather than remaining bounded. This changes the convergence rate of $m_v(t)$ but does not introduce new asymptotic strategic attractors.
\end{remark}

\begin{table}[t]
\centering
\setlength{\tabcolsep}{3pt}
\renewcommand{\arraystretch}{1.4}
\begin{tabular}{ccccc}
\toprule
\textbf{Regime}
&
\textbf{Payoff conditions}
&
$\boldsymbol{g^\infty(s)}$
&
$\boldsymbol{m_s^\infty}$
&
$\boldsymbol{m_v^\infty}$
\\
\midrule
\shortstack{\textbf{1}\\\footnotesize Compromise}
&
$\begin{matrix}
J_{CC}>J_{SC}\\
J_{CS}>J_{SS}
\end{matrix}$
&
$\delta(s-1)$
&
$1$
&
$\displaystyle
m_v^{\mathrm{comp}}$
\\[12pt]
\shortstack{\textbf{2}\\\footnotesize Stubbornness}
&
$\begin{matrix}
J_{CC}<J_{SC}\\
J_{CS}<J_{SS}
\end{matrix}$
&
$\delta(s)$
&
$0$
&
$\displaystyle
m_u$
\\[12pt]
\shortstack{\textbf{3}\\\footnotesize Coexistence}
&
$\begin{matrix}
J_{CC}<J_{SC}\\
J_{CS}>J_{SS}
\end{matrix}$
&
$(\Phi_{\mathcal T_\infty})_\# g_0(s)$
&
$\hat m_s$
&
$\displaystyle
m_u-
\frac{\Sigma_{su}^\infty}{1-\hat m_s}$
\\[12pt]
\shortstack{\textbf{4}\\\footnotesize Bistability}
&
$\begin{matrix}
J_{CC}>J_{SC}\\
J_{CS}<J_{SS}
\end{matrix}$
&
$\begin{cases}
\delta(s-1), & m_s(0)>\hat m_s\\
g_0(s), & m_s(0)=\hat m_s \\
\delta(s), & m_s(0)<\hat m_s
\end{cases}$
&
$\begin{cases}
1, & m_s(0)>\hat m_s\\
\hat m_s, &   m_s(0)=\hat m_s\\
0, & m_s(0)<\hat m_s
\end{cases}$
&
$\begin{cases}
m_v^{\mathrm{comp}}, & m_s(0)>\hat m_s\\
m_u-\displaystyle
\frac{\Sigma_{su}(0)}{1-\hat m_s}, &  m_s(0)=\hat m_s \\ 
m_u, & m_s(0)<\hat m_s
\end{cases}$\\
\bottomrule
\end{tabular}
\caption{Asymptotic strategy distribution, mean strategy, and mean opinion for the two-strategy model with constant payoffs.}
\label{tab:regimes_summary}
\end{table}

\paragraph{Asymptotic joint distribution.}
The evolution of the strategy--conviction marginal $\eta(s,u,t)$ is completely determined by~\eqref{eq:eta-pushforward-two}. To obtain the full
distribution $f$, we introduce the conditional opinion
density
\[
    \rho(v\mid s,u,t):=\frac{f(v,s,u,t)}{\eta(s,u,t)},
\]
well-defined wherever $\eta>0$. Since $\eta$ is independent of $v$, a direct computation using the VFP equation~\eqref{eq:FP_gain} and the transport equation~\eqref{eq:eta-transport-general} for $\eta$
gives \begin{equation*}
\partial_t\rho
+\frac{1}{\tau_\sigma}s(1-s)G(m_s)\,\partial_s\rho
=\frac{1}{\tau_v}\!\left(
-\partial_v(\mathcal{A}^v[f]\rho)
+\frac{\nu^2}{2}\partial_{vv}(D^2\rho)
\right).
\end{equation*}
The left-hand side is the material derivative along the strategic characteristic $S(t)=\Phi_{\mathcal{T}(t)}(S_0)$. Therefore, along each characteristic, $\rho$ satisfies the non-autonomous Fokker--Planck equation in $v$, given by:
\[
\frac{d}{dt}\rho(v\mid S(t),u,t)
=\frac{1}{\tau_v}\!\left(
-\partial_v\!\bigl((a(t)-v)\rho\bigr)
+\frac{\nu^2}{2}\partial_{vv}(D^2\rho)
\right),
\quad
a(t):=S(t)m_v(t)+(1-S(t))u,
\]
where $S(t)$ and $m_v(t)$ are known exactly from the preceding analysis.  In all regimes for which the convergence of $S(t)\to s^\infty$ and $m_v(t)\to m_v^\infty$ is exponential, one can adapt the relative entropy argument in \cite{FURIOLI20192065} to provide convergence of $\rho$ to the limiting Gibbs density for the autonomous Fokker--Planck equation, $\rho(v\mid S(t),u,t)\to\rho^\infty(v\mid s^\infty,u)$
in $L^1(I)$, where $\rho^\infty$ is the Gibbs density with center:
\begin{equation}\label{eq:rho_cond_Gibbs}
\rho^\infty(v\mid s,u)
=\frac{C(s,u)}{D(v)^2}
\exp\!\left(
\frac{2}{\nu^2}\int_0^v
\frac{V_\infty(s,u)-z}{D(z)^2}\,\mathrm{d}z
\right), \qquad V_\infty(s,u):=s m_v^\infty+(1-s)u.
\end{equation}

Integrating out the strategy variable yields the asymptotic opinion profile of the subpopulation with conviction $u$. Let $\eta^\infty(ds\mid u)$ be
the conditional distribution of the strategy given $u$ in the asymptotic regime. Then
\begin{equation*}
\rho^\infty(v\mid u)
=
\int_{\Delta_1}
\rho^\infty(v\mid s,u)\,
\eta^\infty(\mathrm{d}s\mid u),
\end{equation*}
for $k(u)>0$. Thus, the asymptotic opinion profile at fixed conviction is a mixture of the conditional Gibbs densities associated with the different strategies represented in the limiting
population. The mixture is determined by the asymptotic strategy--conviction marginal $\eta^\infty$, while the mean opinion $m_v^\infty$ enters the attraction center $V_\infty(s,u)$.

We now distinguish two cases based on whether $\mathcal{T}_\infty$ is finite or infinite.

\emph{Regimes~1, 2, and off-threshold Regime~4 ($\mathcal{T}_\infty=\pm\infty$).}
All characteristics starting from $S_0\in(0,1)$ converge to the same $s^\infty\in\{0,1\}$, so $\rho(v\mid S(t),u,t)\to\rho^\infty(v\mid s^\infty,u)$ along every characteristic. In the weak sense of measures:
\[
f(v,s,u,t)\rightharpoonup
\rho^\infty(v\mid s^\infty,u)\,\delta(s-s^\infty)\,k_0(u)
=:f^\infty(v,s,u).
\]
In particular in Regime~1, $V_{\infty}(1,u)=m_v^{\rm comp}$ which is independent of $u$ and all agents share the same opinion profile regardless of conviction. In Regime~2, $V_{\infty}(0,u)=u$ and each agent's opinion is attracted toward their private conviction.

\emph{Regime~3 ($\mathcal{T}_\infty\in\mathbb{R}$).}
The projective map $\Phi_{\mathcal{T}_\infty}$ is a diffeomorphism of $(0,1)$. The characteristic arriving at $s$ as time $t\to\infty$ originated from $S_0=\Phi_{-\mathcal{T}_\infty}(s)$, so
\begin{align}\label{eq:finfty-factored}
f(v,s,u,t)\to
\rho^\infty(v\mid s,u)\,\eta^\infty(s,u)
=:f^\infty(v,s,u),
\end{align}
where $\eta$ is absolutely continuous whenever $\eta_0$ is. In this case $\rho^\infty(v\mid u)$ is a nontrivial mixture of Gibbs densities with a continuum of attraction centers $V_{\infty}(s,u)=sm_v^{\infty}+(1-s)u$ indexed by $(s,u)$.

In all cases, the asymptotic factorization $f^\infty=\rho^\infty\eta^\infty$ holds in the appropriate sense, with $\eta^\infty$ determined by the push-forward representation and $\rho^\infty$ given by~\eqref{eq:rho_cond_Gibbs}.

\subsection{Two-strategy model with state-dependent payoffs}\label{sec:state-dependent}
We retain the two behavioral strategies  \(\Pi=\{\pi_C,\pi_S\},\) corresponding to the compromise and stubbornness strategies, respectively, with compromise parameters $\theta=(1,0)$, and now allow the payoff entries to depend on the states of both interacting agents: $J = \left(J_{ij}(v,u,v_*,u_*)\right)_{i,j \in \{C,S\}}$.

\paragraph{Benchmark payoff specification.} We consider the quadratic payoff family  \eqref{eq:payoffs-micro}, in which each strategy incurs a  cost associated with the type of disagreement that it is designed to reduce. The compromise strategy penalizes discrepancy from the opponent's opinion, while the stubbornness strategy penalizes discrepancy from the agent's own conviction. 

Specifically, we impose $a_{CC}=a_{CS}=0$ and $b_{SC}=b_{SS}=0$, so that compromise payoffs do not  account for the distance from the private conviction, and stubbornness payoffs do not  account for the opponent's opinion. We further set $b_{CC} = b_{CS} = b > 0$ and $a_{SC} = a_{SS} = a > 0$, so that the cost of each type of discrepancy is independent of the opponent's strategy. The parameters $c_{CS},c_{SC}\geq0$ capture additional costs in cross-strategy interactions, representing, for instance, social or psychological frictions arising when the two strategies interact. 

The resulting payoff matrix is 
\[
J =
\begin{pmatrix}
J_{CC}(v,v_*) & J_{CS}(v,v_*) \\
J_{SC}(v,u) & J_{SS}(v,u)
\end{pmatrix}
=
\begin{pmatrix}
-b(v-v_*)^2 &
-b(v-v_*)^2-c_{CS} \\
-a(v-u)^2-c_{SC} &
-a(v-u)^2
\end{pmatrix}.
\]

Integrating the payoff against the opponent over the opponent's state distribution, the  expected payoffs for each strategy are:
\begin{align*}
\mathcal{P}_C[f](v, t) &= -b(v - m_v(t))^2 - b\Sigma_v(t)- c_{CS}(1 - m_s(t)), \\
\mathcal{P}_S[f](v, u, t) &= -a(v - u)^2 - c_{SC} m_s(t).
\end{align*} 
The corresponding payoff gain,
\(
\mathcal{G}[f]
:=\mathcal{P}_C[f]-\mathcal{P}_S[f],
\)
measuring the payoff advantage of $\pi_C$ over $\pi_S$, is
\begin{equation}\label{eq:gain-timedep}
\mathcal{G}[f](v,u,t)
=
a(v-u)^2
-b(v-m_v(t))^2
+K(t),
\end{equation}
where we have introduced
\begin{equation*}
K(t)
:=
c_{SC}m_s(t)-c_{CS}(1-m_s(t))
-b \Sigma_v(t).
\end{equation*}

The gain \eqref{eq:gain-timedep} reflects the competition between  two behavioral mechanisms; the term $a(v-u)^2$ favors the compromise strategy when the public opinion departs from the private conviction, while the term $-b(v-m_v(t))^2$ favors stubbornness for opinions far from the population mean. 
The time-dependent shift $K(t)$ incorporates the influence of the mean strategy and the opinion variance, while the population mean $m_v(t)$ determines the position of the quadratic term in the gain.

\begin{remark}[Key difference from the constant-payoff case.] Unlike the constant-payoff case, the gain \eqref{eq:gain-timedep} now depends on both the individual state $(v,u)$ and the macroscopic state of the population through $m_v(t)$, $\Sigma_v(t)$ and $m_s(t)$. In particular, the strategy marginal $g(s,t)$ no longer satisfies an autonomous equation, and the exact push-forward representation of Theorem~\ref{thm:projective} is not available. The analysis therefore shifts from exact solvability to a qualitative characterization of  the geometry of the segregation boundary and the structure of stationary states. \end{remark}

\paragraph{Time-dependent segregation boundary.}

At each time $t$, the sign of $\mathcal{G}[f]$ partitions the opinion--conviction phase space $I^2$ into two evolving regions:
\begin{equation*}
  \mathcal{R}_C(t) = \{(v,u)\in I^2:\mathcal{G}[f](v,u,t)>0\},
  \qquad
  \mathcal{R}_S(t) = \{(v,u)\in I^2:\mathcal{G}[f](v,u,t)<0\}.
\end{equation*} The region $\mathcal{R}_C(t)$  identifies states for which the current payoff structure favors $\pi_C$, whereas $\mathcal{R}_S(t)$ identifies those for which $\pi_S$ is favored. These regions characterize the instantaneous payoff advantage of the two strategies and should not be interpreted as regions occupied exclusively by pure C or S agents. Indeed, at each time $t$ every agent updates their strategy based on the sign of \eqref{eq:gain-timedep}, so that the regions determine the \emph{direction} of the strategic drift for each agent, not the asymptotic composition of the population. 

The boundary defined by $\mathcal{G}[f](v,u,t)=0$  is the zero-level set of the quadratic form:
\begin{equation}\label{eq:hyperbola}
 a(v-u)^2 - b(v-m_v(t))^2 = -K(t), 
\end{equation} which describes a conic in the $(v,u)$-plane. Its shape depends on the sign of $K(t)$, and only its intersection with the admissible domain $I^2$ is physically relevant. In particular, for $K(t)\neq0$, the boundary is a non-degenerate hyperbola in the $(v,u)$-plane , whereas for $K(t)=0$ it degenerates into two intersecting lines. 

\subparagraph{Case $K(t)<0$: hyperbola opening along $v-u$.} 
The radicand $b(v-m_v(t))^2+|K(t)|>0$ is strictly positive for all $(v,u)\in I^2$, so the boundary is well-defined throughout the domain. The two regions are:
\[
\mathcal{R}_S(t)
=
\left\{
(v,u)\in I^2:
|v-u|
<
\sqrt{
{\left((b(v-m_v(t))^2+|K(t)|\right)}/{a}
}
\right\},
\]
\[
\mathcal{R}_C(t)
=
\left\{
(v,u)\in I^2:
|v-u|
>
\sqrt{
{\left(b(v-m_v(t))^2+|K(t)|\right)}/{a}
}
\right\}.
\]
Thus the stubbornness-favoring region $\mathcal R_S$ contains the entire diagonal $\{v=u\}$. The vertical separation between the two branches has minimum width, attained at $v=m_v(t)$, given by
\begin{equation*}
\Delta u_{\min}(t)
=
2\sqrt{|{K(t)}|/{a}}
.
\end{equation*}

\subparagraph{Case $K(t)>0$: hyperbola opening along $v-m_v(t)$.}  
For a fixed $v$, the boundary exists only where $|v-m_v(t)|\geq\sqrt{K(t)/b}$; elsewhere $\mathcal{G}>0$ for all $u\in I$, so the entire vertical section belongs to $\mathcal{R}_C(t)$. The two regions are
\begin{align*}
    \mathcal{R}_C(t)
    &=
    \{(v,u)\in I^2:\,b(v-m_v(t))^2-a(v-u)^2<K(t)\},
    \\
    \mathcal{R}_S(t)
    &=
    \{(v,u)\in I^2:\,b(v-m_v(t))^2-a(v-u)^2>K(t)\}.
\end{align*}
The compromise-favoring region contains a central strip around $v=m_v(t)$, within which all values of $u$ favor $\pi_C$. Along the diagonal $\{v=u\}$, compromise is favored whenever $|v-m_v(t)|<\sqrt{K(t)/b}$.

\subparagraph{Case $K(t)=0$: degenerate boundary.}
The quadratic form factors and the boundary degenerates into two intersecting lines
\begin{equation*}
    u=v\mp\sqrt{b/a}\,(v-m_v(t)),
\end{equation*}
which coincide with the asymptotes of the hyperbolas for $K(t)\neq0$. The ratio $b/a$ determines the orientation of the boundary, while $K(t)$ controls its position and scale.

\paragraph{Moment equations and non-closure.}
The state-dependent gain couples the strategic variable $s$ to the opinion--conviction variables $(v,u)$ at every order of the moment hierarchy, and no finite-dimensional closure, in general, is available.

The mean opinion equation~\eqref{eq:moment-v} retains its form, but the covariance $\Sigma_{su}(t)$ now evolves non-trivially. Even under factorized initial conditions $f_0=g_0(s)h_0(v,u)$, the state-dependent gain can induce correlations between strategy and conviction by selectively driving agents with different states toward different strategic choices.

The mean strategy evolves as 
\begin{equation*}
\tau_\sigma \dot m_s
=
a\int_\Omega s(1-s)(v-u)^2f\, \mathrm{d}\omega-b\int_\Omega s(1-s)(v-m_v)^2f\,\mathrm{d}\omega
+K(t)\bigl(m_s-m_s^{(2)}\bigr),
\end{equation*} where the first two integrals involve mixed moments of the joint distribution $(v,s,u)$.

Introducing the strategy--opinion covariance $\Sigma_{sv}(t)$, the opinion variance satisfies
\begin{equation*}
\tau_v\,\dot \Sigma_v
=
-2\Sigma_v
+2\,\Sigma_{vu}(t)
-2\int_\Omega
s(v-m_v)(u-m_u)\,f d \omega+2(m_v-m_u)\Sigma_{sv}
+\nu^2\int_\Omega D(v)^2\,f \mathrm{d}\omega,
\end{equation*}
and the strategy--conviction covariance evolves as
\begin{equation*}
\tau_\sigma
\frac{d}{dt}\Sigma_{su}(t)
=
\int_\Omega
(u-m_u)s(1-s)\mathcal{G}[f]\, f \mathrm{d}\omega,
\end{equation*}
where the time derivative of $(u-m_u)$ vanishes since $u$ is frozen. Substituting~\eqref{eq:gain-timedep} into each equation shows that the right-hand sides involve higher-order mixed moments of $(v,s,u)$, so the hierarchy does not close at any finite order. This contrasts with the constant-payoff case, where the strategy dynamics decouples from $(v,u)$ and the macroscopic description reduces to the exact scalar ODE~\eqref{eq:strategy-cumulative-gains}.

\paragraph{Candidate stationary states.}

We characterize zero-flux stationary states $f^\infty(v,s,u)$ by imposing the simultaneous vanishing of the opinion and strategy fluxes in the stationary equation~\eqref{eq:FP-two}:
\begin{equation*}
\left(
s(m_v-v)+(1-s)(u-v)
\right)f^{\infty}
-\frac{\nu^2}{2}
\partial_v\left(D(v)^2f^{\infty}\right)=0,
\qquad
s(1-s)\mathcal{G}[f^\infty](v,u)f^\infty(v,s,u)
=0.
\end{equation*}

The second condition implies that any stationary mass in the interior $s\in(0,1)$ must be supported on the zero-gain set $\{\mathcal{G}[f^\infty]=0\}$, while the mass in $s\in\{0,1\}$ is compatible with zero strategic flux since $s(1-s)=0$. This generates three potential classes of candidate stationary states.

To describe the stationary opinion profiles, we introduce the normalized zero-flux density
\begin{equation}\label{eq:gibbs-profile}
p(v;\lambda)
=
\frac{C(\lambda)}{D(v)^2}
\exp\!\left(
\frac{2}{\nu^2}
\int^v
\frac{\lambda-w}{D(w)^2}\,\mathrm{d}w
\right),
\end{equation}
where $C(\lambda)$ is the normalization constant and
$\lambda\in(-1,1)$. This density profile is normalized to $1$ in $I$ and gives $\int_I vp(v;\lambda)\mathrm{d}v=\lambda$, provided that $D(v)$ satisfies condition (D1) of Proposition~\ref{prop:domain_invariance}. 

\subparagraph{1. Pure states.} 
The entire population adopts a single strategy.  

\emph{Pure Compromise ($s^\infty=1$):} The distribution takes the form $f^\infty(v,s,u) = \rho_C^\infty(v,u)\delta(s-1)$, where $\rho_C^\infty(v,u)=p(v;{m_v^{\infty}})k(u)$. 
Since the stationary opinion profile $p(v;{m_v^{\infty}})$ has mean $m_v^{\infty}$, the self-consistency condition for the population mean is automatically satisfied. Thus, any $m_v^\infty\in(-1,1)$ is admissible, and the pure compromise state admits a one-parameter family of
stationary profiles.
    
\emph{Pure Stubborn ($s^\infty=0$):} The distribution is $f^\infty(v,s,u) = \rho_S^\infty(v,u)\delta(s)$, where $\rho_S^\infty(v,u) =p(v;u)k(u)$. Since the conditional stationary opinion profile $p(v;u)$ has mean $u$, the global mean opinion satisfies
\(m_v^\infty=m_u\).

Each pure state is an exact stationary state for any sign of $\mathcal{G}[f^\infty]$, since $s(1-s)$ vanishes at $s\in\{0,1\}$. The sign of $\mathcal{G}[f^\infty]$ determines the direction of the strategic drift near each pure state. The pure compromise state is stable with respect to strategic perturbations if $\mathcal{G}[f^\infty]>0$ a.e.\ on the support of $\rho_C^\infty$, and can be unstable otherwise. Similarly, the pure stubborn state is stable with respect to strategic perturbations if $\mathcal{G}[f^\infty]<0$ a.e.\ on the support of $\rho_S^\infty$, and can be unstable otherwise.

\subparagraph{2. Internal coexistence.} A stationary distribution with positive absolutely continuous mass for $s\in(0,1)$ must have its support contained in $\{\mathcal{G}[f^\infty]=0\}$. For the present payoff specification, this set is generically a one-dimensional curve in the two-dimensional $(v,u)$-space and cannot carry positive Lebesgue measure. Interior coexistence within the class of regular densities is therefore excluded, but this argument does not exclude singular stationary measures supported on the zero-gain set.

\subparagraph{3. Segregated states.} Within the class of regular stationary states with no persistent
interior strategic mass, the stationary distribution may concentrate
on the two boundaries of the strategic simplex and take the form
\begin{equation}\label{eq:stationary-mixed}
f^\infty(v,s,u)
=
\rho_C^\infty(v,u)\delta(s-1)
+
\rho_S^\infty(v,u)\delta(s),
\end{equation}
Using the stationary profiles defined in \eqref{eq:gibbs-profile}, the two components can be written
as
\[
\rho_C^\infty(v,u)
=p(v;{m_v^\infty})
k_C(u),
\qquad
\rho_S^\infty(v,u)
=
p(v;u)k_S(u),
\]
where $k_C(u)$ and $k_S(u)$ are the conviction marginals of the compromise and stubborn subpopulations, respectively.
Since the conviction variable is conserved, they satisfy
\(
k_C(u)+k_S(u)=k(u).
\)

Defining the masses of the two subpopulations by
\[
M_C:=\int_I k_C(u)\,\mathrm{d}u,
\qquad
M_S:=\int_I k_S(u)\,\mathrm{d}u,
\qquad
M_C+M_S=1,
\] the mean strategy is then $m_s^\infty=M_C$, while the
self-consistency condition for the mean opinion reads
\begin{equation*}
m_v^\infty
=
M_Cm_v^\infty+\int_I u\,k_S(u)\,\mathrm{d}u.
\end{equation*}
In particular, if $M_S>0$, this condition gives
\[
m_v^\infty
=
\frac{1}{M_S}\int_I u\,k_S(u)\,\mathrm{d}u,
\]
so that the global mean opinion coincides with the mean conviction of the stubborn subpopulation.

The remaining macroscopic quantities must satisfy
\[
\Sigma_v^\infty
=
\int_{I^2}(v-m_v^\infty)^2
\bigl(\rho_C^\infty(v,u)+\rho_S^\infty(v,u)\bigr)
\,\mathrm{d}v\,\mathrm{d}u
\]
and
\(
K^\infty
=
c_{SC}M_C-c_{CS}M_S-b\Sigma_v^\infty.
\)
Together with $k_C(u)+k_S(u)=k(u)$, these relations provide the self-consistency conditions for the segregated stationary
state. 

The sign of the stationary gain determines the direction of the strategic drift near each boundary state. In particular, the pointwise conditions
\[
\mathcal{G}[f^\infty](v,u)>0
\quad\text{on }\operatorname{supp}\rho_C^\infty,
\qquad
\mathcal{G}[f^\infty](v,u)<0
\quad\text{on }\operatorname{supp}\rho_S^\infty
\]
ensure that the strategic drift points toward the corresponding pure strategy on each occupied subpopulation. Whether these conditions can be  simultaneously satisfied depends on the geometry of the stationary segregation boundary $\{\mathcal{G}[f^\infty]=0\}$ relative to the supports of the two subpopulations. They provide pointwise strategic invasion criteria, rather than a complete stability characterization of the coupled system.

\begin{remark}[Factorization structure at stationarity]
Unlike the constant-payoff case, the asymptotic strategy--conviction marginal is not known a priori through an exact push-forward representation. It must instead be determined together with the stationary state. Within the class of regular zero-flux stationary states considered here, the absence of positive absolutely continuous mass in the interior of the strategy interval leads to a representation in terms of the two boundary components $\rho_C^\infty(v,u)$ and $\rho_S^\infty(v,u)$. Their conviction marginals need not coincide separately with the conserved population marginal $k(u)$, but their sum must equal $k(u)$.
\end{remark}

\begin{remark}[Role of the segregation boundary at stationarity]
The zero-gain set $\{\mathcal{G}[f^\infty]=0\}$ is not a support constraint for the components $\rho_C^\infty$ and $\rho_S^\infty$ of a segregated state~\eqref{eq:stationary-mixed}; since the strategic flux vanishes identically at $s\in\{0,1\}$, both components can have support throughout $I^2$. Rather, the sign of $\mathcal{G}[f^\infty]$ determines whether each pure-strategy state is stable to strategic perturbations, as described above. The segregation boundary thus characterizes the direction of strategic drift rather than the support of the stationary distribution.
\end{remark}

\subsection{Three-strategy model with constant payoffs}\label{sec:three-strategy}

We extend the analysis to the three-strategy setting \( \Pi=\{\pi_1,\pi_2,\pi_3\}, \) identifying them  with the compromise, moderate, and stubborn behaviors, respectively, with compromise parameters $\theta=(1,\tfrac12,0)$, consistently with the ordering \eqref{theta_order}.
We write $\sigma=(s,r,q)$, with $q=1-s-r$ determined by the normalization constraint, and represent the strategic simplex in reduced coordinates as \[\Delta_2 = \{(s,r)\in\mathbb{R}^2_+: s+r\leq 1\},\] where $s$ and $r$ denote the frequencies of  $\pi_1$ and  $\pi_2$, respectively.
The effective compromise propensity is $\Theta(s,r)=s+r/2$, and the joint density $f=f(v,s,r,u,t)$  is defined on $\Omega = I\times\Delta_2\times I$.

The VFP equation~\eqref{eq:FP_gain} specializes to
\begin{align*}
\partial_t f = &-\frac{1}{\tau_v}\partial_v \left( \left( \left(s+\frac r2\right)(m_v - v) + \left(1-s-\frac r2\right)(u - v) \right) f - \frac{\nu^2}{2} \partial_v \left( D(v)^2 f \right) \right) \nonumber\\ &- \frac{1}{\tau_\sigma} \left(\partial_s \left( s\left((1-s)\,\mathcal{G}_{1}[f]-r\,\mathcal{G}_{2}[f]\right) f  \right) +\partial_r \left(r\left((1-r)\,\mathcal{G}_{2}[f]-s\,\mathcal{G}_{1}[f]\right)f\right)\right).
\end{align*}
where the relative payoff gains of $\pi_1$ and $\pi_2$, with respect to the reference strategy $\pi_3$ are,
$$\mathcal{G}_{1}[f]:=\mathcal{P}_1[f]-\mathcal{P}_3[f], \qquad \mathcal{G}_{2}[f]:=\mathcal{P}_2[f]-\mathcal{P}_3[f].$$

\paragraph{Constant payoff specification.}
Let $J=(J_{ij})_{i,j=1}^3\in\mathbb{R}^{3\times3}$ be a constant payoff matrix. Since the payoff is independent of the agents' states, the gains depend on $f$ only through the mean strategy frequencies
\[
m_s(t)
=
\int_{\Omega} s\,f(v,s,r,u,t)\,\mathrm{d}\omega,
\qquad
m_r(t)
=
\int_{\Omega} r\,f(v,s,r,u,t)\,\mathrm{d}\omega.
\]
Taking $\pi_3$ as the reference strategy, the gains reduce to the affine functions \begin{align*}
\mathcal{G}_{1}[f]
& = a_{13} + m_s(t)(a_{11} - a_{13}) + m_r(t)(a_{12} - a_{13}) = G_1(m_s(t),m_r(t)),\\
\mathcal{G}_{2}[f] 
& =  a_{23} + m_s(t)(a_{21} - a_{23}) + m_r(t)(a_{22} - a_{23}) = G_2(m_s(t),m_r(t)).
\end{align*}
which are independent of $(v,u)$.

\paragraph{Exact strategy dynamics.}
Since $G_1$ and $G_2$ are independent of $(v,u)$, Theorem~\ref{thm:projective} applies directly. The strategy marginal $$g(s,r,t):=\iint_{I^2}f(v,s,r,u,t)\,\mathrm{d}v\,\mathrm{d}u$$ satisfies the autonomous equation
\begin{equation*}
  \partial_t g
  +
  \frac{1}{\tau_\sigma} \left(
  \partial_s(s((1-s){G}_{1}-r{G}_{2})g)
  +\partial_r(r((1-r){G}_{2}-s{G}_{1})g)
  \right)
  =0,
\end{equation*}
and its unique solution is the push-forward
\begin{equation}\label{eq:g-pushforward-3s}
    g(\cdot,\cdot,t)
    =
    \bigl(\Phi_{\mathcal{T}(t)}\bigr)_\#g_0,
    \qquad
    (m_s(t),m_r(t))
    =
    \mathcal{E}(\mathcal{T}(t)),
\end{equation}
where the projective map $\Phi_{\mathcal{T}}:\Delta_2\to\Delta_2$ and cumulative gains $\mathcal{T}(t)=(\mathcal{T}_1(t), \mathcal{T}_2(t))^\top$ are given by~\eqref{eq:strategy-projective-map} and~\eqref{cumulative_gains}. For absolutely continuous $g_0$, the density reads
\begin{equation*}
    g(s,r,t)
    =
    g_0\!\bigl(\Phi_{-\mathcal{T}(t)}(s,r)\bigr)
    \cdot
    \frac{\mathrm{e}^{-\mathcal{T}_1-\mathcal{T}_2}}
    {\bigl(1-s-r+\mathrm{e}^{-\mathcal{T}_1}s+\mathrm{e}^{-\mathcal{T}_2}r\bigr)^3}.
\end{equation*} By the same argument as
Corollary~\ref{cor:eta-pushforward}, the strategy--conviction marginal $\eta(s,r,u,t)$ satisfies
\begin{equation}\label{eq:eta-three}
    \eta(\cdot,\cdot,\cdot,t)
    =
    \bigl(\Phi_{\mathcal{T}(t)},\mathrm{Id}\bigr)_\#\eta_0,
\end{equation}
so that the covariances $\Sigma_{su}(t)$ and $\Sigma_{ru}(t)$ are determined exactly at all times from $\eta_0$ and $\mathcal{T}(t)$.

\paragraph{Macroscopic moment equations.}
The evolution of the first-order macroscopic moments is given by:
\begin{align}
    \tau_\sigma\,\dot{m}_s(t)
    &=
    G_1(m_s,m_r)\bigl(m_s-m_s^{(2)}\bigr)
    -G_2(m_s,m_r)\,m_{sr},
    \label{eq:ms-dot-three}\\
    \tau_\sigma\,\dot{m}_r(t)
    &=
    G_2(m_s,m_r)\bigl(m_r-m_r^{(2)}\bigr)
    -G_1(m_s,m_r)\,m_{sr},
    \label{eq:mr-dot-three}\\
    \tau_v\,\dot{m}_v(t)
    &=
    \bigl(1-m_\Theta(t)\bigr)\bigl(m_u-m_v(t)\bigr)
    -\Sigma_{su}(t)-\tfrac{1}{2}\Sigma_{ru}(t),
    \label{eq:mv-dot-three}
\end{align}
where  
\[
m_s^{(2)}(t)
:=
\int_{\Omega} s^2 f\,\mathrm{d}\omega,
\qquad
m_r^{(2)}(t)
:=
\int_{\Omega} r^2 f\,\mathrm{d}\omega,
\qquad
m_{sr}(t)
:=
\int_{\Omega} sr\,f\,\mathrm{d}\omega,
\]
are higher-order strategy moments, $m_\Theta(t):=m_s(t)+m_r(t)/2$ is the
mean effective compromise level, and
$\Sigma_{su}(t)$,
$\Sigma_{ru}(t)$ are the
strategy--conviction covariances. 

Although $m_s$ and $m_r$ are determined exactly by the push-forward representation~\eqref{eq:g-pushforward-3s}, equations~\eqref{eq:ms-dot-three}--\eqref{eq:mr-dot-three} show that the moment dynamics involves higher-order strategy moments; however, these are consistently determined by the push-forward and are displayed here to make the coupling explicit. Equation~\eqref{eq:mv-dot-three} is not closed at the level of first-order moments, since it depends on $\Sigma_{su}(t)$ and $\Sigma_{ru}(t)$, themselves determined by $\eta_0$ and $\mathcal{T}(t)$ through~\eqref{eq:eta-three}.

\paragraph{Overview of asymptotic  strategy dynamics.}

The long-time behavior of $g$ is governed by the asymptotic behavior of the relative strategy weights
\[
\mathrm{e}^{\mathcal T_1(t)}:\mathrm{e}^{\mathcal T_2(t)}:1.
\]
When these  converge, the limiting distribution is obtained by applying the corresponding limiting projective map to $g_0$. If they do not converge, the push-forward $(\Phi_{\mathcal T(t)})_\#g_0$ need not converge in the weak-$*$ topology of probability measures.

Unlike the two-strategy case, the three-strategy setting admits a substantially richer variety of dynamical scenarios. Complete topological classifications of three-strategy replicator dynamics, comprising 49 qualitatively distinct cases up to reversal of the flow, are available in~\cite{Bomze1995}. Rather than reproducing this classification, we focus on some representative scenarios that are particularly  relevant for the kinetic model, illustrated each with a canonical example and numerically investigated
in Section~\ref{sec:test5}.

\subparagraph{Single-strategy dominance.}
If one strategy has a strictly positive net gain over all opponents, its cumulative gain diverges and the strategy distribution concentrates on the corresponding vertex of $\Delta_2$. This is the direct analogue of the dominance regimes of Section~\ref{sec:2s-constant}: the convergence rate is determined by the value of $G_i$ at the attracting vertex, analogously to Remark~\ref{rem:convergence-rates}.

\subparagraph{Two-strategy reduction.}
If one strategy becomes asymptotically negligible relative to the other two, the distribution concentrates on a face of $\Delta_2$. The three faces correspond to three distinct two-strategy subgames: the face $\{q=0\}$ involves $\pi_1$ and $\pi_2$; the face $\{r=0\}$ involves $\pi_1$ and $\pi_3$; the face $\{s=0\}$ involves $\pi_2$ and $\pi_3$. Each subgame may exhibit any of the four regimes of
Section~\ref{sec:2s-constant}.

\subparagraph{Three-strategy persistence.}
When none of the three strategies becomes asymptotically negligible, all three strategies remain involved in the long-time dynamics. We consider here three genuinely three-strategy phenomena which have no analogue in the one-dimensional two-strategy setting.

\emph{Interior coexistence (convergent case).}
If $\mathcal T(t),$ converges to a finite limit $\mathcal T^\infty\in\mathbb{R}^2$, all three strategies persist and the distribution approaches $g^\infty=(\Phi_{\mathcal T^\infty})_\#g_0$. The mean strategy converges to the interior equilibrium $(\hat m_s,\hat m_r)\in{\rm int}(\Delta_2)$ satisfying the equal-payoff conditions \begin{equation}\label{G120}G_1(\widehat m_s,\widehat m_r)=0,\qquad G_2(\widehat m_s,\widehat m_r)=0,\end{equation}  whose explicit expression is given by \[
(\hat m_s,\hat m_r) =
\frac{
\left(
a_{12}a_{23}-a_{13}a_{22},
\;
a_{13}a_{21}-a_{11}a_{23}
\right)
}{
(a_{11}-a_{13})(a_{22}-a_{23})
-
(a_{12}-a_{13})(a_{21}-a_{23})
},
\] provided that the denominator is non-zero. Notably, this equilibrium is interior if and only if $\hat m_s>0$, $\hat m_r>0$, and $\hat m_s+\hat m_r<1$. The limiting profile $g^\infty$  depends on $g_0$ through the  projective map $\Phi_{\mathcal T^\infty}$; two different initial distributions sharing the same mean $(\hat m_s,\hat m_r)$ will in general produce different limiting profiles. The convergence rate is determined by the eigenvalues of the Jacobian of $(G_1\circ\mathcal{E},G_2\circ\mathcal{E})$ at  $\mathcal{T}_\infty$, which depend on $g_0$ through $\mathcal{E}'(\mathcal{T}_\infty)$.

\emph{Persistent oscillations (rock-paper-scissors).} If the weights oscillate without converging, no stationary distribution exists. A sufficient condition for closed orbits is that the payoff matrix be zero-sum ($J+J^\top=0$), which implies that the trace of the Jacobian at the interior equilibrium vanishes and the dynamics is conservative. In this case, $g(s,r,t)$ exhibits persistent oscillatory behavior, the mean $m_\Theta(t)$ oscillates around the interior equilibrium value $\hat m_\Theta=\hat m_s+\hat m_r/2<1$, and since its time average is strictly less than~1, we have $\mathcal{I}(t)\to+\infty$ linearly in $t$.

\emph{Heteroclinic behavior.} If the interior equilibrium is unstable and each vertex is unstable in exactly one direction, the strategic composition  mean spiral toward $\partial\Delta_2$ and mass alternates among the three vertices with periods of increasing duration, approaching a heteroclinic cycle. Large-amplitude intermittent oscillations of $m_\Theta(t)$ result; the convergence of $\mathcal{I}(t)$ depends on the rate at which the dominance phases lengthen.

\begin{remark}[Kinetic vs.~classical replicator dynamics]The qualitative regime is determined by the structure of the relative gains on the simplex, including the location of their common zeros and their sign pattern on the interior and boundary. When an interior equilibrium is present, its local stability can be characterized through the Jacobian of the associated classical replicator dynamics. For the kinetic model, the evolution of $m_s$ and $m_r$ involves higher-order moments through equations~\eqref{eq:ms-dot-three}--\eqref{eq:mr-dot-three}. Thus, while the equilibrium locations coincide with those of the classical replicator dynamics, the convergence rate and the shape of the limiting distribution $g^\infty$ depend on the initial distribution $g_0$.
\end{remark}

\paragraph{Asymptotic mean opinion.}

Equation~\eqref{eq:mv-dot-three} can be integrated exactly by introducing the accumulated effective compromise $\mathcal I(t):=\int_0^t(1-m_\Theta(\tau))\,d\tau$. This yields
\begin{equation}\label{eq:mv-exact-3s}
    m_v(t)
    =
    m_u
    +\bigl(m_v(0)-m_u\bigr)\mathrm{e}^{-\mathcal I(t)/\tau_v}
    -\frac{\mathrm{e}^{-\mathcal I(t)/\tau_v}}{\tau_v}
    \int_0^t
    \!\Bigl(\Sigma_{su}(\tau)+\tfrac{1}{2}\Sigma_{ru}(\tau)
    \Bigr)
    \mathrm{e}^{\mathcal I(\tau)/\tau_v}\,\mathrm{d}\tau.
\end{equation}
The asymptotic behavior is therefore governed not only by the limiting strategic composition, but also by $\mathcal I(t)$ and by the strategy--conviction covariance term. 

The asymptotic behavior of $m_v(t)$ depends on three qualitatively distinct scenarios for $m_\Theta(t)$.

\subparagraph{Convergent regimes with $m_\Theta^\infty<1$.} This covers dominance of $\pi_2$ and $\pi_3$, as well as convergent interior coexistence regimes. In this case
\(
\mathcal I(t)\sim(1-m_\Theta^\infty)t\to+\infty
\). If
\(\Sigma_{su} (t)\to\Sigma_{su}^\infty\), and \(\Sigma_{ru}(t)\to\Sigma_{ru}^\infty\), equation~\eqref{eq:mv-exact-3s} gives
\begin{equation*}
    m_v^\infty
    =
    m_u
    -
\frac{\Sigma_{su}^\infty+\frac12\Sigma_{ru}^\infty}
    {1-m_\Theta^\infty}.
\end{equation*}
Under the factorization assumption($\Sigma_{su}\equiv\Sigma_{ru}\equiv0$), this reduces to $m_v^\infty=m_u$ for all such regimes.  

\subparagraph{Convergence to $m_\Theta^\infty=1$.}
In this case, the behavior of $m_v$ depends on the \emph{rate} at which
$1-m_\Theta(t)$ vanishes.  
In particular if, 
\[
\mathcal I_\infty:=\int_0^\infty(1-m_\Theta(\tau))\,d\tau<+\infty,
\]
and the covariance term $\Sigma_{su}+\frac12\Sigma_{ru}$ is integrable on $[0,+\infty)$, then
\begin{equation}\label{eq:mv-infty-comp-3s}
    m_v^\infty
    =
    m_u
    +\bigl(m_v(0)-m_u\bigr)\mathrm{e}^{-\mathcal I_\infty/\tau_v}
    -\frac{\mathrm{e}^{-\mathcal I_\infty/\tau_v}}{\tau_v}
    \int_0^\infty
    \!\Bigl(\Sigma_{su}(\tau)+\tfrac12\Sigma_{ru}(\tau)\Bigr)
    \mathrm{e}^{\mathcal I(\tau)/\tau_v}\,d\tau.
\end{equation}

This situation occurs, in particular, when $\pi_1$ dominates and $1-m_\Theta(t)$ is integrable. The initial-condition contribution is then not fully damped, and in general $m_v^\infty\neq m_u$. This is the three-strategy analogue of Regime~1 in Section~\ref{sec:2s-constant}: in both cases $\mathcal{I}_\infty<+\infty$ and the initial-condition contribution to $m_v(t)$ is not fully damped, so the asymptotic mean opinion retains memory of $m_v(0)$. 
The convergence $m_\Theta(t)\to1$ alone does not guarantee that $\mathcal I_\infty$ is finite. 
If instead $\mathcal I(t)\to+\infty$, the initial-condition contribution vanishes, and the limit of $m_v(t)$ depends on the asymptotic behavior of the covariance term.

\subparagraph{Persistent oscillations of $m_\Theta(t)$.}
Suppose that $m_\Theta(t)$ does not converge and has a long-time average
\[
\overline m_\Theta
:=
\lim_{T\to\infty}\frac1T\int_0^T m_\Theta(t)\,\mathrm{d}t<1.
\] Then
\(
\mathcal I(t)\sim(1-\overline m_\Theta)t
\to+\infty,
\)
so the initial-condition contribution in
\eqref{eq:mv-exact-3s} vanishes.

For conservative oscillations along a periodic orbit, the time average of the strategy frequencies coincides with the interior equilibrium, provided the corresponding interior equilibrium is unique and non-degenerate. In this case, $\overline m_\Theta=\hat m_s+\hat m_r/2<1$.

Under the factorization assumption, the covariance term is identically zero and $m_v(t)\to m_u$, despite the persistent oscillations of the strategic composition. Without factorization, the time-dependent covariances act as a forcing in the mean-opinion equation and may induce persistent oscillations of $m_v(t)$. In this case, convergence to a constant mean opinion cannot be concluded without additional assumptions on the long-time behavior of the covariance term. More specifically, if both $m_\Theta(t)$ and $\Sigma_{su}(t)+\frac12\Sigma_{ru}(t)$ are periodic, the mean-opinion equation has an attracting periodic solution whenever the time average of $1-m_\Theta(t)$ is positive. Consequently, $m_v(t)$ approaches this periodic solution,  which need not oscillate around $m_u$.

\paragraph{Asymptotic opinion distribution.}

Whenever the strategy marginal  converges to  a stationary profile $g^\infty(s,r)$, the argument of Section~\ref{sec:exact-analysis} applies. The conditional  opinion density $\rho(v\mid s,r,u,t)$ satisfies along the projective characteristic a one-dimensional non-autonomous Fokker--Planck equation whose drift converges to $$V_{\infty}(s,r,u)=\Theta(s,r)m_v^\infty+(1-\Theta(s,r))u= \left(s+\tfrac{r}{2}\right)m_v^\infty+\left(1-s-\tfrac{r}{2}\right)u$$ The conditional density
converges to the Gibbs density
\begin{equation*}
    \rho^\infty(v\mid s,r,u)
    =
    \frac{C(s,r,u)}{D(v)^2}
    \exp\left(
        \frac{2}{\nu^2}
        \int_0^v
        \frac{V_{\infty}(s,r,u)-z}{D(z)^2}\,\mathrm{d}z
    \right),
\end{equation*} and the asymptotic joint distribution factorizes as $f^\infty=\rho^\infty(v|s,r,u)\,\eta^\infty(s,r,u)$. The observable marginal opinion profile for agents with conviction $u$ is:
\begin{equation*}
\rho^\infty(v\mid u)
=
\iint_{\Delta_2}
\rho^\infty(v\mid s,r,u)\,
\eta^\infty(\mathrm{d}s,\mathrm{d}r\mid u).
\end{equation*}
The structure of the attraction center $V_{\infty}(s,r,u)$ depends on the asymptotic strategic configuration:

\emph{$\pi_1$ dominates
($g^\infty=\delta(s-1)\delta(r)$):}
$\Theta(1,0)=1$, so $V_{\infty}(1,0,u)=m_v^\infty$, independent of
$u$. All agents share the same Gibbs profile regardless of
conviction, with $m_v^\infty$ given
by~\eqref{eq:mv-infty-comp-3s}.

\emph{$\pi_2$ dominates
($g^\infty=\delta(s)\delta(r-1)$):}
$\Theta(0,1)=1/2$, so
$V_{\infty}(0,1,u)=\frac{1}{2}m_v^\infty+\frac{1}{2}u
=\frac{1}{2}(m_u+u)$, where the last equality uses
$m_v^\infty=m_u$ (since $m_\Theta^\infty=1/2<1$). The
attraction center depends linearly on $u$ with slope $1/2$.

\emph{$\pi_3$ dominates
($g^\infty=\delta(s)\delta(r)$):}
$\Theta(0,0)=0$, so $V_{\infty}(0,0,u)=u$. Each agent's opinion
is attracted to their own conviction.

\emph{Interior coexistence ($g^\infty$ non-degenerate):}
the attraction center varies with $(s,r)$, according to
\[
V_{\infty}(s,r,u)
=
\Theta(s,r)m_v^\infty+(1-\Theta(s,r))u.
\]
For fixed $u$, it is a convex combination of $u$ and $m_v^\infty$,
and therefore
\[    V_{\infty}(s,r,u)
    \in
    \bigl[\min\{u,m_v^\infty\},\max\{u,m_v^\infty\}\bigr],
\]
and $\rho^\infty(v\mid u)$ is a superposition of Gibbs densities
with heterogeneous centers, weighted by $\eta^\infty(ds,dr\mid u)$.

\begin{remark}[Three-strategy vs.\ two-strategy profiles]
The three pure-strategy limits produce three qualitatively distinct attraction patterns: strategy-independent ($V_{\infty}=m_v^\infty$), slope-$1/2$ affine in conviction ($V_{\infty}=\frac12(m_u+u)$), and identity ($V_{\infty}=u$). These interpolate between the two patterns available in the two-strategy case ($V_{\infty}=m_v^{\rm comp}$ for $\pi_C$ dominance and $V_{\infty}=u$ for $\pi_S$ dominance). In coexistence regimes, the mixture weight $g^\infty(s,r)$ combines all three patterns simultaneously, producing a strictly richer family of profiles than the two-strategy coexistence regime, where only the two endpoint patterns are mixed.
\end{remark}

\section{Numerical experiments}\label{sec:numerics}

We present a series of numerical tests to investigate the dynamics of the kinetic model and to assess the theoretical predictions derived in the previous sections. 

The first two tests focus on the asymptotic regimes characterized by a separation of time scales, as in Section~\ref{sec:asymptotic}. Test~1  considers the fast-opinion regime, investigating the relaxation of the conditional opinion distributions towards their Gibbs quasi-equilibria and the validity of the predicted relaxation laws, up to terms of order $\epsilon_v$. Test~2 considers the fast-strategy regime, comparing the moments of the particle system with those predicted by the reduced moment system to assess the accuracy of  the closure and its limitations in describing the opinion variance.

The remaining tests explore different aspects of the coupled opinion--strategy dynamics for comparable time scales, as highlighted in Section~\eqref{sec:exact-analysis}. Test~3 examines the strategy dynamics across the four asymptotic regimes of the two-strategy model with constant payoffs. Test~4 investigates the long-time behavior of the mean opinion and the role of initial strategy--conviction correlations, while Test~5 compares the asymptotic opinion--conviction profiles with their theoretical predictions. Test~6 focuses on the effects of state-dependent payoffs on strategic segregation and on the geometry of the segregation boundary. Finally, Test~7 explores the richer asymptotic dynamics of the three-strategy model.

For Test~3, the strategy marginal is computed using the exact push-forward representation~\eqref{eq:g-pushforward-two}, requiring only the numerical solution of the scalar ODE~\eqref{eq:cumulative-gain-equation-two}, from which the
strategy distribution is reconstructed through the associated projective flow. Tests~1--2 and 4--7 are instead performed using a direct simulation Monte Carlo (DSMC) particle scheme, described below, to simulate the coupled opinion--strategy--conviction dynamics.

\paragraph{Numerical scheme.}
To simulate the evolution of the VFP equation~\eqref{eq:FP_gain}  we employ an asymptotic Monte Carlo stochastic particle method of Nanbu type, based on random pairings  with small increments consistent with the quasi-invariant scaling \eqref{eq:grazing}; see \cite{pareschi2013interacting} for further details.

The distribution $f(v,\widehat \sigma,u,t)$ is approximated by an ensemble of $N$ particles, sampled from the initial distribution $f_0(v,\widehat \sigma,u)$. At the discrete time $t^n=n\Delta t$ the state of particle $i$ is denoted by $(v_i^n,\widehat \sigma_i^n,u_i)$, where $i=1,\ldots,N$, and $\Delta t>0$ is the fixed time increment. At each timestep, the particles are randomly permuted and divided into $N/2$ pairs, without repetition. Both particles in each pair are updated simultaneously using an explicit time discretization of  opinion and strategy dynamics, evaluated at the
current time level.

\emph{Opinion update.}  For a particle $i$ paired with particle $i_*$, the opinion variable is updated according  to the explicit discretization of~\eqref{eq:binary-opinion}, which reads
\begin{equation}\label{scheme_v}
v_i^{n+1}
=v_i^n
+\frac{\Delta t}{\tau_v}\left(
\Theta\left(\sigma_i^n\right)\left(v_{i_*}^n-v_i^n\right)
+\left(1-\Theta\left(\sigma_i^n\right)\right)\left(u_i-v_i^n\right)
\right)+\sqrt{\frac{\Delta t}{\tau_v}}\,\nu D(v_i^n)\xi_i^n,
\end{equation}
where $u_i$ remains constant in time and the stochastic increment $\xi_i^n$ is sampled according to the noise prescription of Proposition~\ref{prop:domain_invariance}, specifically  $\xi_i^n\sim\mathrm{Unif}(-\xi_{\max},\xi_{\max})$ with $\xi_{\max}=(1-\Delta t/\tau_v)/C_D$. 

\emph{Strategy update.}  The $K-1$ reduced strategy coordinates $\widehat \sigma^n_i=(\sigma_{1,i}^n,\ldots,\sigma^n_{K-1,i})$  are advanced by an explicit Euler discretization of the replicator equation~\eqref{eq:binary-game-gain}:
\begin{equation}\label{scheme_s}
\sigma_{k,i}^{n+1}
=
\sigma_{k,i}^{n}+\frac{\Delta t}{\tau_\sigma}\,\sigma_{k,i}^{n}
\left( 
G_k(v_i^n,u_i^n;v_{i_*}^n,u_{i_*}^n,\widehat{\sigma}^n_{i_*})
-
\displaystyle\sum_{\ell=1}^{K-1}
\sigma_{\ell,i}^{n}G_\ell(v_i^n,u_i^n;v_{i_*}^n,u_{i_*}^n,\widehat{\sigma}_{i_*}^n)
\right), \qquad k=1,\dots,K-1
\end{equation}
Here, the relative gains are evaluated at the current strategy of the partner particle. Domain preservation for the strategy update is guaranteed by the condition of Proposition~\ref{prop:domain_invariance}: $\Delta t/\tau_\sigma \leq1/(J_{\max}-J_{\min})$, satisfied for all payoff matrices used in the tests.

The same updates are applied to particle $i_*$, with the particle indices exchanged.

\emph{Parameters.} All particle-based simulations  in Tests~4--7 use
$\tau_v=\tau_\sigma=1$, $D(v)=1-v^2$, noise intensity $\nu=0.25$, and $N=10^6$ particles. The time step is set to $\Delta t=10^{-3}$. Test-specific initial distributions, payoff parameters, and final simulation times are specified in the corresponding tests whenever they vary between configurations. The parameters used in Tests~1--2 are  specified directly in the corresponding subsections.

\begin{remark}[Effective variance]
By Proposition~\ref{prop:domain_invariance} for $D(v)=1-v^2$ ($C_D=2$),  $\mathrm{Var}(\xi_i^n)=\xi_{\max}^2/3$,
the scheme discretizes the VFP equation with the effective intensity
\begin{equation*}
\nu_{\rm eff}:=\nu\,\frac{\xi_{\max}}{\sqrt3}\simeq\frac{\nu}{2\sqrt3},
\end{equation*} which is $\simeq0.072$ for $\nu=0.25$. All the theoretical predictions that depend on the noise intensity are evaluated with $\nu_{\rm eff}$; in Test~1 the noise has unit variance, so that $\nu_{\rm eff}=\nu$. For $D(v)=\sqrt{1-v^2}$ (Tests~1--2) Proposition~\ref{prop:domain_invariance} does not apply since condition (D2) is not satisfied ($C_D=+\infty$) and the condition $|v|\le1$ is checked numerically at every step.\end{remark}

\begin{remark}[Strategy update for small $\tau_\sigma$]
\label{rem:strategy-noise}
For constant payoffs the gain is an affine function of the opponent's strategy, and its mean-field value is $G(m_\sigma^n)$, with $m_\sigma^n$ the empirical mean strategy vector. In the binary scheme~\eqref{scheme_s}, evaluating the gain at a randomly sampled partner adds a sampling fluctuation to each strategy increment. For $K=2$, the conditional variance (up to finite-population corrections)
\[
    \operatorname{Var}(\Delta s_i\mid s_i)
    \simeq
    \frac{\Delta t^2}{\tau_\sigma^2}\,s_i^2(1-s_i)^2\operatorname{Var}_{g^n}(G(s_{i_*})),
\]
is proportional to $\Delta t/\tau_\sigma^2$ per unit time. The accumulated dispersion over a time $T$ is negligible only if $T\,\Delta t\,\operatorname{Var}(G)/\tau_\sigma^2\ll1$. This is satisfied for $\tau_\sigma=1$, but not for $\tau_\sigma=0.02$ and $\Delta t=10^{-3}$: in a preliminary simulation with $\Theta\equiv0.3$ all strategies accumulated at $\{0,1\}$, giving $\operatorname{Var}(s)\simeq0.21=\Theta(1-\Theta)$ instead of the  value $0.0676$ of the continuous replicator flow. In Test~2, where $\tau_\sigma\ll\tau_v$, the gain is therefore evaluated at the population mean for all particles, which gives the deterministic replicator update and preserves the projective evolution of the strategy distribution; the opinion update~\eqref{scheme_v} retains the random pairing of particles.
\end{remark}

\subsection{Test~1. Fast-opinion regime: quasi-equilibrium of the opinion distribution}

This test investigates the fast-opinion regime $\tau_v\ll\tau_\sigma$ of Section~\ref{sec:asymptotic}, for $K=2$ strategies and $D(v)=\sqrt{1-v^2}$. We consider a finite scale separation $\epsilon_v=2\cdot10^{-3}$ and examine whether the opinion distribution rapidly approaches its conditional Gibbs quasi-equilibrium and remains close to it while the strategy--conviction marginal evolves on the
slower time scale.

In the fast-opinion limit, the conditional opinion density given $(s,u)$ relaxes towards the Gibbs distribution \eqref{eq:rho-eq}, whose attraction center is given by \eqref{eq:local-target}. For $K=2$ and $\Theta(s)=s$ this center reads
\[
    \widehat V(s,u;\tau)=s\,m_v^{\rm eq}(\tau)+(1-s)u,
\]
where the self-consistent equilibrium mean \eqref{eq:fast-opinion-mean} takes the form
\begin{equation}\label{eq:test8-mveq}
    m_v^{\rm eq}(\tau)=m_u^{(0)}-\frac{\Sigma_{su}^{(0)}(\tau)}{1-m_s^{(0)}(\tau)} .
\end{equation}
Here $m_s^{(0)}$ and $\Sigma_{su}^{(0)}$ are the mean strategy and the strategy--conviction covariance of the slow
marginal $\eta^0(s,u,\tau)$, while $m_u^{(0)}$ is constant because convictions do not evolve. For constant payoffs,
which is the case considered here, both the full strategy--conviction marginal and $\eta^0$ are transported by the
same projective flow (Corollary~\ref{cor:eta-pushforward}), so that $m_s^{(0)}$ and $\Sigma_{su}^{(0)}$ coincide with
the values $m_s(t)$ and $\Sigma_{su}(t)$ measured on the particle system, with $\tau=t/\tau_\sigma$. In what follows we
therefore write $m_v^{\rm eq}(t)$ and omit the superscript $(0)$.

The test focuses on three aspects: the relaxation of the conditional opinion density towards the Gibbs distribution
on the fast time scale $\tau_v$; the agreement between the empirical mean opinion and the self-consistent value
$m_v^{\rm eq}$; and the persistence of this quasi-equilibrium as the strategy--conviction marginal evolves on the
slow time scale $\tau_\sigma$.

\paragraph{Setup.}We consider the following configuration.

\emph{Initial conditions.}
We consider two initial configurations of the strategy--conviction marginal, differing only in the dependence between
$s$ and $u$. In both cases opinions and convictions are initially independent, with $v\sim\mathrm{Unif}[-0.8,0.8]$ and
$u\sim\mathrm{Unif}[-0.6,0.6]$, so that $m_u=0$ and $\Sigma_u=0.12$.

\emph{Case~1: Initially uncorrelated strategy and conviction.}
The initial strategy is uniformly distributed on $\Delta_1$ and independent of $u$. Consequently $\Sigma_{su}(0)\simeq0$
and $m_v^{\rm eq}(0)\simeq0$. The initial joint distribution is given by
\[
f_0(v,s,u)
=
\frac{1}{1.92}
\mathbf{1}_{[-0.8,0.8]}(v)
\mathbf{1}_{\Delta_1}(s)
\mathbf{1}_{[-0.6,0.6]}(u).
\]

\emph{Case~2: Initially correlated strategy and conviction.}
The initial strategy is generated by
\[
    s=\min\Bigl\{\max\Bigl\{\frac{u+0.6}{1.2}+0.2\Bigl(\xi-\frac12\Bigr),\,0\Bigr\},\,1\Bigr\},
    \qquad \xi\sim\mathrm{Unif}[0,1],
\]
so that $m_s(0)\simeq0.5$ and $\Sigma_{su}(0)\simeq0.099$. The initial quasi-equilibrium mean is therefore
$m_v^{\rm eq}(0)\simeq-0.198$, despite $m_u=0$ and $m_v(0)\simeq0$. To characterize the resulting joint distribution, let
\[
z(u)=\frac{u+0.6}{1.2},
\qquad
p_0(u)=\frac{(0.1-z(u))_+}{0.2},
\qquad
p_1(u)=\frac{(z(u)-0.9)_+}{0.2},
\]
where $(x)_+=\max\{x,0\}$. The conditional distribution of
$s$ given $u$ is then
\[
f_{s|u}(s|u)
=
5\cdot \mathbf{1}_{(z(u)-0.1,z(u)+0.1)\cap(0,1)}(s)
+p_0(u)\delta_0(s)+p_1(u)\delta_1(s).
\]
Consequently, the initial joint distribution is
\[
f_0(v,s,u)
=
\frac{1}{1.92}
\mathbf{1}_{[-0.8,0.8]}(v)
\mathbf{1}_{[-0.6,0.6]}(u)\cdot f_{s|u}(s|u).
\] 

In both cases the payoff matrix and the payoff gain are
\[
    J=\begin{pmatrix}0&3\\2&0\end{pmatrix},\qquad G(m_s)=3-5m_s,
\]
and the equilibrium mean strategy is $\widehat m_s=0.6$. The strategy marginal thus evolves on the slow time scale
$\tau_\sigma$, while the conviction distribution remains fixed. The evolution of $m_s$ and $\Sigma_{su}$ induces a slow
variation of the self-consistent mean $m_v^{\rm eq}$; for the correlated initialization it changes from approximately
$-0.198$ at $t=0$ to $-0.221$ at $t=1$.

\emph{Parameters.}
The parameters are $\tau_v=0.002$, $\tau_\sigma=1$, $\nu=0.25$, $N=2\cdot10^6$, $\Delta t=10^{-5}$ (so that
$\alpha_v=\Delta t/\tau_v=0.005$) and $T=1$. 

\emph{Strategy update.}
The strategy update uses the strategy of a randomly selected partner. Since $\tau_\sigma=1$, the fluctuations
associated with this sampling are small, consistently with Remark~\ref{rem:strategy-noise}.

\paragraph{Predicted relaxation laws.}
For fixed $(s,u)$, the first two moments of the conditional density,
$M(t)=\mathbb E[v\mid s,u]$ and $Q(t)=\mathbb E[v^2\mid s,u]$, satisfy exactly (the drift is linear and the boundary
terms vanish by the zero-flux conditions)
\begin{equation}\label{eq:test8-moments}
    \tau_v\dot M=\lambda(t)-M,\qquad
    \tau_v\dot Q=2\lambda(t)M-(2+\nu^2)Q+\nu^2,\qquad
    \lambda(t)=s\,m_v(t)+(1-s)u .
\end{equation}
Moreover, the mean opinion satisfies $\tau_v\dot m_v=(1-m_s)(m_u-m_v)-\Sigma_{su}$. On the fast scale, where $m_s$ and
$\Sigma_{su}$ are essentially frozen, this gives
\begin{equation}\label{eq:test8-laws}
\begin{aligned}
    m_v(t)&=m_v^{\rm eq}+\bigl(m_v(0)-m_v^{\rm eq}\bigr)\mathrm{e}^{-at},\qquad a=\frac{1-m_s}{\tau_v},\\
    M(t)-\lambda_\infty&=e_M(0)\,\mathrm{e}^{-bt}+e_\lambda(0)\,\frac{b}{b-a}\bigl(\mathrm{e}^{-at}-\mathrm{e}^{-bt}\bigr),\qquad b=\frac1{\tau_v},
\end{aligned}
\end{equation}
where $\lambda_\infty=s\,m_v^{\rm eq}+(1-s)u$, $e_M(0)=M(0)-\lambda_\infty$ and
$e_\lambda(0)=s\,(m_v(0)-m_v^{\rm eq})$. When $m_v^{\rm eq}$ is constant, as in Case~1 ($m_v(0)$ and $m_v^{\rm eq}$ are approximately zero), $e_\lambda\equiv0$ and
$M(t)-\lambda_\infty=\mathrm{e}^{-t/\tau_v}(M(0)-\lambda_\infty)$. The higher moments of the conditional density relax faster:
for the Gibbs operator with $D=\sqrt{1-v^2}$ the $n$-th mode decays at rate $n\bigl(1+(n-1)\nu^2/2\bigr)/\tau_v$, which is
$1/\tau_v$ for $n=1$ (the mean) and at least $(2+\nu^2)/\tau_v$ for $n\ge2$. When $\lambda$ varies in time, as in Case~2, its evolution may
also induce transient contributions to the higher modes.
Nevertheless, after the initial relaxation, the leading
discrepancy is expected to be a residual shift $e_M(t)$ of
the conditional mean. For a small translation of a unimodal density, the corresponding
$L^1$ distance is approximately $2\rho_{\max}|e_M|$.
Using $\rho_{\max}\simeq2.3$ for the profiles considered here,
this provides an estimate of the contribution of the residual
mean shift to the observed $L^1$ discrepancy.

\paragraph{Diagnostics.}
The particle distribution is recorded at the times $t\in\{0,\tau_v,5\tau_v,10\tau_v,\tau_\sigma\}$. At each observation
time we compute the empirical moments $m_s(t)$ and $\Sigma_{su}(t)$ and use them in~\eqref{eq:test8-mveq} to evaluate
the instantaneous quasi-equilibrium mean $m_v^{\rm eq}(t)$.

To assess the conditional relaxation we select the particles in two rectangular windows of the $(s,u)$-plane,
$|s-\bar s|<0.07$ and $|u-\bar u|<0.07$. For Case~1 the window centers are $(\bar s,\bar u)=(0.4,-0.5)$ and $(0.75,0.3)$,
whereas for Case~2 they are $(0.083,-0.5)$ and $(0.75,0.3)$; the first window of Case~2 is chosen along the initial
strategy--conviction correlation, so that it contains a sufficiently large number of particles. We refer to the two
windows as window~1 and window~2. The windows contain $n\simeq3.3\cdot10^4$ particles in Case~1 and
$n\simeq1.6\cdot10^5$ in Case~2 up to $t=10\tau_v$, and between $3.0\cdot10^4$ and $3.8\cdot10^4$, and between
$1.2\cdot10^5$ and $1.3\cdot10^5$ at $t=\tau_\sigma$.

Since $(s,u)$ varies within each window, the theoretical conditional density is not a single Gibbs distribution.
It is obtained by averaging the Gibbs densities of the particles in the window: for each particle $i$ we compute
$\widehat V_i(t)=s_i\,m_v^{\rm eq}(t)+(1-s_i)u_i$ and evaluate the corresponding Gibbs density $\rho_{\widehat V_i}$
with parameter $\nu$, and the predicted density is the empirical mixture
\[
    \rho_{\rm mix}(v,t)=\frac1n\sum_{i=1}^n\rho_{\widehat V_i(t)}(v),
\]
where $n$ is the number of particles in the window. Each component is a rescaled Beta density on $I$
(Remark~\ref{rmk:beta_quasi_eq}), with mean and variance
\[
    \widehat V_i,\qquad \frac{\nu^2\bigl(1-\widehat V_i^2\bigr)}{2+\nu^2}.
\]
The variance of the mixture is therefore
\[
    \operatorname{Var}(\rho_{\rm mix})
    =\frac1n\sum_{i=1}^n\frac{\nu^2(1-\widehat V_i^2)}{2+\nu^2}+\operatorname{Var}_n(\widehat V_i),
\]
where $\operatorname{Var}_n$ is the empirical variance of the particle-dependent centers.

We compare the empirical histogram, computed with 30 bins, with this mixture through the $L^1$ distance between the
bin probabilities, which are computed exactly by integrating each Beta density over the bins, and through the ratio of the
variances. For efficiency, the bin probabilities of the mixture are computed from a random subsample of $3000$ particles of the window, while the means, the variances and the sampling level use all $n$ particles; the effect of the subsampling on the $L^1$ distance is expected to be well below the sampling level. To distinguish discrepancies due to the model from those due to the finite number of particles, we also report the sampling level
\[
    E_{\rm samp}=\sum_{i=1}^{30}\sqrt{\frac{2p_i}{\pi n}},
\]
where $p_i$ is the probability that the predicted mixture assigns to the $i$-th bin. This is the expected $L^1$
distance of the histogram from an exact prediction, with the factor $1-p_i$ neglected, so that it slightly
overestimates it. For the ratio of the variances, the sampling standard error is approximately $\sqrt{1.9/n}$, that is
$0.8\%$ in Case~1 and $0.35\%$ in Case~2, and for the conditional means $\sigma/\sqrt n$, that is about $10^{-3}$ in
Case~1 and $4\cdot10^{-4}$ in Case~2.

\begin{table}[ht]
\centering
\caption{Test~1. Comparison between the empirical opinion distributions and the corresponding Gibbs mixtures, for the
two initial conditions and the two $(s,u)$-windows.}
\label{tab:test8}

\begin{subtable}{\textwidth}
\centering
\caption{$L^1$ distance between the empirical conditional distribution and the Gibbs mixture at the observation times.
The last column is the range of the sampling level $E_{\rm samp}$ over the observation times.}
\label{tab:test8-transient}
\begin{tabular}{c|c|ccccc|c}
\hline
Initial condition & Window & $t=0$ & $t=\tau_v$ & $t=5\tau_v$ & $t=10\tau_v$ & $t=\tau_\sigma$ & $E_{\rm samp}$\\
\hline
Uncorrelated & 1 & 1.1010 & 0.4825 & 0.0203 & 0.0108 & 0.0180 & 0.0156--0.0163\\
             & 2 & 1.0881 & 0.3422 & 0.0167 & 0.0144 & 0.0150 & 0.0147--0.0158\\
\hline
Correlated   & 1 & 1.1751 & 0.7350 & 0.0224 & 0.0068 & 0.0083 & 0.0069--0.0079\\
             & 2 & 1.0773 & 0.4522 & 0.1106 & 0.0148 & 0.0123 & 0.0073--0.0078\\
\hline
\end{tabular}
\end{subtable}

\vspace{0.4cm}

\begin{subtable}{\textwidth}
\centering
\caption{Offset $M-\lambda_\infty$ of the conditional mean (empirical minus predicted by the Gibbs mixture) against the
law~\eqref{eq:test8-laws} evaluated at the window centers. Sampling standard error: $0.0009$ (Case~1), $0.0004$ (Case~2).}
\label{tab:test8-offset}
\begin{tabular}{c|c|rr|rr|rr}
\hline
 & & \multicolumn{2}{c|}{$t=\tau_v$} & \multicolumn{2}{c|}{$t=5\tau_v$} & \multicolumn{2}{c}{$t=10\tau_v$}\\
Initial condition & Window & DSMC & law & DSMC & law & DSMC & law\\
\hline
Uncorrelated & 1 & $\phantom{-}0.1089$ & $\phantom{-}0.1103$ & $\phantom{-}0.0024$ & $\phantom{-}0.0020$ & $-0.0003$ & $0.0000$\\
             & 2 & $-0.0257$ & $-0.0276$ & $-0.0005$ & $-0.0005$ & $\phantom{-}0.0017$ & $0.0000$\\
\hline
Correlated   & 1 & $\phantom{-}0.1820$ & $\phantom{-}0.1826$ & $\phantom{-}0.0057$ & $\phantom{-}0.0057$ & $-0.0002$ & $0.0002$\\
             & 2 & $\phantom{-}0.0987$ & $\phantom{-}0.0981$ & $\phantom{-}0.0233$ & $\phantom{-}0.0229$ & $\phantom{-}0.0027$ & $0.0020$\\
\hline
\end{tabular}
\end{subtable}

\vspace{0.4cm}

\begin{subtable}{\textwidth}
\centering
\caption{Conditional means, ratio of the variances and $L^1$ distance at $t=1$, with the sampling level.}
\label{tab:test8-final}
\begin{tabular}{c|c|rr|c|cc}
\hline
Initial condition & Window & $m_v^{\mathrm{DSMC}}$ & $m_v^{\mathrm{G}}$ & $\mathrm{Var}_{\mathrm{DSMC}}/\mathrm{Var}_{\mathrm{G}}$ & $L^1$ & $E_{\rm samp}$\\
\hline
Uncorrelated & 1 & $-0.3004$ & $-0.2999$ & 1.0107 & 0.0180 & 0.0163\\
             & 2 & $\phantom{-}0.0744$ & $\phantom{-}0.0746$ & 1.0037 & 0.0150 & 0.0147\\
\hline
Correlated   & 1 & $-0.4788$ & $-0.4784$ & 1.0045 & 0.0083 & 0.0079\\
             & 2 & $-0.1034$ & $-0.1032$ & 1.0157 & 0.0123 & 0.0078\\
\hline
\end{tabular}
\end{subtable}
\end{table}

\begin{figure}[p]
\centering
\begin{subfigure}[b]{\textwidth}
\centering
\includegraphics[width=\linewidth]{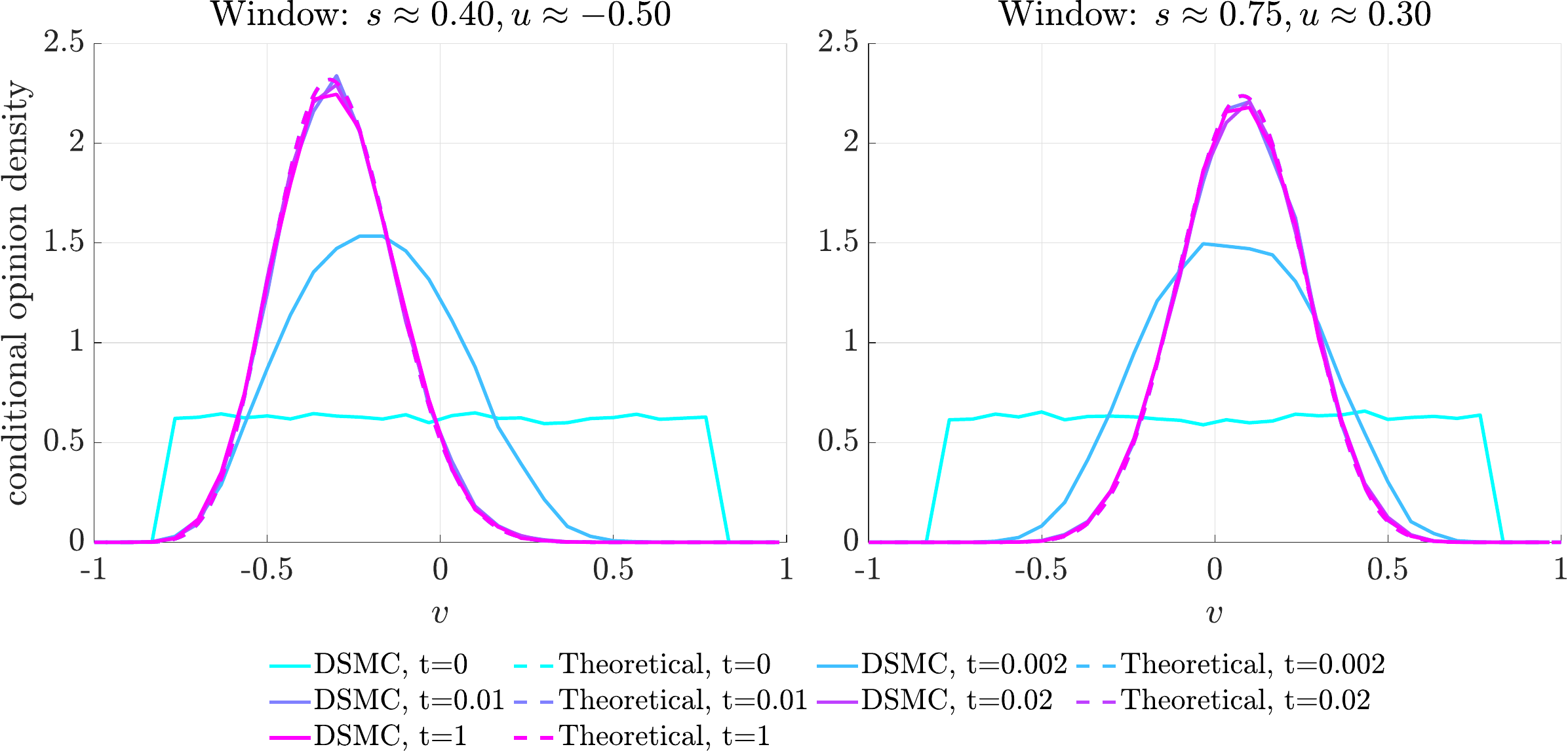}\\
\includegraphics[width=\linewidth]{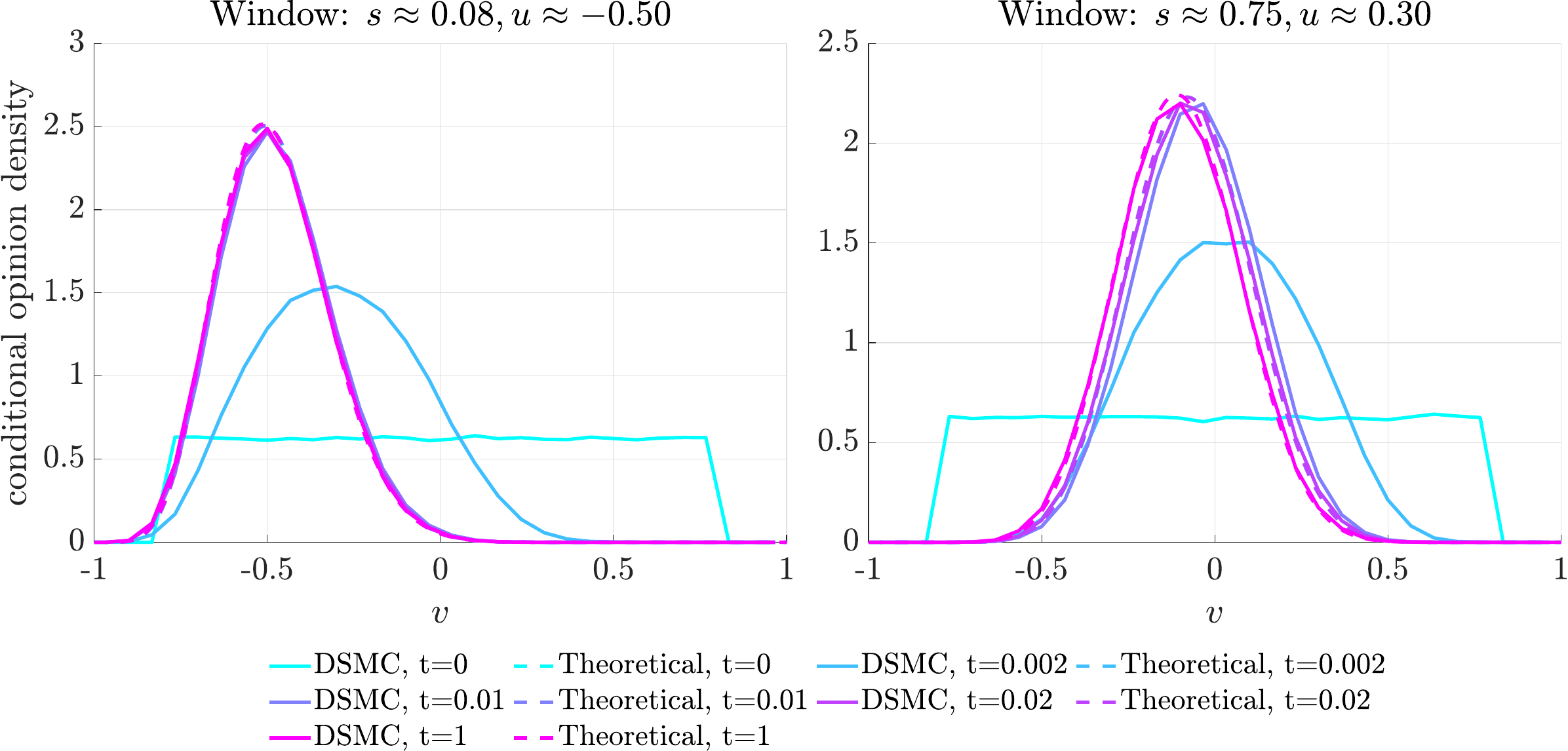}
\caption{Conditional opinion densities in the two selected $(s,u)$-windows for Case~1 (upper panel) and Case~2 (lower panel), at $t=0,\tau_v,5\tau_v,10\tau_v,\tau_\sigma$. Within each panel, the left and right plots correspond to the first and second windows, respectively. Solid lines: DSMC histograms. Dashed lines: mixtures of Gibbs densities with the self-consistent mean~\eqref{eq:test8-mveq}.}
\label{fig:test8a}
\end{subfigure}
\begin{subfigure}[b]{\textwidth}
    \centering
\includegraphics[width=0.49\linewidth]{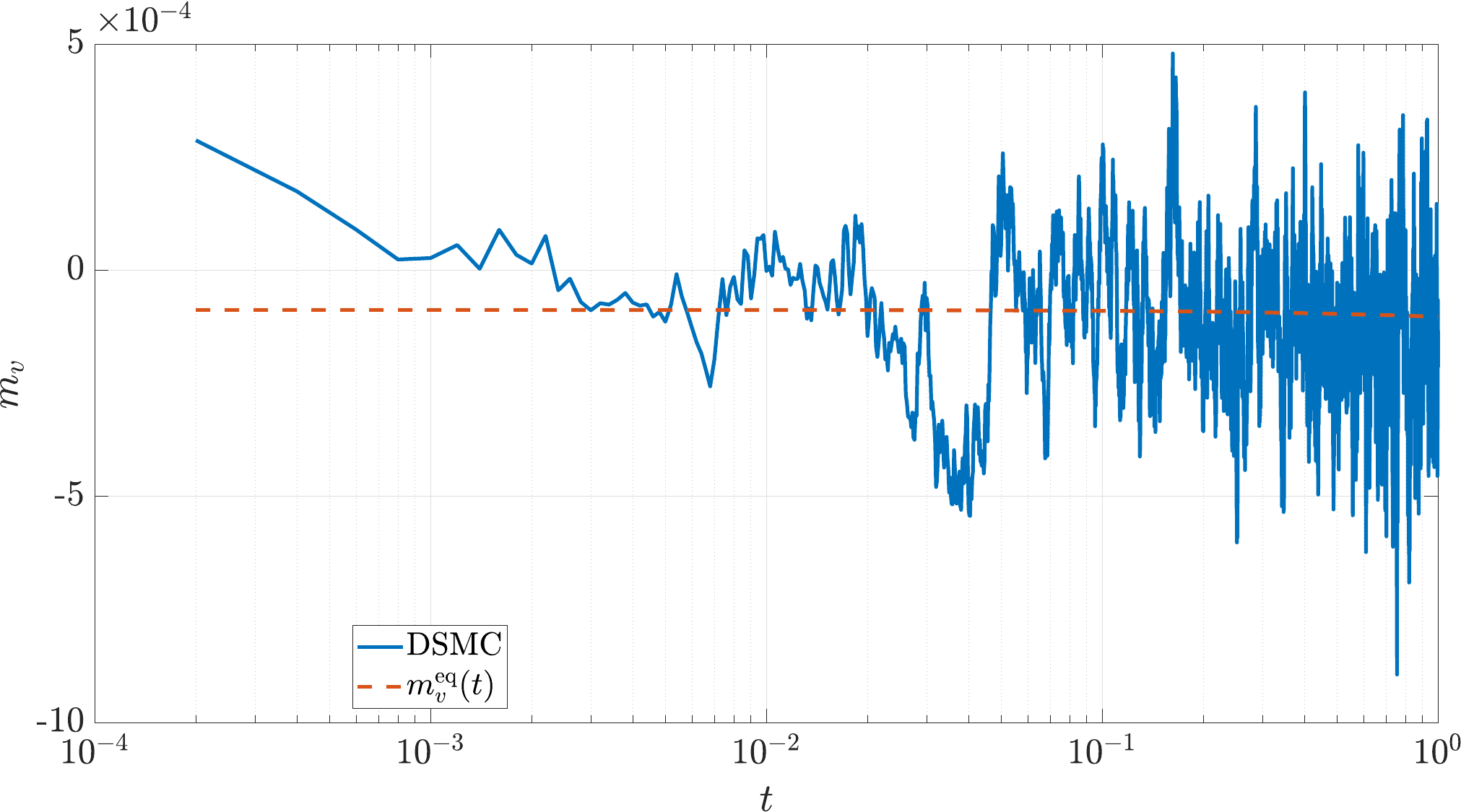} \includegraphics[width=0.49\linewidth]{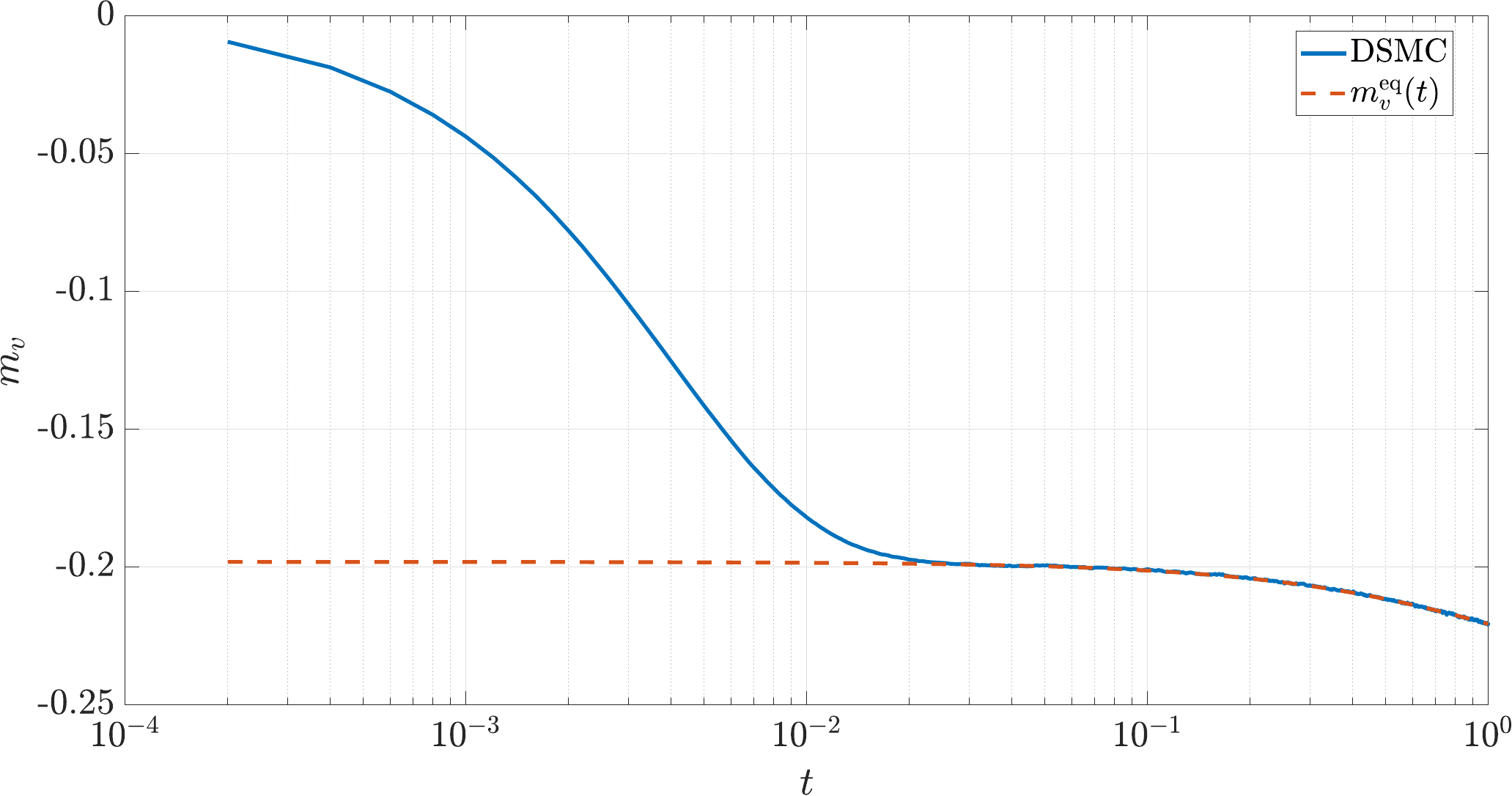}
\caption{Empirical mean opinion $m_v(t)$ (solid lines)  and quasi-equilibrium mean $m_v^{\rm eq}(t)$ (dashed lines), computed from the measured $m_s(t)$ and $\Sigma_{su}(t)$. The left and right panels correspond to Cases~1 and~2, respectively.}
\label{fig:test8b}
\end{subfigure}
\caption{Test~1. Fast-opinion relaxation and persistence of the conditional Gibbs quasi-equilibrium.}
\label{fig:test8}
\end{figure}

\paragraph{Numerical results.}

\emph{Relaxation of the conditional densities.}
Figure~\ref{fig:test8a} compares the empirical histograms with the Gibbs mixtures in the two windows, at the
observation times; the upper and lower panels refer to Cases~1 and~2, respectively. The $L^1$ distances are reported
in Table~\ref{tab:test8-transient}. At $t=0$ the conditional density is the uniform initial one and is far from the
Gibbs mixture (distances between $1.08$ and $1.18$). At $t=\tau_v$ the discrepancy is still large ($0.48$ and $0.34$ in
Case~1, $0.74$ and $0.45$ in Case~2). By $t=5\tau_v$ the distances have decreased to at most $0.023$, with the exception
of window~2 of Case~2, where it is $0.1106$. At $t=10\tau_v$ the distances are between $0.0068$ and $0.0148$, of the same
order as the sampling levels ($0.007$--$0.016$), again with the exception of window~2 of Case~2 ($0.0148$ against $0.0073$).

The residual shifts of the conditional means explain these late-time distances. Table~\ref{tab:test8-offset} compares the
measured offsets $M-\lambda_\infty$ (empirical minus predicted mean) with the laws~\eqref{eq:test8-laws} evaluated at
the window centers, with $m_s\simeq0.5$: measured and predicted offsets agree within about two sampling standard
errors at all three times, in both cases. The large residual of window~2 of Case~2 is due to the fact that the center of
this window weights $m_v$ with $s\simeq0.75$, while $m_v$ relaxes at the rate $a\simeq0.5/\tau_v$, slower than the rate
$b=1/\tau_v$ of the conditional mean. With $e_\lambda(0)=0.75\cdot0.198\simeq0.149$ and $e_M(0)\simeq0.074$,
\eqref{eq:test8-laws} gives $e_M(5\tau_v)\simeq0.023$, in agreement with the measured $0.0233$, and hence an $L^1$
distance of about $2\rho_{\max}e_M\simeq0.10$--$0.11$ (with $\rho_{\max}\simeq2.3$), in agreement with the observed
$0.1106$. At $10\tau_v$ the law gives $e_M\simeq0.002$, that is a contribution of about $0.009$, while the measured offset
is $0.0027$ (about $0.012$); together with the sampling level $0.0073$ this accounts for the observed $0.0148$. For the
other windows the residual shifts at $5\tau_v$ are $0.0057$ (Case~2, window~1), $0.002$ (Case~1, window~1) and
$5\cdot10^{-4}$ (Case~1, window~2), with $L^1$ contributions of about $0.026$, $0.009$ and $0.002$, which together with the
sampling level account for the observed $0.0224$, $0.0203$ and $0.0167$.

The variance of the conditional density is also still relaxing at $5\tau_v$. The second-moment equation
in~\eqref{eq:test8-moments} predicts variance ratios of $1.003$, $1.000$, $1.012$ and $1.004$ for the four windows (in the
order of Table~\ref{tab:test8-transient}), while the measured ones are $1.019$, $1.010$, $1.018$ and $1.008$, with sampling
standard errors of $0.8\%$ (Case~1) and $0.35\%$ (Case~2). The differences, $1.6\%$, $1.0\%$, $0.6\%$ and $0.4\%$, are
within about two standard errors and of the order of the scheme error discussed below.

\emph{Mean opinion.}
Figure~\ref{fig:test8b} displays the empirical mean $m_v(t)$ and the quasi-equilibrium prediction $m_v^{\rm eq}(t)$
computed from the measured $m_s(t)$ and $\Sigma_{su}(t)$; the left and right panels correspond to Cases~1 and~2. In the
uncorrelated case $\Sigma_{su}$ remains negligible, changing from approximately $4\cdot10^{-5}$ initially to
$5\cdot10^{-5}$ at $t=1$, which is compatible with the sampling fluctuation $\sqrt{\Sigma_s\Sigma_u/N}\simeq7\cdot10^{-5}$.
Consequently the quasi-equilibrium mean remains close to zero ($m_v^{\rm eq}\simeq-1\cdot10^{-4}$), the measured $m_v$
fluctuates between $-2\cdot10^{-4}$ and $3\cdot10^{-4}$, and the final DSMC mean $m_v(1)=-0.0002$ agrees with
$m_v^{\rm eq}(1)\simeq-0.0001$ within the sampling error. In the correlated case the covariance remains significant,
changing from $0.09907$ to $0.09797$, while the mean strategy increases from $0.5000$ to $0.5560$. As a consequence the
quasi-equilibrium mean varies from $m_v^{\rm eq}(0)\simeq-0.198$ to $m_v^{\rm eq}(1)\simeq-0.2207$, and the latter value
is accurately reproduced by the DSMC mean, $m_v(1)=-0.2207$.

The relaxation of $m_v$ towards $m_v^{\rm eq}$ in Case~2 follows the first law in~\eqref{eq:test8-laws}: with
$m_s\simeq0.5$ it predicts $m_v(\tau_v)\simeq-0.078$, $m_v(5\tau_v)\simeq-0.182$ and $m_v(10\tau_v)\simeq-0.197$, while the
measured values are $-0.0780$, $-0.1820$ and $-0.1973$. The mean opinion thus converges to $m_v^{\rm eq}$ on the scale
$\tau_v/(1-m_s)$, and not to $m_u=0$, which confirms that the initial correlation between strategy and conviction has a
persistent effect on the macroscopic opinion dynamics, even after the fast opinion relaxation has taken place.

\emph{Persistence of the quasi-equilibrium on the slow scale.}
Table~\ref{tab:test8-final} compares the empirical conditional distributions with the Gibbs mixtures at $t=1$. In all four
windows the empirical and predicted means differ by at most $5\cdot10^{-4}$, below the sampling error of the mean
($\simeq10^{-3}$ in Case~1 and $4\cdot10^{-4}$ in Case~2). The $L^1$ distances are at the level of the sampling error: their
ratio to $E_{\rm samp}$ is $1.10$, $1.02$, $1.05$ and $1.58$ for the four windows. The variance ratios lie between $1.004$ and
$1.016$ at $t=1$ and between $0.998$ and $1.018$ at $10\tau_v$; averaged over these two times they are $1.004$ in window~1
and $1.012$ in window~2, with a standard error of about $0.3\%$. This small excess, larger in the window where agents give
more weight to the interaction partner ($s\simeq0.75$), is of the order of the $\mathcal O(\alpha_v)$ error of the particle
scheme, with $\alpha_v=0.005$. For a particle with center $\lambda$ and weight $s$ on the partner, the explicit relaxation
step and the sampling of the partner in the opinion update give the stationary variance
$\bigl[\nu^2(1-\lambda^2)+\alpha_vs^2\mathrm{Var}(v)\bigr]/(2+\nu^2-\alpha_v)$, where $\mathrm{Var}(v)$ is the variance of the
population, instead of $\nu^2(1-\lambda^2)/(2+\nu^2)$. With $\mathrm{Var}(v)\simeq0.07$ the relative excess is about $0.3\%$ for
$s\simeq0.4$ and $0.6\%$ for $s\simeq0.75$, compatible with the measured values within about two standard errors. These results confirm that the conditional opinion distributions remain close to their Gibbs quasi-equilibria also on the
slow strategy time scale, while the location of these equilibria depends strongly on the strategy--conviction correlation.

\subsection{Test~2. Fast-strategy closure: first- and second-order moments}

This test examines the fast-strategy regime $\tau_\sigma\ll\tau_v$ of Section~\ref{sec:asymptotic} for $K=2$ strategies and $D(v)=\sqrt{1-v^2}$. The particle simulation approximates the full distribution $f$ at the finite ratio $\epsilon_\sigma=2\cdot 10^{-3}$; its moments (written without superscript) are compared with those of the reduced model under the closure~(A3) (written with superscript $(0)$). For the opinion mean $m_v$ and the opinion--conviction covariance $\Sigma_{vu}$ we expect agreement up to $\mathcal O(\epsilon_\sigma)$, whereas the stationary opinion variance $\Sigma_v$ exhibits a persistent $\mathcal O(1)$ closure error whenever both the limiting strategy distribution and the conviction distribution are non-degenerate (Remark~\ref{rem:A3-defect}). The test is designed to illustrate: the relaxation of $m_v$ at rate $(1-\overline\Theta^{\rm eq})/\tau_v$ and the persistence of memory for $\overline\Theta^{\rm eq}=1$; the accuracy of the reduction for $m_v$ and $\Sigma_{vu}$; and the variance defect of~(A3) and its origin in a mixed moment.

\paragraph{Setup.} We consider the following configuration. 

\emph{Initial conditions.} The strategy is uniformly distributed, $s\sim\mathrm{Unif}(\Delta_1)$, and independent of the opinion--conviction pair $(v,u)$, so that $f_0(v,s,u)=h_0(v,u)g_0(s)$ as in assumption~(A1). Opinions and convictions are generated from a common latent variable,
\[
    v=-0.5+0.4z,\qquad u=0.5+0.32z+0.08w,
\]
with $z,w\sim\mathrm{Unif}[-1,1]$ independent. The joint initial density is
\[
    f_0(v,s,u)
    =
    \frac{125}{16}\,
    \mathbf 1_{[-0.9,-0.1]}(v)\,
    \mathbf 1_{[0,1]}(s)\,
    \mathbf 1_{[0.82+0.8v,\;0.98+0.8v]}(u),
\]
and the initial moments are
\[
    m_v(0)=-0.5,\quad m_u=0.5,\quad
    \Sigma_v(0)\simeq0.0534,\quad
    \Sigma_u\simeq0.0363,\quad
    \Sigma_{vu}(0)\simeq0.0427,\quad
    \Sigma_{su}(0)=0.
\]
The large offset $m_v(0)-m_u$ makes the evolution of $m_v$ non-trivial, and the correlated initial pair $(v,u)$ allows us to follow the relaxation of $\Sigma_{vu}$. Since the strategy update acts identically on all agents, the independence of $s$ and $u$ is preserved, so that $\Sigma_{su}(t)=0$.

We consider four constant payoff matrices, corresponding to $\overline\Theta^{\rm eq}\in\{0,0.3,0.6,1\}$ (stubbornness dominance, coexistence with $\hat m_s=0.3$ and $\hat m_s=0.6$, and compromise
dominance),
\[
J^{(\overline\Theta^{\rm eq}=0)}
=
\begin{pmatrix}-2&-1\\0&0\end{pmatrix}, \quad
J^{(\overline\Theta^{\rm eq}=0.3)}
=
\begin{pmatrix}-7&3\\0&0\end{pmatrix}, \quad
J^{(\overline\Theta^{\rm eq}=0.6)}
=
\begin{pmatrix}0&3\\2&0\end{pmatrix}, \quad
J^{(\overline\Theta^{\rm eq}=1)}
=
\begin{pmatrix}2&1\\0&-1\end{pmatrix},
\] thus satisfying assumption (A2).

\emph{Parameters.} The parameters are $\tau_v=1$, $\tau_\sigma=0.002$, $N=2\cdot10^6$, $\Delta t=10^{-4}$ and $T=15$. With the bounded uniform perturbations of the particle scheme, the effective diffusion intensity is $\nu_{\rm eff}\simeq0.0722$ for the nominal value $\nu=0.25$.

\emph{Reduced model.} Under~(A3) the moments of the reduced model satisfy, in the original time variable $t=\tau\cdot\tau_v$, the closed system~\eqref{eq:mean-opinion-fast-strategy}--\eqref{eq:fast-strategy-moments-general}, where the first two equations are solved explicitly by \eqref{sol_mv}-\eqref{sol_svu}. The strategy layer has duration $\mathcal O(\tau_\sigma)$ and does not modify the moments of $(v,u)$ at leading order, so the system is initialized with the moments of the particle system:
$m_v^{(0)}(0)=m_v(0)$, $\Sigma_{vu}^{(0)}(0)=\Sigma_{vu}(0)$, $\Sigma_v^{(0)}(0)=\Sigma_v(0)$.

\emph{Strategy update.} In the particle scheme the gain is evaluated at the population mean, $G(m_s^n)$, for all particles, rather than at the strategy of a random partner; the reasons are discussed in Remark~\ref{rem:strategy-noise}. The opinion update retains the random pairing.

\paragraph{Predictions.}
We first compute the asymptotic moments predicted by the full kinetic model and by the reduced fast-strategy model under closure (A3), in order to provide reference values for the numerical simulations.

\emph{Full model.} For $\overline\Theta^{\rm eq}<1$ we have $\Sigma_{su}\equiv0$ and $\mathcal I(t)\to+\infty$, so that $m_v^\infty=m_u$, while $\Sigma_{vu}^\infty$ and $\Sigma_v^\infty$ are given by~\eqref{eq:Sv-exact-coexistence}, which for $\overline\Theta^{\rm eq}=0$ holds with $\Sigma_s^{\rm eq}=0$. For $\overline\Theta^{\rm eq}=1$ the strategy distribution converges to $\delta(s-1)$ but the opinion mean retains the effect of the initial strategy layer: $m_v^\infty$ is given by~\eqref{eq:mv-infty-comp} with $\Sigma_{su}\equiv0$ and
\[
    \mathcal I_\infty=\frac{\tau_\sigma}{G(1)}\int_0^1 g_0(s)\ln\frac1s\,\mathrm{d}s=\frac{\tau_\sigma}{2}=0.001,
\]
where we used $G\equiv2$ ($G(0)=G(1)=2$), which makes the integral formula valid, and $g_0\equiv1$ on $\Delta_1$. Moreover,
\[
    \Sigma_{vu}^\infty=0,
    \qquad
    \Sigma_{v}^\infty
    =
    \frac{
        \nu_{\rm eff}^2
        \left[1-(m_{v}^\infty)^2\right]
    }{2+\nu_{\rm eff}^2}.
\]

\emph{Reduced model.} The stationary moments $m_v^{(0),\infty}$, $\Sigma_{vu}^{(0),\infty}$, $\Sigma_v^{(0),\infty}$ are those of Proposition~\ref{prop:fast-strategy-stationary}, with $\nu\to\nu_{\rm eff}$; in particular $m_v^{(0),\infty}=m_v(0)$ for $\overline\Theta^{\rm eq}=1$. The variance $\Sigma_s^{\rm eq}$ is computed by quadrature from the projective push-forward $g^{\rm eq}=(\Phi_{\mathcal T_\infty})_\#g_0$. The predictions are collected in Table~\ref{tab:test6}.

\begin{table}[ht]
\centering
\small
\renewcommand{\arraystretch}{1.2}
\begin{tabular}{c ccc ccc c}
\toprule
& \multicolumn{3}{c}{Full model}
& \multicolumn{3}{c}{Reduced model (A3)}
& \\
\cmidrule(lr){2-4}
\cmidrule(lr){5-7}
$\overline\Theta^{\rm eq}$
& $m_v^\infty$
& $\Sigma_{vu}^\infty$
& $\Sigma_v^\infty$
& $m_v^{(0),\infty}$
& $\Sigma_{vu}^{(0),\infty}$
& $\Sigma_v^{(0),\infty}$
& $\Sigma_s^{\rm eq}$
\\
\midrule
$0$
& $0.500058$
& $0.036295$
& $0.038148$
& $0.500058$
& $0.036295$
& $0.038148$
& $0$
\\
$0.3$
& $0.500058$
& $0.025406$
& $0.022133$
& $0.500058$
& $0.025406$
& $0.019686$
& $0.067593$
\\
$0.6$
& $0.500058$
& $0.014518$
& $0.010612$
& $0.500058$
& $0.014518$
& $0.007740$
& $0.079349$
\\
$1$
& $-0.498913$
& $0$
& $1.950\cdot10^{-3}$
& $-0.499913$
& $0$
& $1.948\cdot10^{-3}$
& $0$
\\
\bottomrule
\end{tabular}
\caption{Test~2. Asymptotic moment predictions for the full kinetic model and the reduced fast-strategy model under closure~(A3). The predictions use the empirical values $m_v(0)=-0.499913$, $m_u=0.500058$, $\Sigma_u=0.036295$ and $\nu_{\rm eff}=0.072162$. For $\overline\Theta^{\rm eq}=1$, the full-model values account for the finite strategy-relaxation time
$\tau_\sigma=0.002$.}
\label{tab:test6}
\end{table}

\begin{figure}[p]
\centering
\begin{subfigure}[ht]{\textwidth}
\centering
\includegraphics[width=\linewidth]{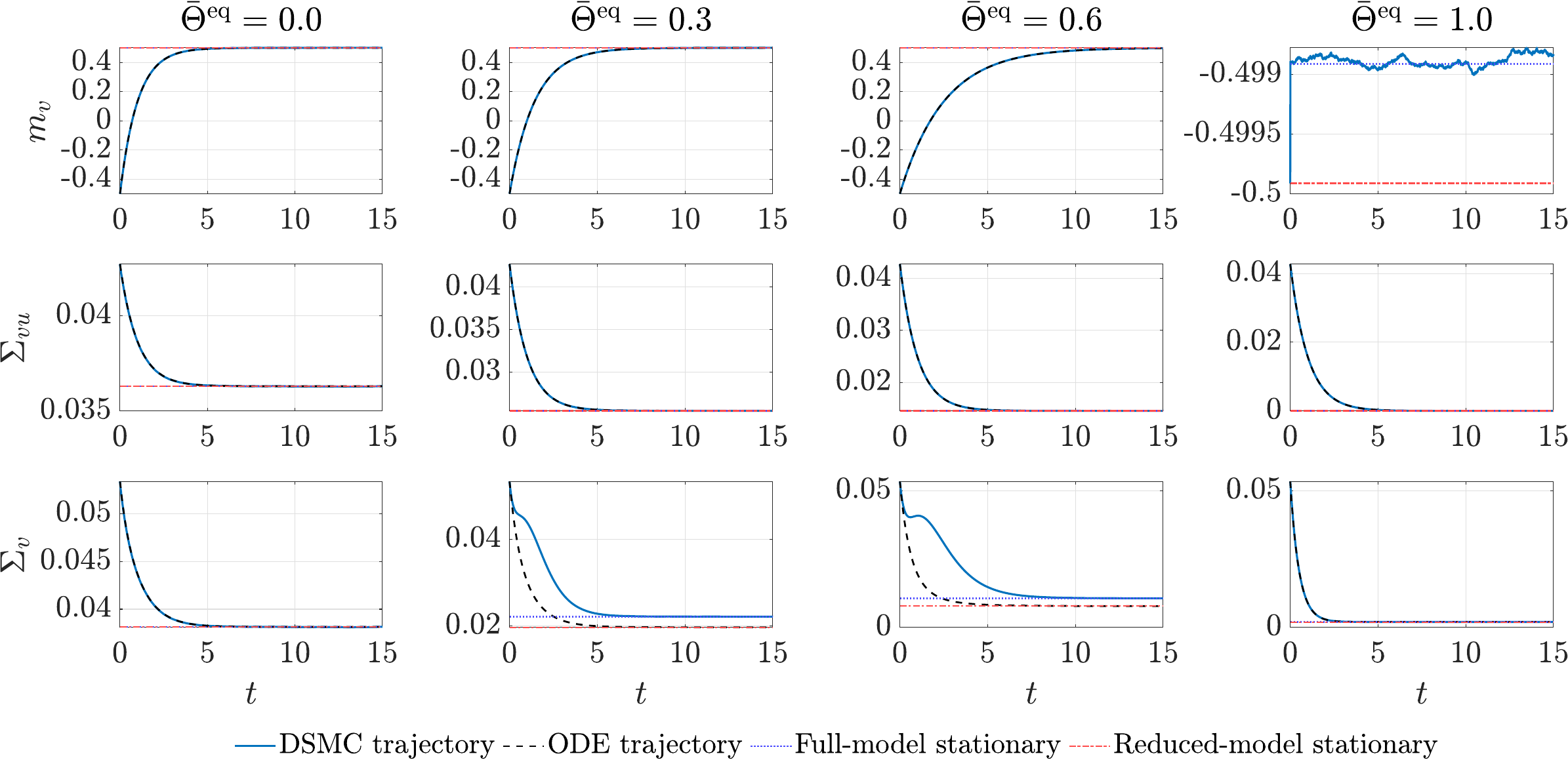}
\caption{Evolution of the opinion mean $m_v$ (top), opinion--conviction covariance $\Sigma_{vu}$ (middle) and opinion variance $\Sigma_v$ (bottom) for $\overline\Theta^{\rm eq}\in\{0,0.3,0.6,1\}$ (left to right). DSMC trajectories are compared with the reduced solution and with the stationary values of the full and reduced models.}
\label{fig:test6a}
\end{subfigure}
\begin{subfigure}[ht]{\textwidth}
\centering
\includegraphics[width=0.9\linewidth]{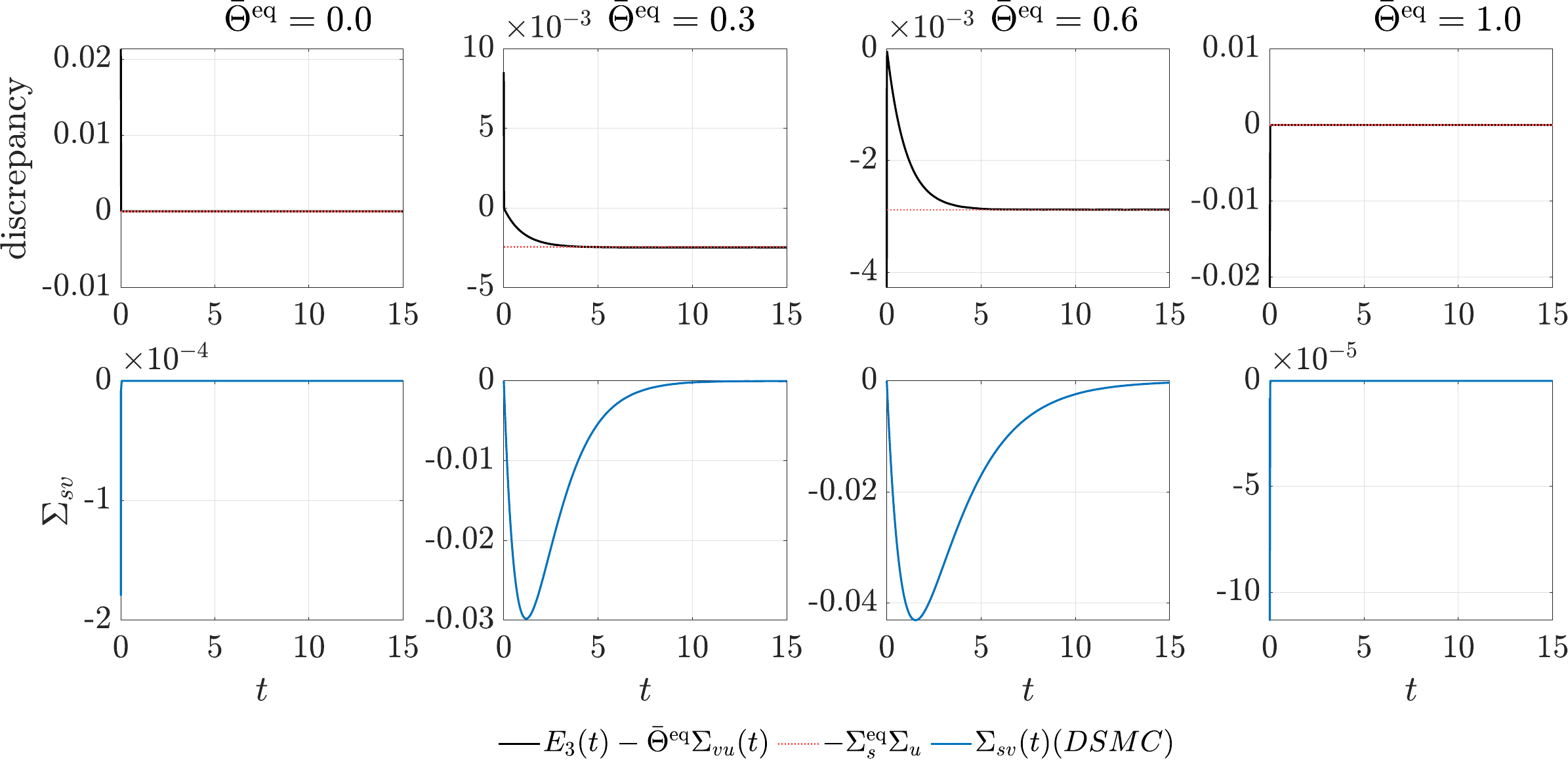}
\caption{Diagnostics of the closure~(A3). Top: mixed-moment defect $E_3-\overline\Theta^{\rm eq}\Sigma_{vu}$, with its predicted stationary value $-\Sigma_s^{\rm eq}\Sigma_u$ (dotted). Bottom: strategy--opinion covariance $\Sigma_{sv}$, which vanishes asymptotically in all four cases.}
    \label{fig:test6b}
\end{subfigure}
\caption{Test~2. Fast-strategy reduction: comparison of the full and reduced moment dynamics and assessment of the factorization closure~(A3).}
\label{fig:test6}
\end{figure}
\paragraph{Numerical results.}

Figure~\ref{fig:test6a} compares the DSMC evolution of $m_v$, $\Sigma_{vu}$ and $\Sigma_v$ with the predictions of the reduced moment system and the corresponding stationary values of the full and reduced
models.

\emph{Mean opinion and opinion--conviction covariance.}
In all four cases the DSMC trajectories of $m_v$ and $\Sigma_{vu}$ follow the reduced solution after a layer of duration $\mathcal O(\tau_\sigma)$, in which the strategy distribution relaxes towards $g^{\rm eq}$. For $\overline\Theta^{\rm eq}<1$ both models converge to the same limits of $m_v$ and $\Sigma_{vu}$, and the layer discrepancy is removed by the relaxation of the opinion dynamics. The mean opinion relaxes at rate $(1-\overline\Theta^{\rm eq})/\tau_v$, namely $1$, $0.7$ and $0.4$ for $\overline\Theta^{\rm eq}=0,\,0.3,\,0.6$, so that for $\overline\Theta^{\rm eq}=0.6$ a small residual of approximately
$|m_v(0)-m_u|\mathrm{e}^{-6}\simeq2.5\cdot10^{-3}$ is still present at $T=15$; $\Sigma_{vu}$ relaxes at rate $1/\tau_v$ in all cases. For $\overline\Theta^{\rm eq}=1$ the behavior is different: the reduced model conserves $m_v^{(0)}\equiv m_v(0)$, while the DSMC mean converges to  the full-model value accounting for the finite strategy layer. This displacement, of order $(m_u-m_v(0))\,\mathcal I_\infty/\tau_v=\mathcal O(\tau_\sigma)$, persists at stationarity, because $1-m_s(t)$ is integrable and the mean is frozen once the strategies have reached $s=1$; it vanishes as $\tau_\sigma/\tau_v\to0$. The covariance $\Sigma_{vu}$ still converges to zero, consistently with both models.

\emph{Opinion variance and closure defect.}
For $\overline\Theta^{\rm eq}=0$ the full and reduced predictions coincide, since the limiting strategy distribution is concentrated at $s=0$. In the two coexistence cases the DSMC variance converges to the full-model value~\eqref{eq:Sv-exact-coexistence}, which is strictly larger than the value obtained under~(A3); the difference agrees with the correction~\eqref{eq:Sv-defect}, namely $0.00245$ and $0.00287$ for $\overline\Theta^{\rm eq}=0.3$ and $0.6$ (relative gaps of about $12.4\%$ and $37.1\%$). In the two coexistence cases, the variance discrepancy persists as $\tau_\sigma/\tau_v\to0$, since $\Sigma_s^{\rm eq}$ is independent of $\tau_\sigma$. For $\overline\Theta^{\rm eq}=1$, instead, $\Sigma_s^{\rm eq}=0$ and the difference between the full and reduced predictions is solely due to the finite strategy-relaxation transient. This difference vanishes as $\tau_\sigma/\tau_v\to0$. 

Figure~\ref{fig:test6b} shows the origin of the discrepancy. In the upper panels, the quantity $E_3(t)-\overline\Theta^{\rm eq}\Sigma_{vu}(t)$ converges to the predicted value $-\Sigma_s^{\rm eq}\Sigma_u$ in the coexistence cases and to zero in the dominance cases, while, in ther lower panels, $\Sigma_{sv}$ decays to zero in all four cases. The persistent variance discrepancy is therefore entirely associated with the mixed moment $E_3$, and not with a non-vanishing stationary covariance
$\Sigma_{sv}$.

\emph{Control of the numerical scheme.} In the two coexistence cases the empirical strategy variances at the final time are $\Sigma_s(T)=0.06846$ and $0.07939$, to be compared with the continuous-flow predictions $\Sigma_s^{\rm eq}$ (relative differences $1.28\%$ and $0.05\%$). The small differences are consistent with the finite-step Euler discretization of the replicator flow, and are larger for the stiffer payoff ($|G(1)|=7$). The empirical $\Sigma_{su}$ remains negligible, consistently with the initial independence of strategies and convictions ($\max_{t}|\Sigma_{su}(t)|<4.4\cdot10^{-6}$ in the
coexistence cases). The residual oscillations of the DSMC trajectories are due to finite-particle fluctuations in the opinion dynamics; they are most visible for $\overline\Theta^{\rm eq}=1$, where the mean opinion has no restoring mechanism after the strategy layer.

Overall, the fast-strategy reduction accurately captures the opinion mean and the opinion--conviction covariance, up to $\mathcal O(\tau_\sigma)$ for $\overline\Theta^{\rm eq}=1$, while the closure~(A3) leaves a persistent error on the opinion variance when the limiting strategy distribution is non-degenerate.

\subsection{Test~3. Two-strategy dynamics across the four asymptotic regimes}

This test illustrates the evolution of the strategy marginal $g(s,t)$ across the four asymptotic regimes identified in Section~\ref{sec:2s-constant} for the two-strategy model with constant payoffs. We consider Regime~1 (compromise dominance), Regime~2 (stubbornness dominance), Regime~3 (stable coexistence), and the threshold state of Regime~4 (bistability). The off-threshold cases of Regime~4 are omitted, since their long-time behavior is qualitatively analogous to that of Regime~1 and Regime~2 depending on the basin of attraction.

\paragraph{Setup.} The four payoff matrices are 
\begin{equation}\label{J_test}
J^{(1)}=
\begin{pmatrix}
2 & 1\\
0 & -1
\end{pmatrix},
\qquad
J^{(2)}=
\begin{pmatrix}
-1 & -2\\
1 & 0
\end{pmatrix},
\qquad
J^{(3)}=
\begin{pmatrix}
0 & 3\\
2 & 0
\end{pmatrix},
\qquad
J^{(4)}=
\begin{pmatrix}
2 & -2\\
0 & 0
\end{pmatrix}.
\end{equation} For each matrix, the values of $G(0)=J_{CS}-J_{SS}$ and $G(1)=J_{CC}-J_{SC}$ confirm the intended regime: $J^{(1)}$ gives $G(0)=G(1)=2>0$ (Regime~1); $J^{(2)}$ gives $G(0)=G(1)=-2<0$ (Regime~2); $J^{(3)}$ gives $G(0)=3>0$ and $G(1)=-2<0$, with stable zero $\hat m_s=G(0)/(G(0)-G(1))=3/5=0.6$ (Regime~3); $J^{(4)}$ gives $G(0)=-2<0$ and $G(1)=2>0$, with unstable zero $\hat m_s=0.5$ (Regime~4).

The initial distribution is the symmetric $\mathrm{Beta}(2,2)$ density $g_0(s)=6s(1-s)$, centered at $m_s(0)=0.5$. The strategy marginal is evolved via the exact push-forward representation~\eqref{eq:g-pushforward-two}: at each time $t$, $g(s,t)=(\Phi_{\mathcal{T}(t)})_\#g_0$ is obtained by transporting $g_0$ through the projective map $\Phi_{\mathcal{T}(t)}$ of~\eqref{eq:two-strategy-flow}, where $\mathcal{T}(t)$ is the solution of~\eqref{eq:cumulative-gain-equation-two}.
\begin{figure}[ht]
\centering
\includegraphics[width=1\linewidth]{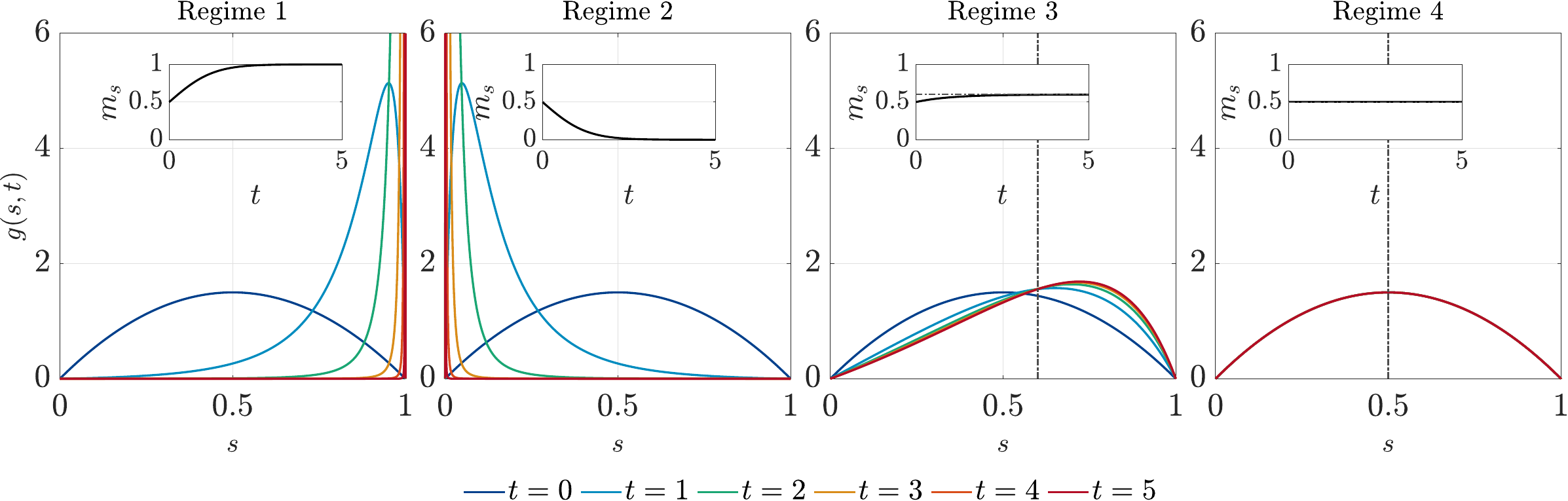}
\caption{Test 3. Evolution of the strategy marginal $g(s,t)$ across the four asymptotic  regimes. Distributions are shown at $t=0,1,2,3,4,5$, the inset in each panel shows the corresponding evolution of $m_s(t)$. In the Regime~3 and Regime~4 panels, the vertical black dash-dotted line indicates $\hat m_s$.}
\label{fig:strategy-density-regimes}
\end{figure}

\paragraph{Numerical results.}
Figure~\ref{fig:strategy-density-regimes} shows the strategy marginal at several time snapshots, together with the evolution of $m_s(t)$ in the inset of each panel.

In Regime~1, $G(m_s)>0$ throughout, so $\mathcal{T}(t)\to+\infty$ and the population progressively concentrates at the compromise strategy: $g(s,t)\to\delta(s-1)$ and $m_s(t)\to1$.

In Regime~2, $G(m_s)<0$ throughout, so $\mathcal{T}(t)\to-\infty$ and the population is driven to stubbornness: $g(s,t)\to\delta(s)$
and $m_s(t)\to0$.

In Regime~3, $G'=-5<0$ and $G$ has a stable zero at $\hat m_s=0.6$. The mean strategy converges to this interior equilibrium, while the strategy distribution does not collapse to a Dirac mass: the asymptotic profile $g^\infty=(\Phi_{\mathcal{T}_\infty})_\#g_0$ retains a non-trivial shape inherited from $g_0$, deformed by the nonlinear projective flow.

In Regime~4, the initial mean $m_s(0)=\hat m_s=0.5$ coincides exactly with the unstable equilibrium. Consequently $G(m_s(0))=0$, the cumulative gain satisfies $\mathcal{T}(t)\equiv0$ by uniqueness of equation~\eqref{eq:cumulative-gain-equation-two}, and the projective map reduces to the identity: $g(s,t)=g_0(s)$ and $m_s(t)\equiv0.5$ for all $t$. The numerical simulation confirms
this threshold behavior.

\subsection{Test~4. Asymptotic behavior and memory effects in the mean opinion}

This test investigates the evolution of the mean opinion $m_v(t)$ in Regime~1 and Regime~2, using the payoff matrices $J^{(1)}$ and $J^{(2)}$ from~\eqref{J_test}. The test is designed to isolate two qualitatively distinct phenomena: the persistence or loss of memory of the initial mean opinion, and the role of transient strategy--conviction correlations.

\paragraph{Setup.}
The conviction distribution is fixed as $k(u) = \tfrac 12\,\mathbf{1}_I(u)$ giving $m_u=0$, and remains constant in time. For the initial public opinion, we consider three uniform distributions
\begin{equation*}
    h_0(v)
    =
    \mathbf{1}_A(v)/|A|,
    \quad
    A\in\bigl\{[-1,0],\,I,\,[0,1]\bigr\},
\end{equation*}
giving $m_v(0)\in\{-0.5,0,0.5\}$. We consider two initial joint configurations of strategy and conviction, differing only in their initial correlation.

\textbf{Case~1 ($s$ and $u$ initially independent):}
\begin{equation*}
    f_0(v,s,u) = h_0(v)\,\mathbf{1}_{\Delta_1}(s)\,k(u),
    \qquad
    m_s(0)=0.5,
    \qquad
    \Sigma_{su}(0)=0.
\end{equation*}
Since the payoffs are constant, the strategy dynamics is independent of $u$. Hence the initial independence of $s$ and $u$ is preserved, and $\Sigma_{su}(t)=0$ for all $t\geq0$.

\textbf{Case~2 ($s$ and $u$ initially correlated):}
\begin{equation*}
    f_0(v,s,u)
    =
    h_0(v)\,
    \delta\!\left(s-\tfrac 12(u+1)\right)k(u),
\end{equation*}
Here $s=\frac{1}{2}(u+1)\in[0,1]$ for all $u\in[-1,1]$, so the strategy is a deterministic increasing function of conviction. Since $u\sim\mathrm{Unif}(-1,1)$, this gives $s\sim\mathrm{Unif}(0,1)$ and $m_s(0)=0.5$, while $\Sigma_{su}(0)=\mathrm{Cov}(\tfrac{1}{2}(u+1),u) =\tfrac{1}{2}\mathrm{Var}(u)=\tfrac{1}{6}$. The two cases thus share the same $m_s(0)=0.5$ and the same marginal $g_0(s)=\mathbf{1}_{[0,1]}(s)$, but differ in their initial strategy--conviction correlation.

The simulations are performed  up to the final time $T=5$.

\begin{table}[t]
\centering
\small
\renewcommand{\arraystretch}{1.1}
\setlength{\tabcolsep}{5pt}
\begin{tabular}{ccc|ccc}
\toprule
\multicolumn{3}{c|}{\textbf{Case~1} ($\Sigma_{su}(0)=0$)}
&\multicolumn{3}{c}{\textbf{Case~2} ($\Sigma_{su}(0)=1/6$)}\\
\cmidrule(r){1-3}\cmidrule(l){4-6}
$m_v(0)$
&$m_v^\infty$, Regime~1
&$m_v^\infty$, Regime~2
&$m_v(0)$
&$m_v^\infty$, Regime~1
&$m_v^\infty$, Regime~2\\
\midrule
$-0.5$ &$\approx-0.3$ &$0$ &$-0.5$ &$\approx-0.5$ &$0$\\
$0$    &$0$           &$0$ &$0$ &$\approx-0.2$ &$0$\\
$0.5$  &$\approx0.3$  &$0$ &$0.5$  &$\approx0.1$  &$0$\\
\bottomrule
\end{tabular}
\caption{Test~4. Analytical predictions of $m_v^\infty$ for the twelve simulations in Test~4. In Regime~2, all initial conditions give $m_v^\infty=m_u=0$ regardless of $\Sigma_{su}(0)$.}
\label{tab:test2}
\end{table}

\paragraph{Analytical predictions.}
In both Regime~1 and Regime~2, the strategy distribution converges to a pure state, so $\Sigma_{su}^\infty=0$ regardless of initial correlation. The covariance $\Sigma_{su}(t)$ is therefore purely transient, but it contributes to the integral term in~\eqref{eq:mv-exact}. For Case~1, $\Sigma_{su}(t)\equiv0$ and formula~\eqref{eq:mv-exact} gives directly:
\begin{itemize}
\item Regime~1: $m_v^\infty=m_v^{\rm comp} =(m_v(0)-m_u)\mathrm{e}^{-\mathcal{I}_\infty/\tau_v}$ (since $m_u=0$); in general $m_v^\infty\neq0$.
\item Regime~2: $m_v^\infty=m_u=0$ for all three initial conditions.
\end{itemize} For Case~2, the transient covariance $\Sigma_{su}(t)>0$ generates an additional negative shift in $m_v^\infty$ via the integral term in~\eqref{eq:mv-exact}: the positive initial correlation between $s$ and $u$ causes the population to over-represent high-$s$ agents with positive $u$, driving $m_v$ toward a value below the Case~1 prediction. The resulting predictions are summarized in Table~\ref{tab:test2}.

\paragraph{Numerical results.} \begin{figure}
\centering
\includegraphics[width=1\linewidth]{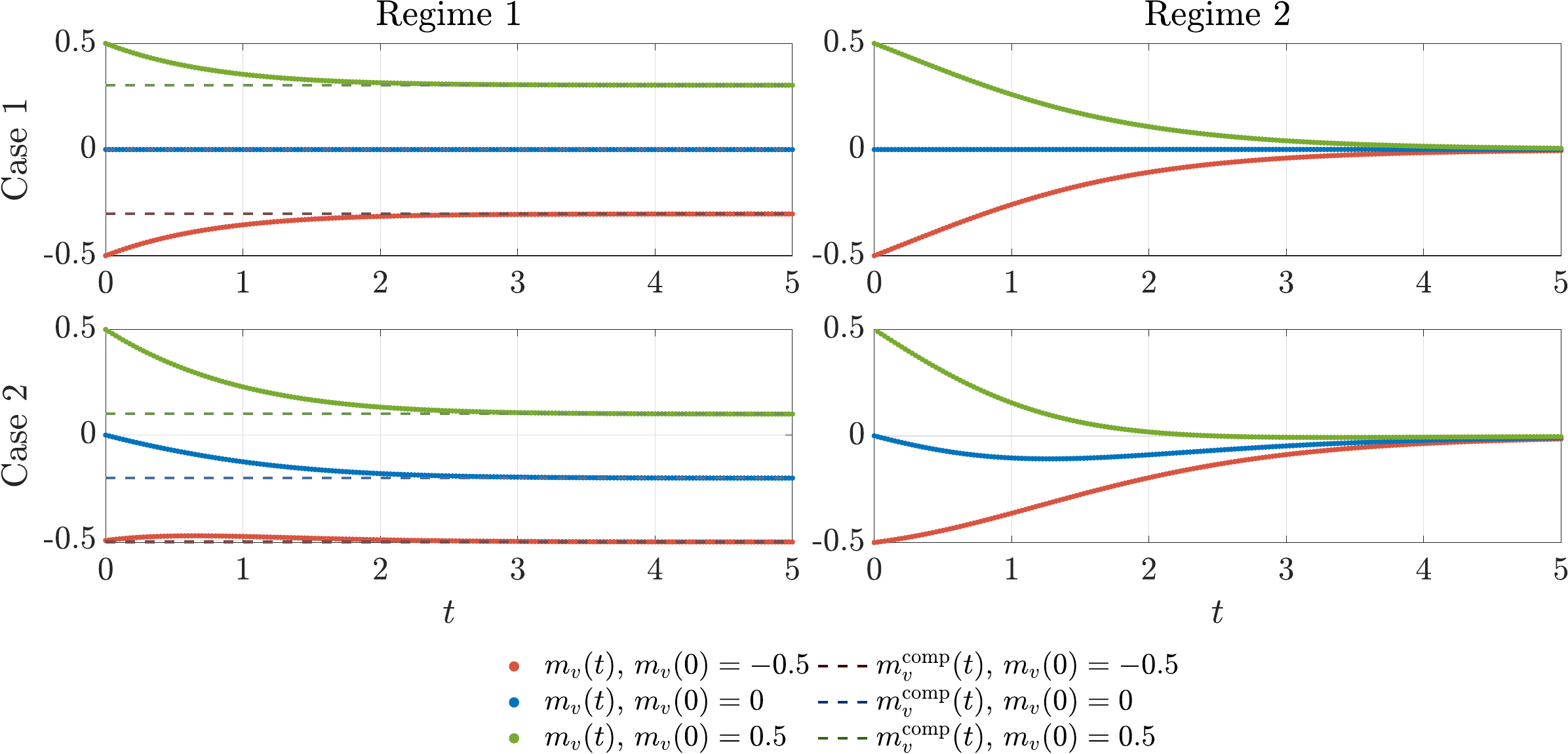}
\caption{Test~4. Evolution of $m_v(t)$ for three initial opinions in Regime~1 (left column) and Regime~2 (right column). Rows correspond to Case~1 (top, $\Sigma_{su}(0)=0$) and Case~2 (bottom, $\Sigma_{su}(0)=\tfrac{1}{6}$). Colors: blue, green, red for $m_v(0)=-0.5,0,0.5$. In Regime~1, horizontal dashed lines indicate the analytical asymptotic values $m_v^{\rm comp}$ from~\eqref{eq:mv-infty-comp}.}
\label{fig:test2}
\end{figure}

Figure~\ref{fig:test2} summarizes the numerical results. In Regime~1, the trajectories converge to the analytical values $m_v^{\rm comp}$ for both initial configurations. The different limiting values across the three initial opinions reflect the persistence of the initial memory. In Case~2, the additional negative shift due to the transient covariance $\Sigma_{su}(t)$ is visible: all three asymptotic values are lower than in Case~1, consistently with the predictions of Table~\ref{tab:test2}.

In Regime~2, all trajectories converge to $m_v^\infty=m_u=0$, regardless of the initial configuration. The transient dynamics is nonetheless affected by $\Sigma_{su}(t)$: in Case~2, the positive covariance acts as a temporary negative forcing in the mean-opinion equation, producing a downward shift in the trajectories before they converge to zero.

These results highlight the asymmetry between the two regimes. In Regime~1, the finite accumulated effective compromise $\mathcal{I}_\infty<+\infty$ preserves memory of both $m_v(0)$ and the initial correlation. In Regime~2, the divergence $\mathcal{I}(t)\to+\infty$ progressively erases both, yielding a universal limit $m_v^\infty=m_u$ independent of initial conditions.

\subsection{Test~5. Long-time profiles of the opinion--conviction distribution}

In this test, we investigate the long-time behavior of the distribution $f(v,s,u,t)$ in order to characterize its asymptotic profile in Regimes~1, 2, and~3, using the payoff matrices~\eqref{J_test}. For each regime, we examine both the joint marginal distribution $\rho(v,u,T)$ and the marginal opinion distribution $h(v,T)$. The former describes the asymptotic coupling between expressed opinions and private convictions; the latter shows the resulting observable distribution of opinions and provides a direct comparison with the theoretical prediction.

\paragraph{Setup.}
The initial opinions are uniformly distributed over $I$, while the initial convictions follow a symmetric Beta distribution with parameters $(0.5,0.5)$, rescaled onto $I$. The initial strategy is uniformly distributed over $\Delta_1$. Thus, the initial joint distribution factors as $f_0=h_0(v)g_0(s)k(u)$ with
\begin{equation}\label{eq:ic}
f_0(v,s,u)
=\frac{\mathbf{1}_{I}(v)}{2}
\cdot \mathbf{1}_{\Delta_1}(s) \cdot\frac{
1}{\pi\sqrt{1-u^2}}.
\end{equation} 

Both $s$ and $u$ are initialized independently, so that $\Sigma_{su}(0)=0$. In particular, under constant payoffs the factorization of the strategy--conviction marginal is preserved for all $t\geq0$, and $\Sigma_{su}(t)=0$ throughout the simulation. Since $m_v(0)=0=m_u$, the mean-opinion equation implies $m_v(t)=0$ throughout the evolution, and hence $m_v^\infty=0$ in all three regimes.

The evolution is computed up to the final time $T=50$, which is large enough to observe the long-term behavior of the system.

\paragraph{Theoretical prediction for $h^\infty(v)$.}
Under factorization $\eta^\infty(s,u)=g^\infty(s)k(u)$, the asymptotic marginal opinion distribution is
\begin{equation}\label{eq:hinfty-integral}
    h^\infty(v)
    =
    \iint_{\Delta_1\times I}
    \rho^\infty(v\mid s,u)\,g^\infty(s)\,k(u)\,\mathrm{d}s\,\mathrm{d}u,
\end{equation}
where $\rho^\infty(v\mid s,u)$ is the Gibbs density \eqref{eq:rho_cond_Gibbs} with attraction center $V_{\infty}(s,u)=(1-s)u$, since $m_v^\infty=0$. The theoretical profile is obtained from \eqref{eq:hinfty-integral} in each regime, as detailed below.

\begin{itemize}
\item\emph{Regime~1:} $g^\infty=\delta(s-1)$, so $V_{\infty}(1,u)=0$ for all $u$, and $h^\infty(v)=\rho^\infty(v\mid V_{\infty}=0)$, the single Gibbs density with center~$0$. For $D(v)=1-v^2$:
\begin{equation}\label{eq:hinfty-regime1}
    h^\infty(v)
    \propto
    \frac{1}{(1-v^2)^2}
    \exp\!\left(-\frac{v^2}{\nu^2(1-v^2)}\right).
\end{equation}

\item\emph{Regime~2:} $g^\infty=\delta(s)$, so $V_{\infty}(0,u)=u$, and $h^\infty(v)=\int_{I}\rho^\infty(v\mid V_{\infty}=u)k(u)\,du$
is a continuous mixture of Gibbs densities indexed by conviction.

\item\emph{Regime~3:} $g^\infty=(\Phi_{\mathcal{T}_\infty})_\#g_0$
with $\hat m_s=0.6$, and
\begin{equation}\label{eq:hinfty-regime3}
    h^\infty(v)
    =
\iint_{\Delta_1\times I}\rho^\infty(v\mid V_{\infty}=(1-s)u)\,
    g^\infty(s)\,k(u)\,ds\,du,
\end{equation}
where $\mathcal{T}_\infty$ is obtained from the numerical
solution of~\eqref{eq:cumulative-gain-equation-two}.
\end{itemize}
\begin{figure}[ht]
\centering

\begin{subfigure}[b]{\textwidth}
    \centering
    \includegraphics[width=\linewidth]{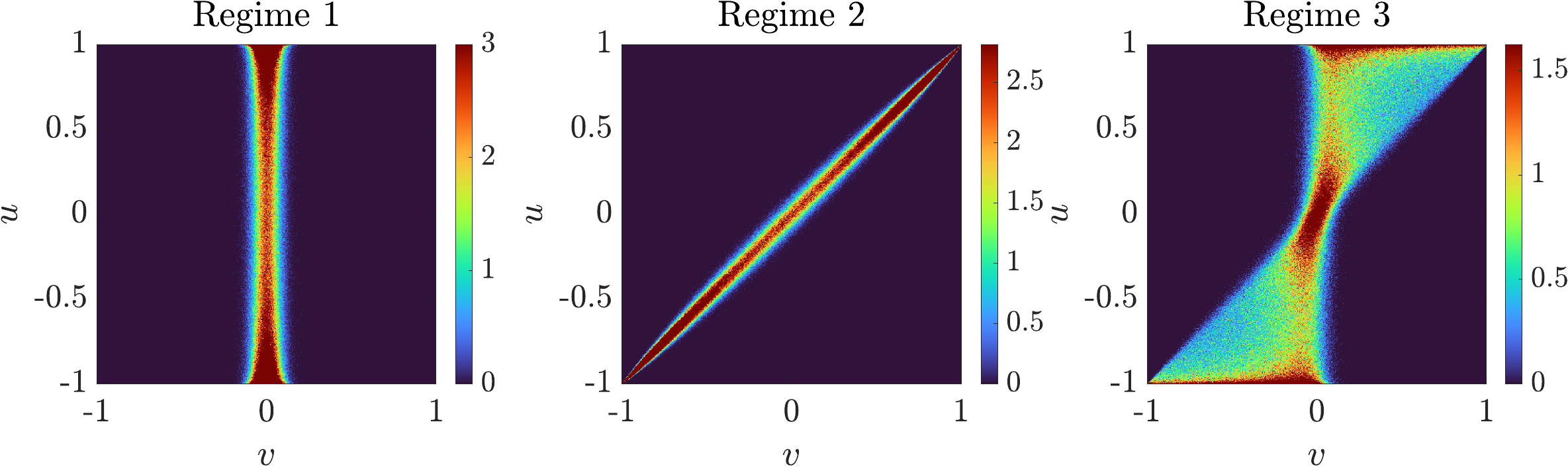}
    \caption{Joint  distribution $\rho(v,u,T)$ at $T=50$.}
    \label{fig:test3a}
\end{subfigure}

\begin{subfigure}[b]{\textwidth}
    \centering
    \includegraphics[width=0.9\linewidth]{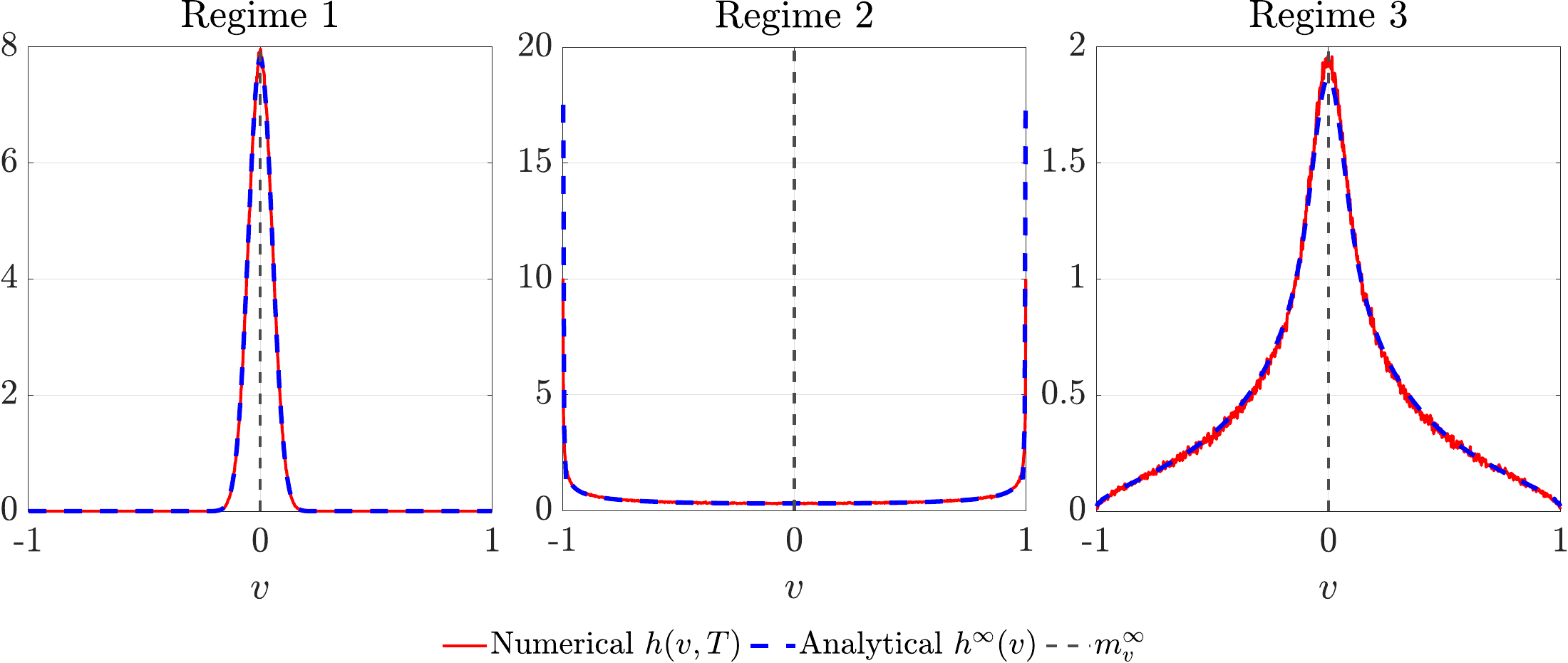}
    \caption{Marginal opinion distribution $h(v,T)$ at
    $T=50$. Red solid: empirical distribution from DSMC.
    Blue dashed: theoretical 
    $h^\infty(v)$  from~\eqref{eq:hinfty-integral}.}
    \label{fig:test3b}
\end{subfigure}

\caption{Test~5. Long-time opinion--conviction distributions for Regimes~1, 2, and~3.}
\label{fig:test3}
\end{figure}

\paragraph{Numerical results.}
Figure~\ref{fig:test3a} shows the joint distribution $\rho(v,u,T)$ at $T=50$, and Figure~\ref{fig:test3b} compares the empirical $h(v,T)$ with the theoretical prediction $h^\infty(v)$.

In \emph{Regime~1}, the joint distribution concentrates in a vertical band around $v=0$, independently of $u$: since $g^\infty=\delta(s-1)$, all agents are attracted to the same opinion target $V_{\infty}=0$ regardless of conviction. The marginal $h(v,T)$ is well approximated by~\eqref{eq:hinfty-regime1}, which is symmetric and concentrated near $v=0$ with tails suppressed by the exponential factor.

In \emph{Regime~2}, the joint distribution develops a ridge around the diagonal $v=u$: since $g^\infty=\delta(s)$, each agent is attracted to their own conviction, while the nonzero diffusion maintains a spread around the corresponding attraction center. The marginal $h(v,T)$ approximates $\int_{I}\rho^\infty(v|u)k(u)\,du$, which inherits the U-shaped character of the $\mathrm{Beta}(\tfrac{1}{2},\tfrac{1}{2})$ conviction distribution: opinion mass is enhanced near $v\approx\pm1$.

In \emph{Regime~3}, we have $\hat m_s=0.6$, so that $V_{\infty}(\hat m_s,u)=(1-\hat m_s)u=0.4u$ provides a representative average attraction center. However, since $g^\infty$ retains a nontrivial distribution, the individual centers $V_\infty(s,u)=(1-s)u$ vary with $s$.The joint distribution is intermediate between the two dominance regimes: broader than in Regime~1 and less U-shaped than in Regime~2, mass is broadly spread over $I^2$, with some concentration around the central region and a fan-like extension toward larger values of $|u|$. The numerical distribution $h(v,T)$ provides a good approximation of theoretical profile~\eqref{eq:hinfty-regime3}.

The close agreement between the empirical and theoretical profiles in all three regimes supports the predicted asymptotic structure \eqref{eq:finfty-factored} and the marginal formula \eqref{eq:hinfty-integral}.

\subsection{Test~6. State-dependent payoffs: geometry of the segregation boundary}

We investigate the long-time behavior of the two-strategy model with the state-dependent payoffs of Section~\ref{sec:state-dependent}, focusing on the two qualitatively distinct geometries of the segregation boundary~\eqref{eq:hyperbola} corresponding to $K(T)<0$ and $K(T)>0$. We consider two parameter configurations, denoted Case~A and Case~B, which share the same initial distribution and all model parameters except the cross-strategy costs $c_{CS}$ and $c_{SC}$. This isolates the effect of the cost asymmetry on the strategic segregation and on the geometry of the gain-induced regions in the $(v,u)$-plane.

\paragraph{Setup.}
Both cases use $a=b=1$ and the initial distribution $f_0$ is given by \eqref{eq:ic} with $v\sim\mathrm{Unif}(I)$, $s\sim\mathrm{Unif}(\Delta_1)$, and $u\sim\mathrm{Beta}(\tfrac12,\tfrac12)$ on $I$, giving $m_v(0)=m_u=0$ and $\Sigma_{su}(0)=0$.
The cross-strategy costs are:
\[
    \text{Case~A:}\quad c_{CS}=0.1,\;c_{SC}=0.05;
    \qquad
    \text{Case~B:}\quad c_{CS}=0.05,\;c_{SC}=0.5.
\]
In Case~A the cost of compromise-against-stubborn exceeds that of stubborn-against-compromise ($c_{CS}>c_{SC}$), while in Case~B the reverse holds.  
We consider a sufficiently large final time $T=1000$, chosen so that the transient dynamics has decayed and the numerical solution can be considered representative of the asymptotic regime.

\paragraph{Numerical results.}

At $T=1000$, the strategy marginal $g(s,T)$ is numerically concentrated near $s=0$ and $s=1$ with no appreciable interior mass, consistent with the segregated stationary structure~\eqref{eq:stationary-mixed}.

\begin{figure}[t]
\centering
\begin{subfigure}[b]{0.5\textwidth}
\centering
\includegraphics[width=\linewidth]{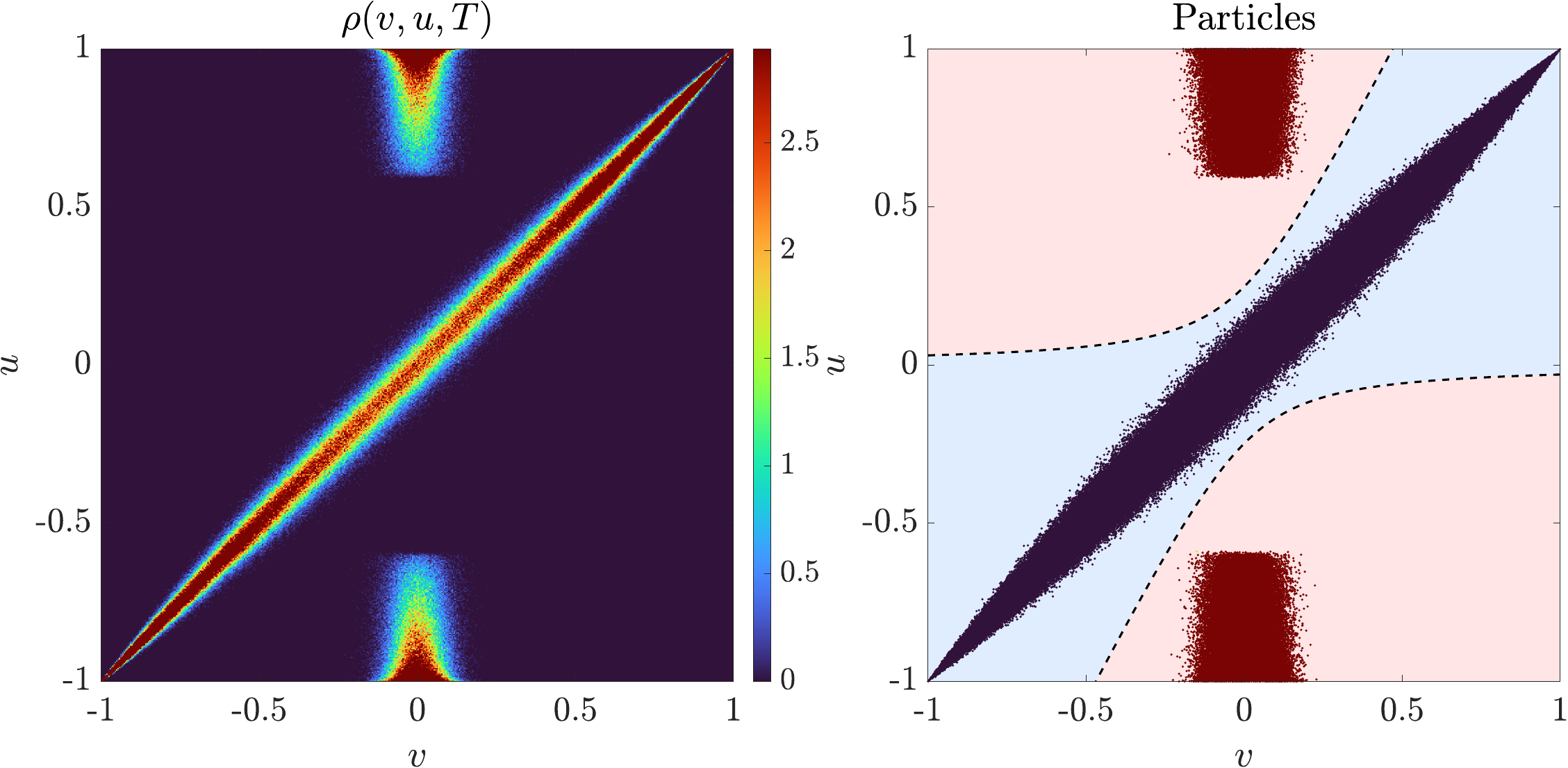}
\caption{Case A  at $T=1000$.}
\label{fig:test4a}
\end{subfigure}\begin{subfigure}[b]{0.5\textwidth}
\centering
\includegraphics[width=\linewidth]{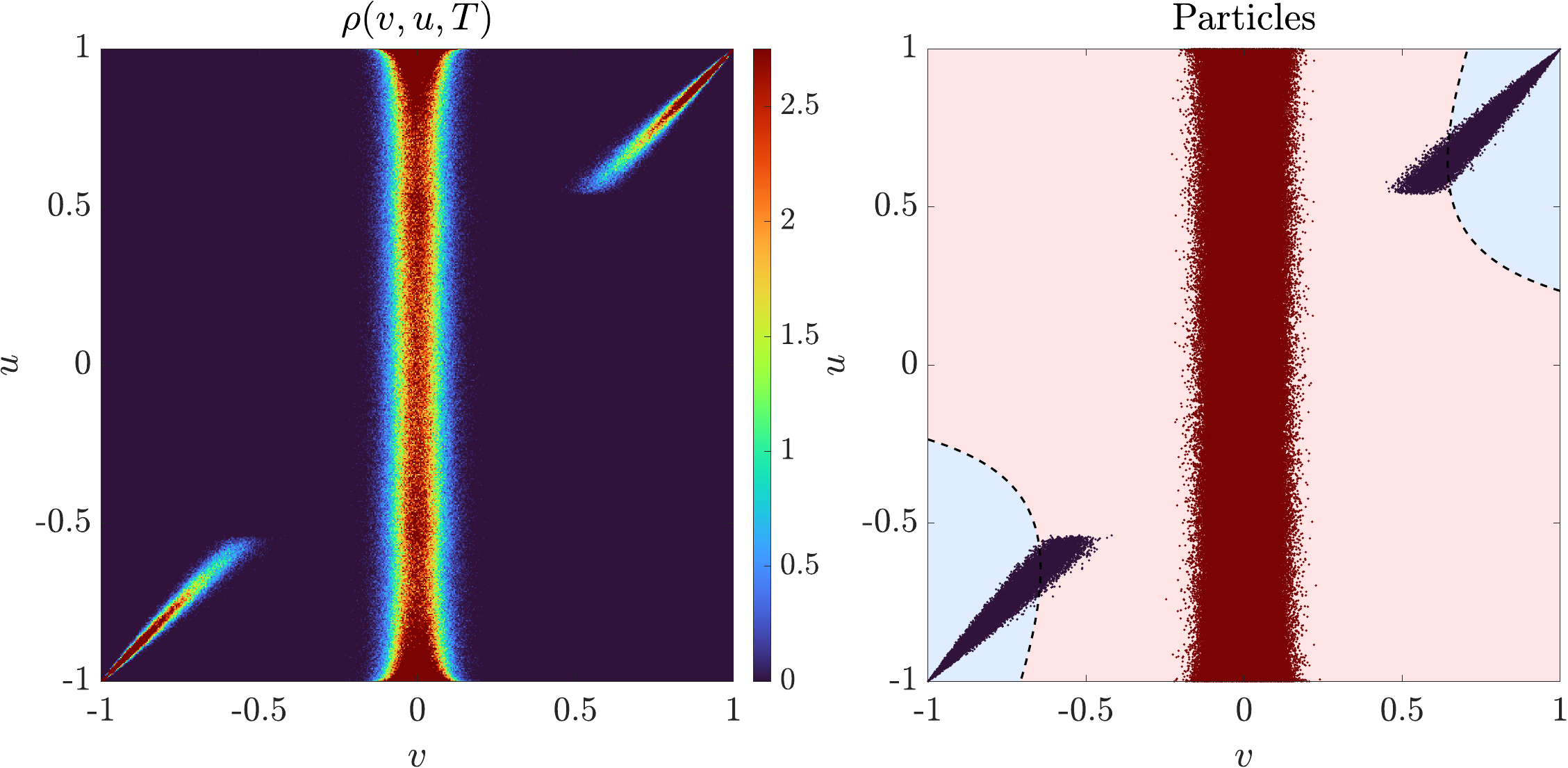}
\caption{Case  B at $T=1000$.}
\label{fig:test4b}
\end{subfigure}
\caption{Test 6.  Long-time state-dependent payoff dynamics. Left panel of each subfigure: joint marginal $\rho(v,u,T)$. Right panel: particle scatter colored by final strategy ($\pi_C$ red, $\pi_S$ blue); dashed curve: segregation boundary $\{\mathcal{G}[f](v,u,T)=0\}$; light red/blue: regions $\mathcal{R}_C(T)$ / $\mathcal{R}_S(T)$. }
\label{fig:test4}
\end{figure}

For each case, Figure~\ref{fig:test4} shows two complementary representations of the final state: the joint marginal $\rho(v,u,T)$ in the opinion--conviction plane (left panel), and a particle scatter plot colored by final strategy ($\pi_C$ red, $\pi_S$ blue), with the instantaneous segregation boundary $\{\mathcal{G}[f](v,u,T)=0\}$ overlaid as a black dashed curve and the regions $\mathcal{R}_C(T)$ and $\mathcal{R}_S(T)$ shaded in light red and light blue (right panel).

\emph{Case~A.}
At $T=1000$, the shift satisfies $K(T)<0$, so the segregation boundary is a hyperbola opening in the $v-u$ direction with the stubbornness-favoring region $\mathcal{R}_S(T)$ containing the diagonal $\{v=u\}$. This structure is clearly reflected in the particle distribution: agents near the diagonal adopt $\pi_S$ (blue), while those far from the diagonal adopt $\pi_C$ (red). The two compromise-dominated clusters are centered near $v\approx0$ with convictions $u\in[-1,-0.6]\cup[0.6,1]$, lying well outside the diagonal. This is consistent with the gain $(a=b=1)$
\[
    \mathcal{G}[f](v,u,T)
    =
    (v-u)^2-(v-m_v(T))^2+K(T),
\]
for agents near $v\approx m_v(T)\approx0$ with $u$ far from $v$, the term $(v-u)^2$ is large and dominates the negative $K(T)$, giving positive gain for $\pi_C$. Conversely, near the diagonal $u\approx v$, the term $(v-u)^2$ is small and $K(T)<0$ makes the gain negative, favoring $\pi_S$.

\emph{Case~B.}
At $T=1000$, $K(T)>0$, so the hyperbola opens in the $v-m_v(T)$ direction, and the compromise-favoring region forms a central strip around the population mean. For $a=b=1$, the strip covers all $v$ satisfying $|v-m_v(T)|<\sqrt{K(T)}$; numerically $\sqrt{K(T)}\approx0.64$, so agents with $|v|\lesssim0.64$ favor $\pi_C$ for all $u\in I$. The particle distribution confirms this: a vertical band centered at $v\approx0$ is predominantly red ($\pi_C$), while the outer regions $|v|\gtrsim0.64$ are predominantly blue ($\pi_S$), particularly near the diagonal $\{u\approx v\}$ where $(v-u)^2$ is small. A small number of blue particles ($\pi_S$) appear inside the theoretical $\mathcal{R}_C(T)$ region. This does not contradict the theory: $\mathcal{R}_C(T)$ describes the direction of the instantaneous strategic drift, but particles that have already reached $s=0$ experience zero strategic flux (since $s(1-s)=0$ at $s=0$) and cannot re-enter the interior even if the gain at their state is positive.

The two cases thus illustrate the two geometries predicted in Section~\ref{sec:state-dependent}: a negative $K$ produces a stubbornness-favoring band around the diagonal, while a positive $K$ produces a compromise-favoring strip around the population mean.

\subsection{Test~7. Asymptotic dynamics of the three-strategy model}
\label{sec:test5}

We consider four representative constant payoff matrices to illustrate the main qualitatively distinct asymptotic behaviors of the three-strategy kinetic model: dominance of the moderate strategy $\pi_2$, stable interior coexistence, persistent oscillations of RPS type, and heteroclinic boundary dynamics. Stable interior coexistence also occurs in the two-strategy model, as shown in Regime~3 of Section~\ref{sec:2s-constant}. By contrast, dominance of an intermediate strategy, persistent RPS-type oscillations, and heteroclinic cycling have no direct analogue in the two-strategy model. The purpose of the test is to validate the push-forward solution~\eqref{eq:strategy-projective-pushforward} and to relate the strategic dynamics to the resulting evolution of the opinion distribution.

Case~A illustrates $\pi_2$-dominance, chosen over $\pi_1$ or $\pi_3$ dominance because the latter two reproduce profiles already characterized in Section~\ref{sec:2s-constant}
(Regimes~1 and~2), whereas $\pi_2$-dominance yields
$m_\Theta^\infty=\theta_2=\tfrac 12$ and a slope-$\tfrac12$ opinion--conviction profile with no two-strategy analogue.

\paragraph{Setup.}
The initial strategy distribution is
$\mathrm{Dirichlet}(1,5,1)$ on $\Delta_2$, with  initial means $m_s(0)=1/7$, $m_r(0)=5/7$, $m_q(0)=1/7$. This asymmetric initialization places the initial mass predominantly near the vertex $\pi_2$. The initial opinions are uniform on $I$ and the convictions follow a $\mathrm{Beta}(1/2,1/2)$ distribution rescaled to $I$. The initial joint distribution factors as
\begin{equation*}
    f_0(v,s,r,u)
    = \frac{\mathbf{1}_I(v)}{2}\cdot 30r^4\mathbf{1}_{\Delta_2}(s,r)\cdot
    \frac{\mathbf{1}_I(u)}{\pi\sqrt{1-u^2}}.
\end{equation*}
Since $s$, $r$, and $u$ are initialized independently, we have $\Sigma_{su}(0)=\Sigma_{ru}(0)=0$. Under constant payoffs this factorization is preserved for all $t\geq0$ (Corollary~\ref{cor:eta-pushforward}), so $\Sigma_{\Theta u}(t)\equiv0$. Since $m_v(0)=m_u=0$, the mean-opinion equation gives $m_v(t)=0$ for all $t\ge0$, including in the oscillatory and heteroclinic regimes. The final time is $T=100$.

\paragraph{Payoff matrices and theoretical predictions.} The four payoff matrices and the corresponding relative gains are:

\textbf{Case~A -- $\pi_2$ dominance.}
\begin{equation*}
    J^{(A)}=\begin{pmatrix}1&1&1\\3&3&3\\0&0&0\end{pmatrix},
    \qquad
    G_1\equiv1>0,
    \qquad
    G_2\equiv3>0.
\end{equation*}
Both gains are strictly positive on $\Delta_2$, so $\pi_2$ strictly dominates. Since $G_2>G_1$, the cumulative gain satisfies $\mathcal{T}_2(t)\gg\mathcal{T}_1(t)\to+\infty$, giving $g^\infty=\delta(s)\delta(r-1)$, $(m_s^\infty,m_r^\infty)=(0,1)$, and $m_\Theta^\infty=1/2$.

\textbf{Case~B -- Stable interior coexistence.}
\begin{equation*}
    J^{(B)}=\begin{pmatrix}0&-1&2\\2&0&-1\\-1&2&0\end{pmatrix},
    \qquad
    \begin{aligned}
        G_1(m_s,m_r)&=2-m_s-5m_r,\\
        G_2(m_s,m_r)&=-1+4m_s-m_r.
    \end{aligned}
\end{equation*}
The equal-payoff conditions~\eqref{G120} yield the unique interior equilibrium $(\hat m_s,\hat m_r)=(1/3,1/3)$, which is a stable focus for the classical replicator dynamics. The cumulative gains $(\mathcal T_1,\mathcal T_2)$ converge to a finite limit $\mathcal T^\infty\in\mathbb R^2$, so $g^\infty=(\Phi_{\mathcal T^\infty})_\#g_0$
with $(m_s^\infty,m_r^\infty)=(1/3,1/3)$ and
$m_\Theta^\infty=1/2$.

\textbf{Case~C -- Conservative oscillations (RPS).}
\begin{equation*}
    J^{(C)}=\begin{pmatrix}0&-1&1\\1&0&-1\\-1&1&0\end{pmatrix},
    \qquad
    \begin{aligned}
        G_1(m_s,m_r)&=1-3m_r,\\
        G_2(m_s,m_r)&=-1+3m_s.
    \end{aligned}
\end{equation*}
The matrix $J^{(C)}$ is zero-sum ($J^{(C)}+(J^{(C)})^\top=0$), which implies that the interior equilibrium $(\hat m_s,\hat m_r)=(\tfrac13,\tfrac13)$ is a center and the classical replicator trajectories form closed orbits. For the kinetic model, the numerical solution of the cumulative-gain equations exhibits a closed orbit with estimated period $T_p\approx16$, for the chosen initial distribution,  and $g(s,r,t)$ oscillates periodically without converging.

\textbf{Case~D -- Heteroclinic cycle.}
\begin{equation*}
    J^{(D)}=\begin{pmatrix}0&-2&1\\1&0&-2\\-2&1&0\end{pmatrix},
    \qquad
    \begin{aligned}
        G_1(m_s,m_r)&=1+m_s-4m_r,\\
        G_2(m_s,m_r)&=-2+5m_s+m_r.
    \end{aligned}
\end{equation*}
The interior equilibrium $(\hat m_s,\hat m_r)=(1/3,1/3)$ is an unstable focus for the classical replicator dynamics. Each vertex of $\Delta_2$ is unstable in exactly one direction toward the next in the cyclic ordering $\pi_1\to\pi_2\to\pi_3\to\pi_1$, forming a heteroclinic cycle on
$\partial\Delta_2$. The projective weights spiral outward, and the alternation of near-dominance phases has periods of increasing duration.

\begin{figure}[t]
\centering
\begin{subfigure}{\textwidth}
\centering
\includegraphics[width=\linewidth]{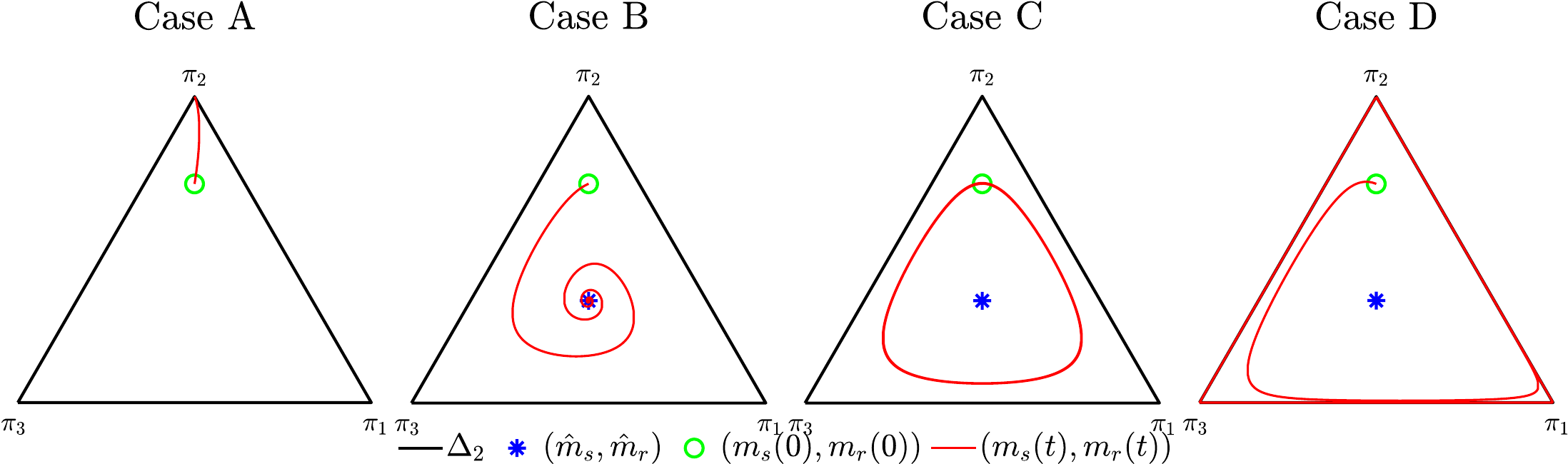}
\caption{Trajectories $(m_s(t),m_r(t))$ on $\Delta_2$ for Cases~A--D (left to right). Green marker: initial condition $(m_s(0),m_r(0))=(1/7,5/7)$. Blue star: interior equilibrium $(1/3,1/3)$ (Cases~B, C, D). Red curve: trajectory. Vertices: $\pi_1$ (bottom right), $\pi_2$ (top), $\pi_3$ (bottom left).}
\label{fig:test5a}
\end{subfigure}
\begin{subfigure}{\textwidth}
\centering
\includegraphics[width=\linewidth]{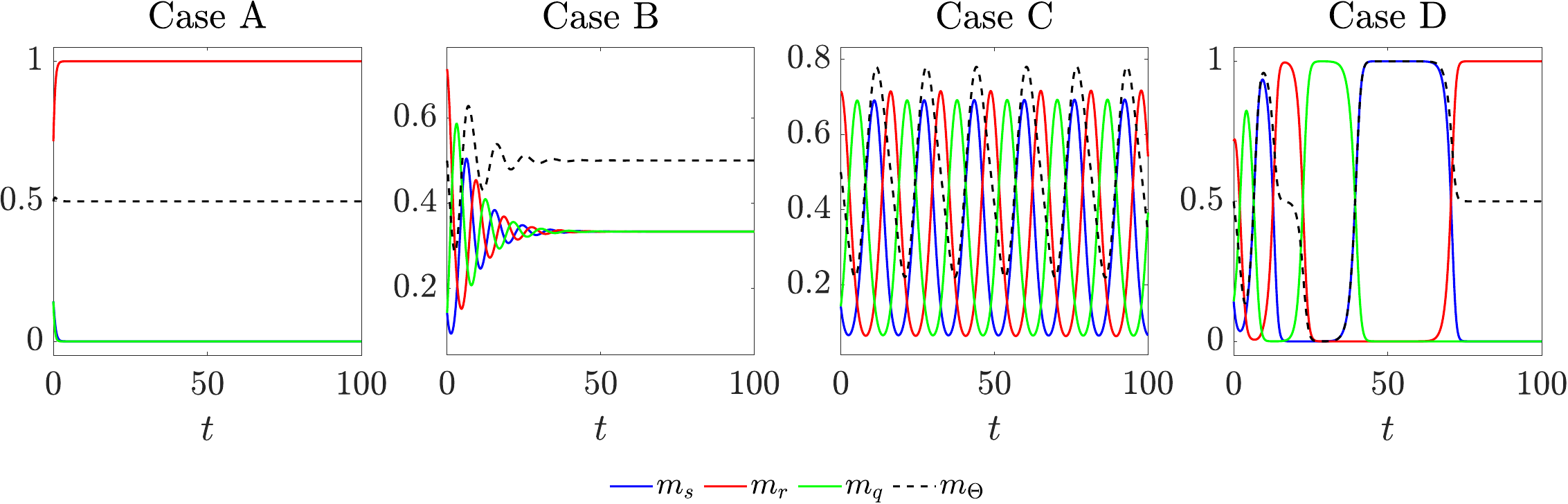}
\caption{Time evolution of $m_s(t)$, $m_r(t)$, $m_q(t)=1-m_s-m_r$, and $m_\Theta(t)=m_s+m_r/2$ for Cases~A--D (left to right).}
\label{fig:test5b}
\end{subfigure}
\caption{Test~7. Strategic dynamics of the three-strategy model. The trajectories and time-dependent strategy means highlight the different asymptotic behaviors in Cases~A--D.}
\label{fig:test5-strategy-overview}
\end{figure}

\begin{figure}[t]
\centering
\includegraphics[width=\linewidth]{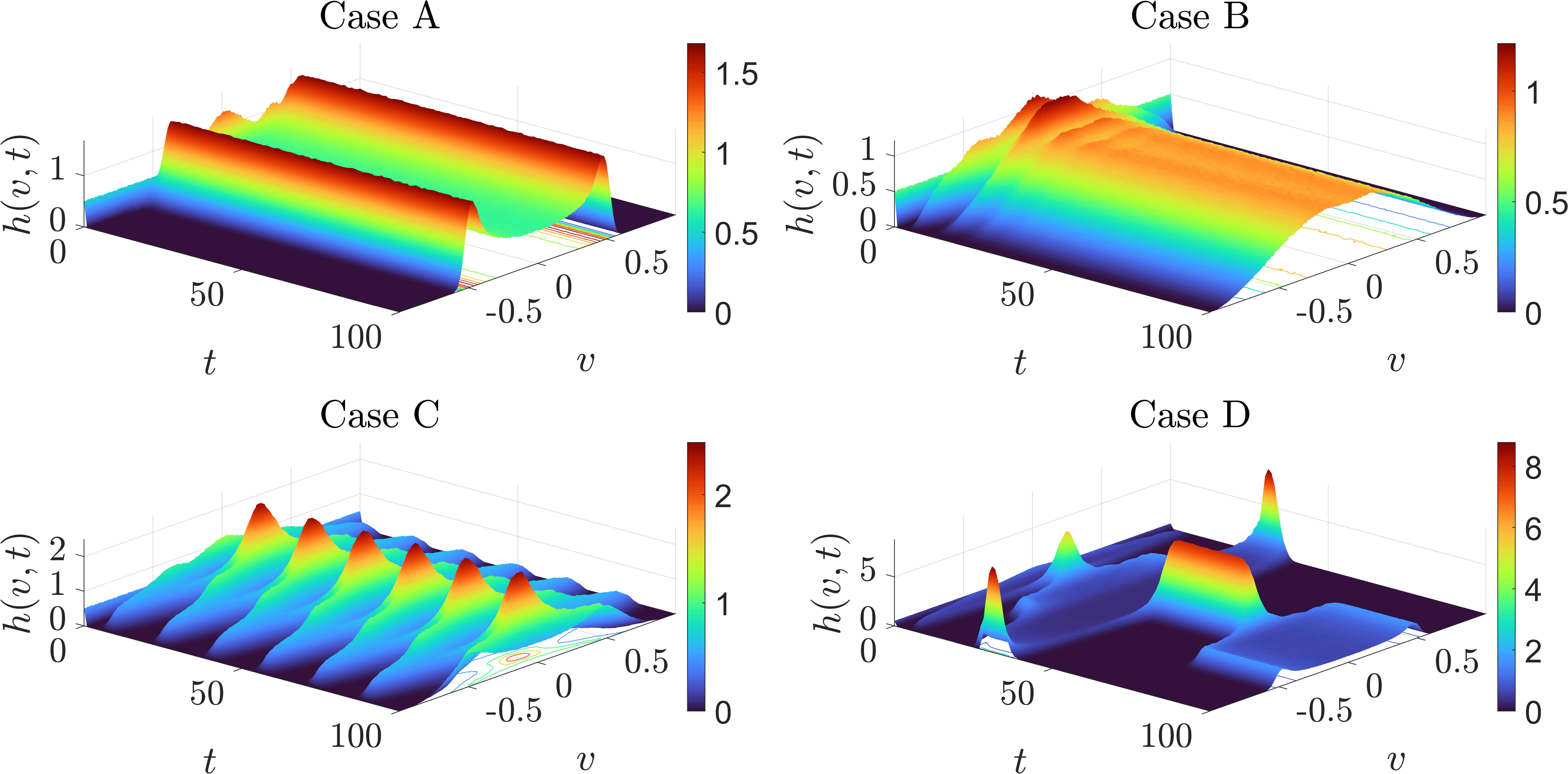}
\caption{Test~7. Surface $h(v,t)$ for Cases~A--D. Case~A: monotone convergence to a profile centered near $v=0$. Case~B: damped convergence to a stationary intermediate profile. Case~C: periodic oscillation between narrow (large $m_\Theta$) and broad (small $m_\Theta$) profiles. Case~D: large-amplitude intermittent oscillations reflecting the heteroclinic dynamics.}
\label{fig:test5d}
\end{figure}

\begin{figure}[t]
\centering
\begin{subfigure}{\textwidth}
\centering
\includegraphics[width=\linewidth]{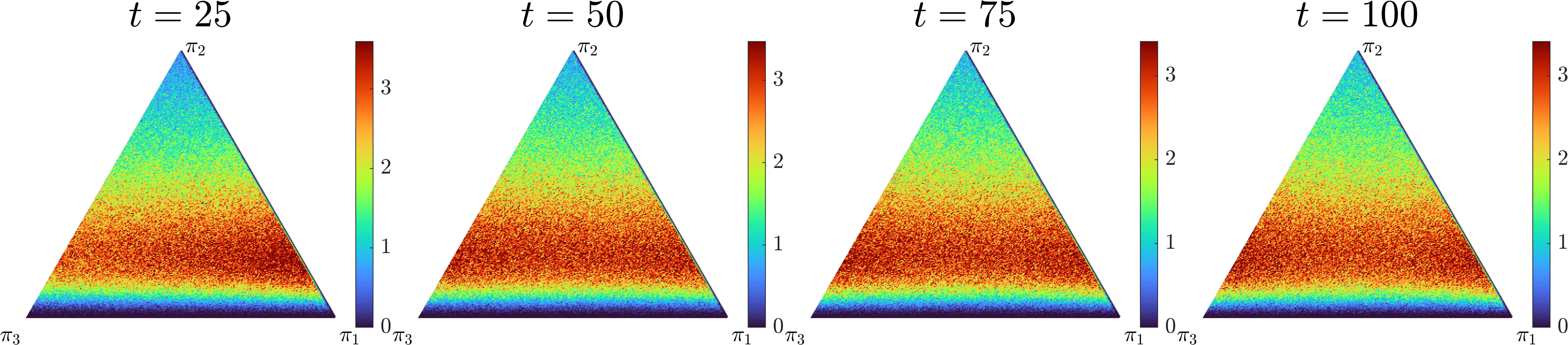}
\caption{Case~B. Empirical density of $g(s,r,t)$ on $\Delta_2$ at $t=25,50,75,100$,  converging to the non-degenerate profile $g^\infty=(\Phi_{\mathcal T^\infty})_\#g_0$
with $m_s^\infty=m_r^\infty=1/3$.}
\label{fig:test5cB}
\end{subfigure}
\begin{subfigure}{\textwidth}
\centering
\includegraphics[width=0.6\linewidth]{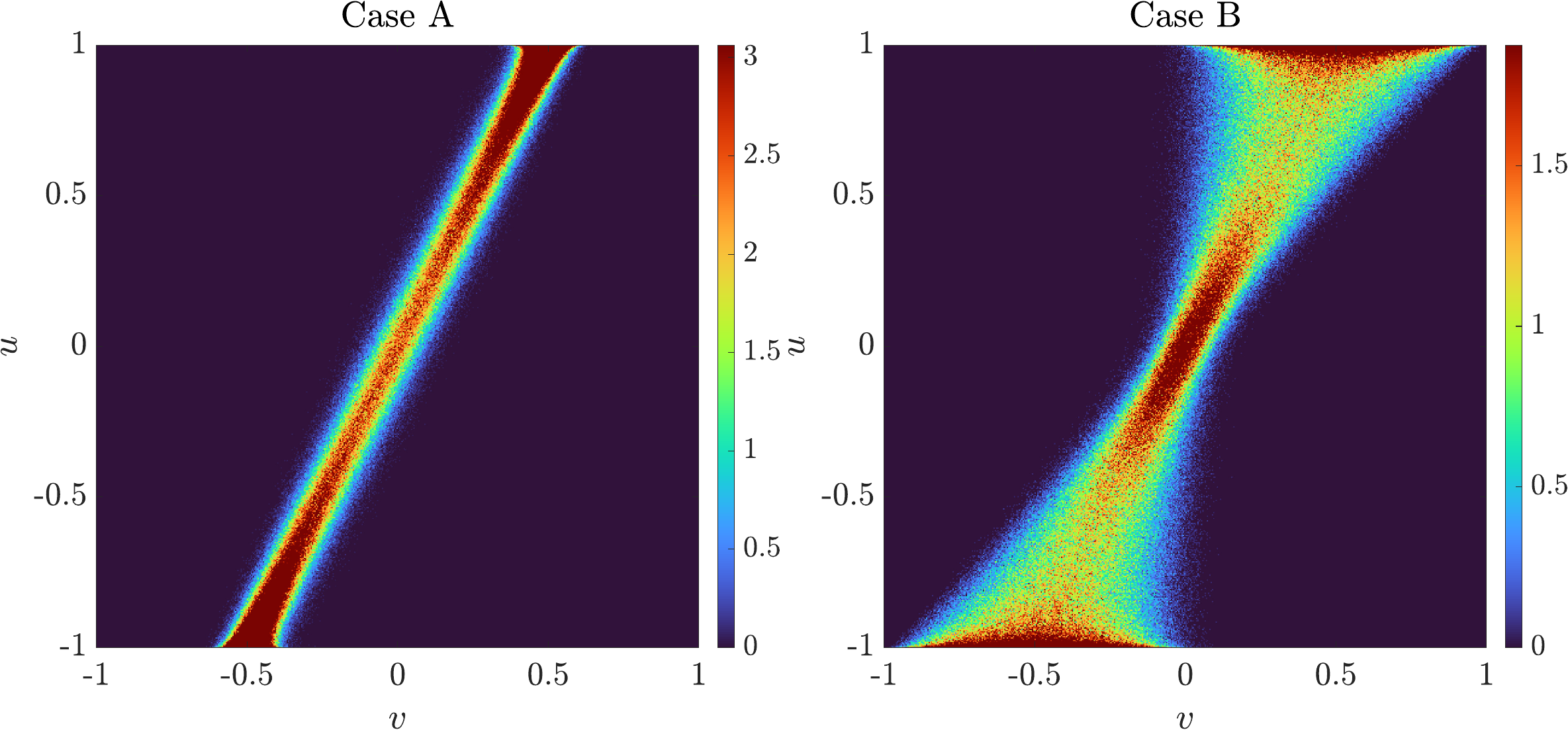}
\caption{Joint marginal $\rho(v,u,T)$ at $T=100$ for
Cases~A (left) and B (right). In Case~A, mass concentrates near
the line $v=u/2$ (dashed), reflecting the attraction center
$V_{\infty}(0,1,u)=u/2$. In Case~B, the distribution is more diffuse,
reflecting the non-degenerate mixture of Gibbs profiles induced
by the non-degenerate $g^\infty$.}
\label{fig:test5eAB}
\end{subfigure}
\caption{Test~7. Asymptotic behavior in the convergent cases: the non-degenerate strategy profile in Case~B and the corresponding opinion--conviction distributions in Cases~A--B.}
\label{fig:test5-convergent}
\end{figure}

\begin{figure}[t]
\centering
\begin{subfigure}{\textwidth}
\centering
\includegraphics[width=\linewidth]{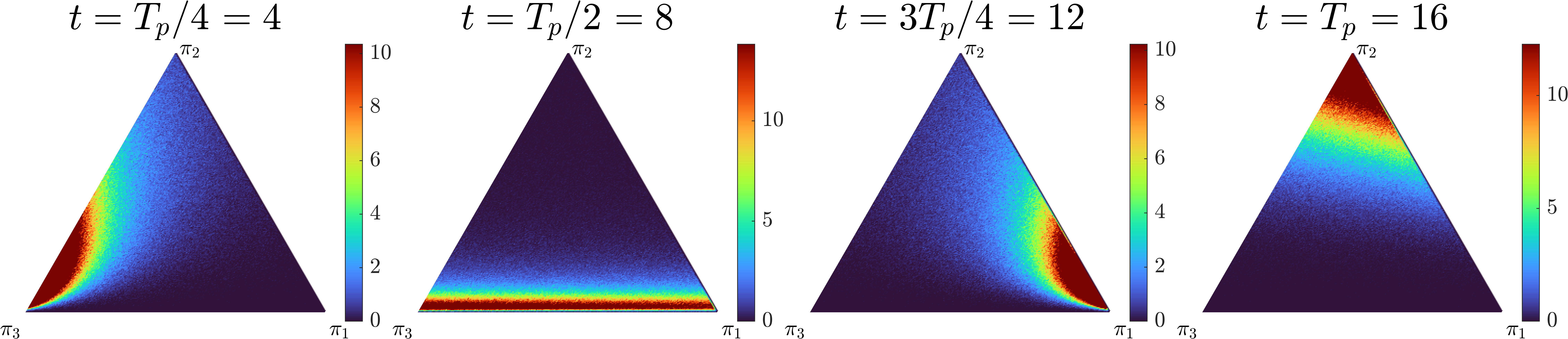}
\caption{Empirical density of $g(s,r,t)$ on $\Delta_2$ at $t=T_p/4,T_p/2,3T_p/4,T_p$ with. The distribution rotates periodically around $(1/3,1/3)$ without converging.}
\label{fig:test5cC}
\end{subfigure}
\begin{subfigure}{\textwidth}
\centering
\includegraphics[width=\linewidth]{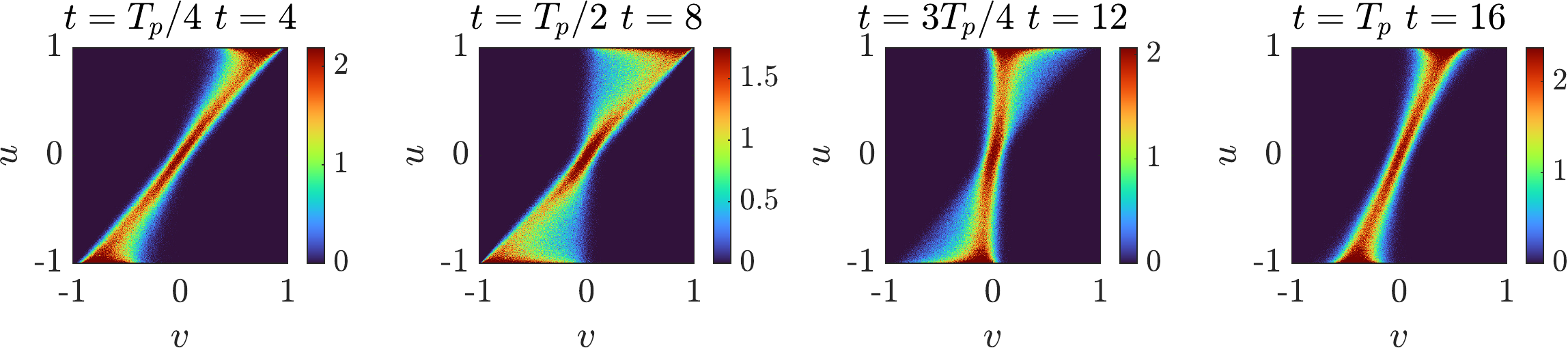}
\caption{Snapshots of $\rho(v,u,t)$ at the same times. The opinion--conviction pattern rotates
periodically: near $\pi_2$ the mass clusters along $v=u/2$; near $\pi_3$ along the diagonal $v=u$; near $\pi_1$ in a vertical band.}
\label{fig:test5eC}
\end{subfigure}
\caption{Test~7, Case~C. Periodic strategy dynamics and the corresponding opinion--conviction patterns, with $T_p\approx16$.}
\label{fig:test5-periodic}
\end{figure}

\begin{figure}[t]
\centering
\begin{subfigure}{\textwidth}
\centering
\includegraphics[width=0.5\linewidth]{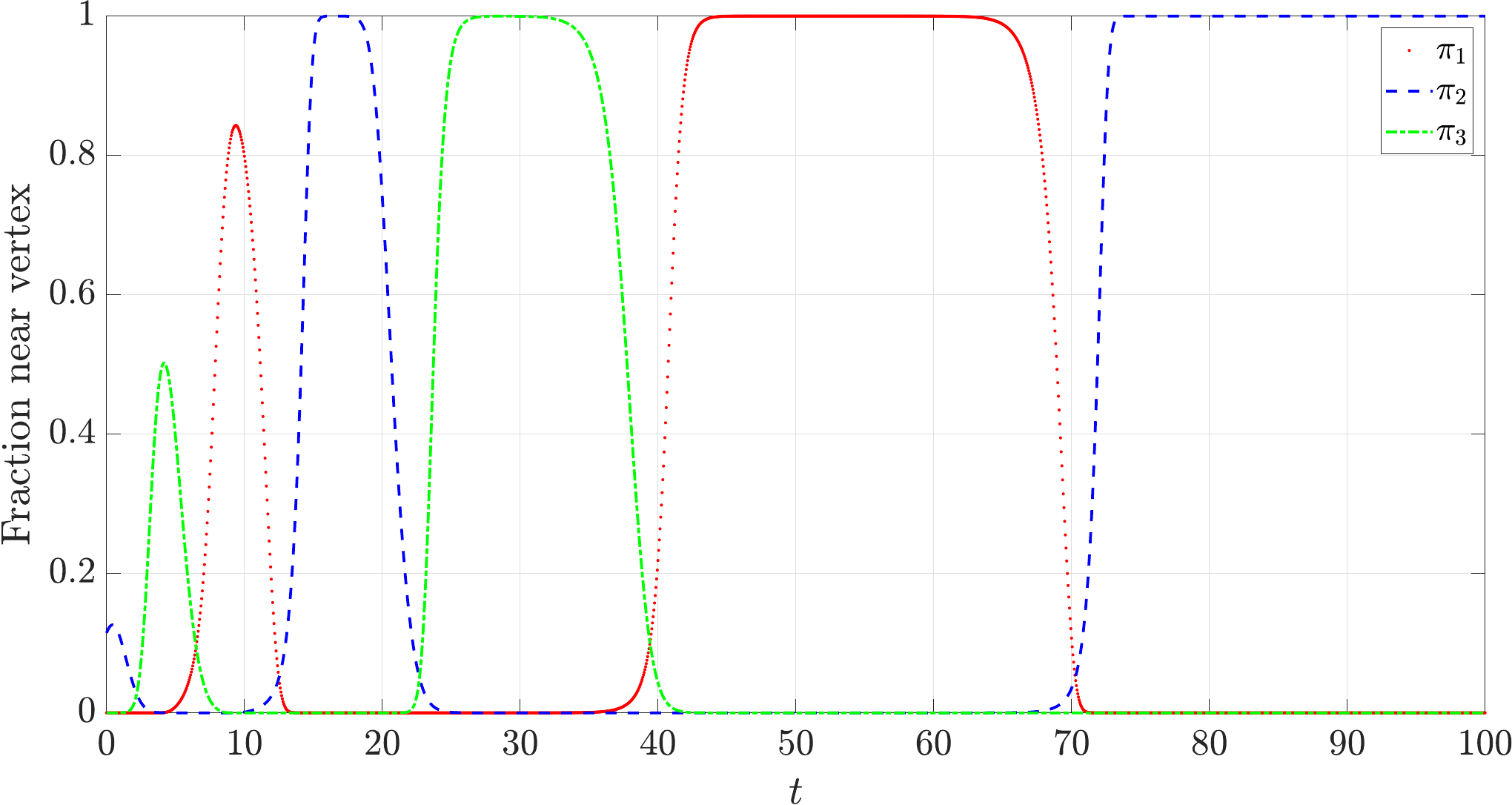}
\caption{Fraction of particles within $\epsilon=0.1$ of each vertex: near $\pi_1$ (red), $\pi_2$ (blue), $\pi_3$ (green). Successive peaks of increasing width confirm the heteroclinic cycle $\pi_1\to\pi_2\to\pi_3\to\pi_1$.}
\label{fig:test5cD}
\end{subfigure}
\begin{subfigure}{\textwidth}
\centering
\includegraphics[width=\linewidth]{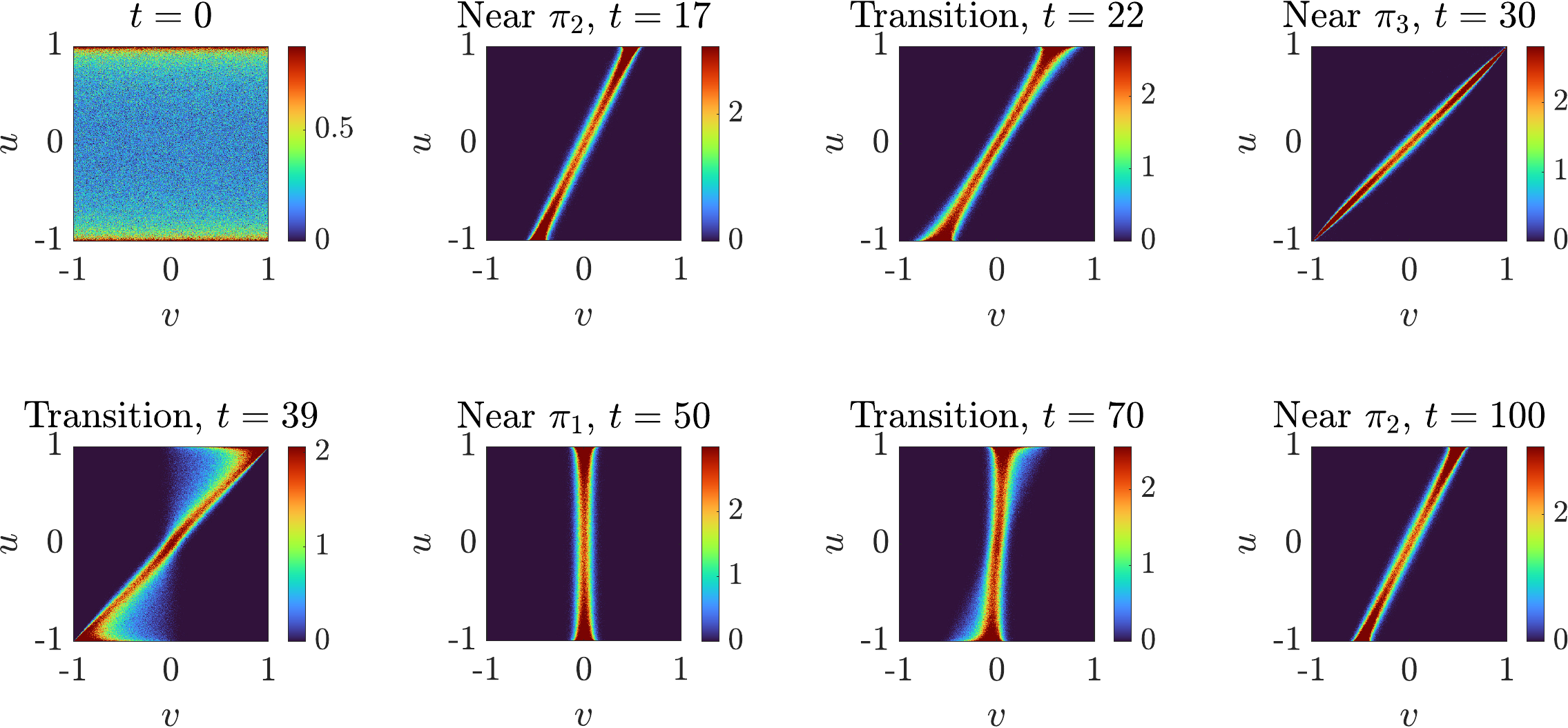}
\caption{Snapshots of $\rho(v,u,t)$ at $t=0$,
$17$, $22$, $30$, $39$, $50$, $70$, $100$, corresponding to the
initial state and to successive near-dominance and transition
phases. Each near-dominance phase reproduces the corresponding
pure-strategy pattern (vertical band for $\pi_1$, slope-$1/2$
line for $\pi_2$, diagonal for $\pi_3$); transition instants show intermediate geometries.}
\label{fig:test5eD}
\end{subfigure}
\caption{Test~7, Case~D. Heteroclinic dynamics and the alternation of opinion--conviction patterns during the successive near-dominance phases.}
\label{fig:test5-heteroclinic}
\end{figure}
\paragraph{Numerical results.}

Figure~\ref{fig:test5a} shows the trajectories $(m_s(t),m_r(t))$ on  $\Delta_2$. Figure~\ref{fig:test5b} reports the time evolution of $m_s$, $m_r$, $m_q=1-m_s-m_r$, and $m_\Theta=m_s+m_r/2$.     Figures~\ref{fig:test5cB} and~\ref{fig:test5cC} display the strategy distribution $g(s,r,t)$ on $\Delta_2$ at selected times for Cases~B and~C respectively. Figure~\ref{fig:test5cD} reports, for Case~D, the fractions of the population within $\epsilon=0.1$ of each vertex. Figure~\ref{fig:test5d} shows the surface $h(v,t)$ for all four cases. Figure~\ref{fig:test5eAB} shows the joint marginal $\rho(v,u,T)$ for Cases~A and~B at $T=100$; Figures~\ref{fig:test5eC} and~\ref{fig:test5eD} show snapshots of $\rho(v,u,t)$  for Cases~C and~D.

\emph{Case~A.}
The trajectory of $(m_s,m_r)$ converges monotonically to the vertex $(0,1)$ rapidly due to the asymmetric initial condition $m_r(0)=\tfrac57$. As $g^\infty=\delta(s)\delta(r-1)$, the attraction center reduces to $V_{\infty}(0,1,u)=\frac12 m_u+\frac12 u=\frac12 u$ (using $m_v^\infty=m_u=0$): the expressed opinion is attracted to the average of the population mean and the individual conviction. The joint distribution $\rho(v,u,T)$ in Figure~\ref{fig:test5eAB} is concentrated near the line $v=u/2$, confirming this intermediate coupling between opinion and conviction.

\emph{Case~B.}
The trajectory exhibits a damped spiral converging to $(1/3,1/3)$, which is the key qualitative distinction from the two-strategy coexistence regime where no spiral is possible. The strategy distribution $g(s,r,t)$ in Figure~\ref{fig:test5cB} evolves along a rotating transient and converges to the non-degenerate profile $g^\infty=(\Phi_{\mathcal T^\infty})_\#g_0$ with $m_s^\infty=m_r^\infty=1/3$. Although $m_\Theta^\infty=1/2$ coincides with Case~A, the non-degenerate $g^\infty$ produces a mixture of Gibbs densities with heterogeneous centers, yielding a broader opinion profile $h^\infty(v)$, as visible in Figure~\ref{fig:test5d}. This is confirmed by $\rho(v,u,T)$ in
Figure~\ref{fig:test5eAB}, which is more diffuse than in Case~A and
does not concentrate on a single line.

\emph{Case~C.}
The trajectory of $(m_s,m_r)$ follows a closed orbit around
$(1/3,1/3)$ with period $T_p\approx16$. The strategy distribution
$g(s,r,t)$ rotates periodically on $\Delta_2$, as shown in
Figure~\ref{fig:test5cC} at times $T_p/4$, $T_p/2$, $3T_p/4$,
$T_p$. The oscillating $m_\Theta(t)$ drives a corresponding
periodic oscillation of $h(v,t)$: the profile narrows when
$m_\Theta$ is large (agents move toward the population mean) and
broadens toward $k(u)$ when $m_\Theta$ is small (agents move toward
their convictions). Snapshots of $\rho(v,u,t)$ in
Figure~\ref{fig:test5eC} confirm this: at $t=T_p/4$ the
distribution concentrates near the $\pi_2$ vertex and the mass
clusters around $v=u/2$; at $t=T_p/2$ near $\pi_3$, it spreads
along the diagonal $v=u$; at $t=3T_p/4$ near $\pi_1$, it
concentrates in a vertical band. This periodic rotation of the
opinion--conviction pattern is the observable fingerprint of the
underlying RPS strategic dynamics and has no analogue in the
two-strategy model.

\emph{Case~D.}
The trajectory spirals outward from $(1/3,1/3)$ toward
$\partial\Delta_2$. The strategy fractions in
Figure~\ref{fig:test5b} exhibit intermittent dominance in the
ordering $\pi_1\to\pi_2\to\pi_3\to\pi_1$, with each phase lasting progressively long. Figure~\ref{fig:test5cD} quantifies this via
the vertex fractions, whose successive peaks of increasing width
confirm the heteroclinic structure. The opinion distribution
$h(v,t)$ in Figure~\ref{fig:test5d} undergoes large-amplitude
intermittent oscillations: it narrows sharply near $\pi_1$
($m_\Theta\approx1$, consensus-like) and broadens toward $k(u)$
near $\pi_3$ ($m_\Theta\approx0$, conviction-aligned).
Figure~\ref{fig:test5eD} shows snapshots of $\rho(v,u,t)$ at eight
selected times. Near a $\pi_i$-dominance phase the pattern matches
the corresponding pure-strategy attractor: near $\pi_1$,
$\rho(v,u,t)$ concentrates in a vertical band; near $\pi_2$, along
the line $v=u/2$; near $\pi_3$, along the diagonal $v=u$.
Transition instants show intermediate patterns. This three-way
alternation of opinion--conviction geometries is the most
distinctive novel feature of the three-strategy model.

\FloatBarrier
\section{Conclusions}
\label{sec:conclusions}

We introduced a kinetic model for opinion dynamics in which each agent carries a public opinion, a fixed private conviction and a mixed behavioral strategy that evolves through a payoff-driven replicator-type rule and sets the balance between compromise and persistence towards the conviction. From the binary interactions we derived, in the quasi-invariant regime, a Vlasov--Fokker--Planck equation for the joint distribution of opinions, convictions and strategies, and we studied it both beyond and within scale separation.

For constant payoffs the strategy marginal is the push-forward of its initial datum through a projective map driven by cumulative relative gains. The strategic part of the problem is therefore exactly solvable, without a monokinetic closure and at any ratio of the opinion and strategy time scales; the same holds for the joint strategy--conviction marginal, which gives exact access to the covariances that enter the mean opinion equation. Stationary strategy distributions are in general not Dirac masses and retain information on the initial heterogeneity. For two strategies this yields a complete classification into compromise dominance, stubbornness dominance, stable coexistence and bistability. When compromise dominates, all agents share an opinion profile independent of their conviction, centered at a mean opinion that retains memory of the initial mean and of the transient strategy--conviction correlation; when stubbornness dominates, each agent is attracted to its own conviction and the mean opinion relaxes to the mean conviction; in the coexistence regime the stationary opinion profile is a mixture of conditional Gibbs densities over the strategy distribution. A third, moderate strategy produces behaviors with no two-strategy analogue, namely interior coexistence, persistent rock-paper-scissors oscillations and heteroclinic cycles, whose signature is a periodic or intermittent rotation of the opinion--conviction patterns. For quadratic state-dependent payoffs the strategy marginal is no longer autonomous, but the sign of the payoff gain partitions the opinion--conviction plane through a time-dependent hyperbolic segregation boundary, and the admissible stationary states are concentrated on the pure strategies.

In the two regimes generated by the time scales, fast opinions lead to a conditional Gibbs quasi-equilibrium with a self-consistent mean opinion, in which the strategy--conviction marginal evolves slowly, while fast strategies lead to a reduced opinion dynamics that depends on the selected strategic profile only through its average compromise propensity. The comparison with the exact two-strategy analysis shows that the factorization closure behind the second reduction reproduces the stationary mean opinion and the opinion--conviction covariance, but not the stationary opinion variance when the limiting strategy distribution is non-degenerate: the error does not vanish with the ratio of the time scales, because each agent retains its own strategy. Particle simulations confirm the analytical predictions in all the regimes considered.

The derivation of the kinetic equation in the quasi-invariant regime and the reductions in the two asymptotic regimes are formal, and the fast-strategy reduction relies on explicit additional assumptions; we do not prove convergence of the full solution to the reduced dynamics, nor a rigorous mean-field limit. The exact analysis requires state-independent payoffs, the three-strategy case is explored through representative scenarios rather than a classification, and for state-dependent payoffs we considered a benchmark quadratic family and characterized candidate stationary states, without proving their existence and stability.

There are several natural extensions to the current work. Proving the mean-field limit and the convergence in the fast regimes, which for two strategies and constant payoffs may be within reach thanks to the triangular structure of the problem, and studying the existence and stability of segregated states for state-dependent payoffs. The mismatch found for the second moments also suggests a reduced description that keeps the individual strategy as a parameter instead of replacing it by its average. Further, one could add noise or mutation to the strategy update, which would regularize the degenerate boundary behavior of the replicator flow and favor interior strategic compositions. Of further interest would be to treat the case of several populations with different payoff structures, for instance a fraction of inherently more stubborn agents who gain more from retaining their convictions. This requires asymmetric games, in which the replicator dynamics of each population is governed by its own payoff matrix.

Replacing the finite set of pure strategies by a continuum of compromise propensities, so that the mixed strategy becomes a probability measure on $[0,1]$ describing how likely an agent is to compromise rather than to remain stubborn is another natural extension of the present work. For state-independent payoffs one expects the finite-dimensional projective flow to be replaced by an exponential tilting of the initial measure, as for replicator dynamics on infinite strategy spaces \cite{oechssler2001evolutionary}. Further, one can consider other strategy dynamics such as imitation or best-response rules in place of the replicator rule, and payoffs that depend on the opinion state beyond the quadratic family considered here. Lastly, allowing a third time scale related to slowly evolving convictions, as well as the inclusion of networks or spatial interaction, and heterogeneous susceptibility to social influence, including credibility biases and population structures in social learning \cite{billiard2020convergence}.

\section*{Acknowledgments}
This work has been written within the activities of GNCS groups of INdAM (Italian National Institute of High Mathematics). GA acknowledges the support of the PRIN 2026 - 2026PYFEZ9, ``Emergent Structures in Kinetic Theory: from Gases to Plasmas" funded with the contribution of MUR. DNR has received funding from the National Science Centre, Poland, grant number 2023/50/A/ST1/00447.

\section*{Contributions}
All authors contributed to conceptualization, formal analysis, investigation, methodology, and writing--review and editing. MP additionally contributed to software development and numerical experiments. All authors reviewed and approved the final manuscript.

\bibliographystyle{abbrv}
\bibliography{biblio}

\end{document}